\documentclass[a4paper,leqno]{article}
\usepackage{amsfonts}
\usepackage{amsmath}
\usepackage{amssymb}
\numberwithin{equation}{section}

\makeatletter
\def\diagram{\leftwidth=\z@ \rightwidth=\z@ \topheight=\z@
\botheight=\z@ \setbox\@picbox\hbox\bgroup}
\def\enddiagram{\egroup\wd\@picbox\rightwidth\unitlength
\ht\@picbox\topheight\unitlength \dp\@picbox\botheight\unitlength
\hskip\leftwidth\unitlength\box\@picbox}
\def\bfig{\begin{diagram}}
\def\efig{\end{diagram}}
\newcount\wideness \newcount\leftwidth \newcount\rightwidth
\newcount\highness \newcount\topheight \newcount\botheight
\def\ratchet#1#2{\ifnum#1<#2 \global #1=#2 \fi}
\def\putbox(#1,#2)#3{%
\horsize{\wideness}{#3} \divide\wideness by 2 {\advance\wideness
by #1 \ratchet{\rightwidth}{\wideness}} {\advance\wideness by -#1
\ratchet{\leftwidth}{\wideness}} \vertsize{\highness}{#3}
\divide\highness by 2 {\advance\highness by #2
\ratchet{\topheight}{\highness}} {\advance\highness by -#2
\ratchet{\botheight}{\highness}} \put(#1,#2){\makebox(0,0){$#3$}}}
\def\putlbox(#1,#2)#3{%
\horsize{\wideness}{#3} {\advance\wideness by #1
\ratchet{\rightwidth}{\wideness}} {\ratchet{\leftwidth}{-#1}}
\vertsize{\highness}{#3} \divide\highness by 2 {\advance\highness
by #2 \ratchet{\topheight}{\highness}} {\advance\highness by -#2
\ratchet{\botheight}{\highness}}
\put(#1,#2){\makebox(0,0)[l]{$#3$}}}
\def\putrbox(#1,#2)#3{%
\horsize{\wideness}{#3} {\ratchet{\rightwidth}{#1}}
{\advance\wideness by -#1 \ratchet{\leftwidth}{\wideness}}
\vertsize{\highness}{#3} \divide\highness by 2 {\advance\highness
by #2 \ratchet{\topheight}{\highness}} {\advance\highness by -#2
\ratchet{\botheight}{\highness}}
\put(#1,#2){\makebox(0,0)[r]{$#3$}}}

\def\adjust[#1]{} 
\newcount \coefa
\newcount \coefb
\newcount \coefc
\newcount\tempcounta
\newcount\tempcountb
\newcount\tempcountc
\newcount\tempcountd
\newcount\xext
\newcount\yext
\newcount\xoff
\newcount\yoff
\newcount\gap%
\newcount\arrowtypea
\newcount\arrowtypeb
\newcount\arrowtypec
\newcount\arrowtyped
\newcount\arrowtypee
\newcount\height
\newcount\width
\newcount\xpos
\newcount\ypos
\newcount\run
\newcount\rise
\newcount\arrowlength
\newcount\halflength
\newcount\arrowtype
\newdimen\tempdimen
\newdimen\xlen
\newdimen\ylen
\newsavebox{\tempboxa}%
\newsavebox{\tempboxb}%
\newsavebox{\tempboxc}%
\newdimen\w@dth
\def\setw@dth#1#2{\setbox\z@\hbox{$#1$}\w@dth=\wd\z@
\setbox\@ne\hbox{$#2$}\ifnum\w@dth<\wd\@ne \w@dth=\wd\@ne \fi
\advance\w@dth by 1.2em}
\def\t@^#1_#2{\def\n@one{#1}\def\n@two{#2}\mathrel{\setw@dth{#1}{#2}
\mathop{\hbox to \w@dth{\rightarrowfill}}\limits
\ifx\n@one\empty\else ^{\box\z@}\fi \ifx\n@two\empty\else
_{\box\@ne}\fi}}
\def\t@@^#1{\@ifnextchar_ {\t@^{#1}}{\t@^{#1}_{}}}
\def\to{\@ifnextchar^ {\t@@}{\t@@^{}}}
\def\t@left^#1_#2{\def\n@one{#1}\def\n@two{#2}\mathrel{\setw@dth{#1}{#2}
\mathop{\hbox to \w@dth{\leftarrowfill}}\limits
\ifx\n@one\empty\else ^{\box\z@}\fi \ifx\n@two\empty\else
_{\box\@ne}\fi}}
\def\t@@left^#1{\@ifnextchar_ {\t@left^{#1}}{\t@left^{#1}_{}}}
\def\toleft{\@ifnextchar^ {\t@@left}{\t@@left^{}}}
\def\two@^#1_#2{\def\n@one{#1}\def\n@two{#2}\mathrel{\setw@dth{#1}{#2}
\mathop{\vcenter{\hbox to \w@dth{\rightarrowfill}\kern-1.7ex
                 \hbox to \w@dth{\rightarrowfill}}%
       }\limits
\ifx\n@one\empty\else ^{\box\z@}\fi \ifx\n@two\empty\else
_{\box\@ne}\fi}}
\def\tw@@^#1{\@ifnextchar_ {\two@^{#1}}{\two@^{#1}_{}}}
\def\two{\@ifnextchar^ {\tw@@}{\tw@@^{}}}
\def\tofr@^#1_#2{\def\n@one{#1}\def\n@two{#2}\mathrel{\setw@dth{#1}{#2}
\mathop{\vcenter{\hbox to \w@dth{\rightarrowfill}\kern-1.7ex
                 \hbox to \w@dth{\leftarrowfill}}%
       }\limits
\ifx\n@one\empty\else ^{\box\z@}\fi \ifx\n@two\empty\else
_{\box\@ne}\fi}}
\def\t@fr@^#1{\@ifnextchar_ {\tofr@^{#1}}{\tofr@^{#1}_{}}}
\def\tofro{\@ifnextchar^ {\t@fr@}{\t@fr@^{}}}

\def\mon{\mathop{\m@th\hbox to
      14.6\P@{\lasyb\char'51\hskip-2.1\P@$\arrext$\hss
$\mathord\rightarrow$}}\limits} 
\def\leftmono{\mathrel{\m@th\hbox to
14.6\P@{$\mathord\leftarrow$\hss$\arrext$\hskip-2.1\P@\lasyb\char'50%
}}\limits} 
\mathchardef\arrext="0200       

\def\settypes(#1,#2,#3){\arrowtypea#1 \arrowtypeb#2 \arrowtypec#3}
\def\settoheight#1#2{\setbox\@tempboxa\hbox{#2}#1\ht\@tempboxa\relax}%
\def\settodepth#1#2{\setbox\@tempboxa\hbox{#2}#1\dp\@tempboxa\relax}%
\def\settokens[#1`#2`#3`#4]{%
     \def\tokena{#1}\def\tokenb{#2}\def\tokenc{#3}\def\tokend{#4}}
\def\setsqparms[#1`#2`#3`#4;#5`#6]{%
\arrowtypea #1 \arrowtypeb #2 \arrowtypec #3 \arrowtyped #4
\width #5 \height #6 }
\def\setpos(#1,#2){\xpos=#1 \ypos#2}

\def\settriparms[#1`#2`#3;#4]{\settripairparms[#1`#2`#3`1`1;#4]}%
\def\settripairparms[#1`#2`#3`#4`#5;#6]{%
\arrowtypea #1 \arrowtypeb #2 \arrowtypec #3 \arrowtyped #4
\arrowtypee #5 \width #6 \height #6 }
\def\resetparms{\settripairparms[1`1`1`1`1;500]\width 500}
\resetparms
\def\mvector(#1,#2)#3{
\put(0,0){\vector(#1,#2){#3}}%
\put(0,0){\vector(#1,#2){26}}%
}
\def\evector(#1,#2)#3{{
\arrowlength #3
\put(0,0){\vector(#1,#2){\arrowlength}}%
\advance \arrowlength by-30
\put(0,0){\vector(#1,#2){\arrowlength}}%
}}
\def\horsize#1#2{%
\settowidth{\tempdimen}{$#2$}%
#1=\tempdimen \divide #1 by\unitlength }
\def\vertsize#1#2{%
\settoheight{\tempdimen}{$#2$}%
#1=\tempdimen
\settodepth{\tempdimen}{$#2$}%
\advance #1 by\tempdimen \divide #1 by\unitlength }
\def\putvector(#1,#2)(#3,#4)#5#6{{%
\ifnum3<\arrowtype \putdashvector(#1,#2)(#3,#4)#5\arrowtype \else
\ifnum\arrowtype<-3 \putdashvector(#1,#2)(#3,#4)#5\arrowtype \else
\xpos=#1 \ypos=#2 \run=#3 \rise=#4 \arrowlength=#5 \ifnum
\arrowtype<0
    \ifnum \run=0
        \advance \ypos by-\arrowlength
    \else
        \tempcounta \arrowlength
        \multiply \tempcounta by\rise
        \divide \tempcounta by\run
        \ifnum\run>0
            \advance \xpos by\arrowlength
            \advance \ypos by\tempcounta
        \else
            \advance \xpos by-\arrowlength
            \advance \ypos by-\tempcounta
        \fi
    \fi
    \multiply \arrowtype by-1
    \multiply \rise by-1
    \multiply \run by-1
\fi \ifcase \arrowtype
\or \put(\xpos,\ypos){\vector(\run,\rise){\arrowlength}}%
\or \put(\xpos,\ypos){\mvector(\run,\rise)\arrowlength}%
\or \put(\xpos,\ypos){\evector(\run,\rise){\arrowlength}}%
\fi\fi\fi }}
\def\putsplitvector(#1,#2)#3#4{
\xpos #1 \ypos #2 \arrowtype #4 \halflength #3 \arrowlength #3
\gap 140 \advance \halflength by-\gap \divide \halflength by2
\ifnum\arrowtype>0
   \ifcase \arrowtype
   \or \put(\xpos,\ypos){\line(0,-1){\halflength}}%
       \advance\ypos by-\halflength
       \advance\ypos by-\gap
       \put(\xpos,\ypos){\vector(0,-1){\halflength}}%
   \or \put(\xpos,\ypos){\line(0,-1)\halflength}%
       \put(\xpos,\ypos){\vector(0,-1)3}%
       \advance\ypos by-\halflength
       \advance\ypos by-\gap
       \put(\xpos,\ypos){\vector(0,-1){\halflength}}%
   \or \put(\xpos,\ypos){\line(0,-1)\halflength}%
       \advance\ypos by-\halflength
       \advance\ypos by-\gap
       \put(\xpos,\ypos){\evector(0,-1){\halflength}}%
   \fi
\else \arrowtype=-\arrowtype
   \ifcase\arrowtype
   \or \advance \ypos by-\arrowlength
       \put(\xpos,\ypos){\line(0,1){\halflength}}%
       \advance\ypos by\halflength
       \advance\ypos by\gap
       \put(\xpos,\ypos){\vector(0,1){\halflength}}%
   \or \advance \ypos by-\arrowlength
       \put(\xpos,\ypos){\line(0,1)\halflength}%
       \put(\xpos,\ypos){\vector(0,1)3}%
       \advance\ypos by\halflength
       \advance\ypos by\gap
       \put(\xpos,\ypos){\vector(0,1){\halflength}}%
   \or \advance \ypos by-\arrowlength
       \put(\xpos,\ypos){\line(0,1)\halflength}%
       \advance\ypos by\halflength
       \advance\ypos by\gap
       \put(\xpos,\ypos){\evector(0,1){\halflength}}%
   \fi
\fi }
\def\putmorphism(#1)(#2,#3)[#4`#5`#6]#7#8#9{{%
\run #2 \rise #3 \ifnum\rise=0
  \puthmorphism(#1)[#4`#5`#6]{#7}{#8}#9%
\else\ifnum\run=0
  \putvmorphism(#1)[#4`#5`#6]{#7}{#8}#9%
\else
\setpos(#1)%
\arrowlength #7 \arrowtype #8 \ifnum\run=0 \else\ifnum\rise=0
\else \ifnum\run>0
    \coefa=1
\else
   \coefa=-1
\fi \ifnum\arrowtype>0
   \coefb=0
   \coefc=-1
\else
   \coefb=\coefa
   \coefc=1
   \arrowtype=-\arrowtype
\fi \width=2 \multiply \width by\run \divide \width by\rise
\ifnum \width<0  \width=-\width\fi \advance\width by60 \if l#9
\width=-\width\fi
\putbox(\xpos,\ypos){#4}
{\multiply \coefa by\arrowlength
\advance\xpos by\coefa \multiply \coefa by\rise \divide \coefa
by\run \advance \ypos by\coefa
\putbox(\xpos,\ypos){#5} }%
{\multiply \coefa by\arrowlength
\divide \coefa by2 \advance \xpos by\coefa \advance \xpos by\width
\multiply \coefa by\rise \divide \coefa by\run \advance \ypos
by\coefa
\if l#9%
   \putrbox(\xpos,\ypos){#6}%
\else\if r#9%
   \putlbox(\xpos,\ypos){#6}%
\fi\fi }%
{\multiply \rise by-\coefc
\multiply \run by-\coefc \multiply \coefb by\arrowlength \advance
\xpos by\coefb \multiply \coefb by\rise \divide \coefb by\run
\advance \ypos by\coefb \multiply \coefc by70 \advance \ypos
by\coefc \multiply \coefc by\run \divide \coefc by\rise \advance
\xpos by\coefc \multiply \coefa by140 \multiply \coefa by\run
\divide \coefa by\rise \advance \arrowlength by\coefa
\ifcase\arrowtype
\or \put(\xpos,\ypos){\vector(\run,\rise){\arrowlength}}%
\or \put(\xpos,\ypos){\mvector(\run,\rise){\arrowlength}}%
\or \put(\xpos,\ypos){\evector(\run,\rise){\arrowlength}}%
\fi}\fi\fi\fi\fi}}

\newcount\numbdashes \newcount\lengthdash \newcount\increment
\def\howmanydashes{
\numbdashes=\arrowlength \lengthdash=40 \divide\numbdashes by
\lengthdash \lengthdash=\arrowlength \divide\lengthdash by
\numbdashes
\increment=\lengthdash \multiply\lengthdash by 3
\divide\lengthdash by 5 }
\def\putdashvector(#1)(#2,#3)#4#5{%
\ifnum#3=0 \putdashhvector(#1){#4}#5 \else \ifnum#2=0
\putdashvvector(#1){#4}#5\fi\fi}
\def\putdashhvector(#1,#2)#3#4{{%
\arrowlength=#3 \howmanydashes
\multiput(#1,#2)(\increment,0){\numbdashes}%
{\vrule height .4pt width \lengthdash\unitlength} \arrowtype=#4
\xpos=#1 \ifnum\arrowtype<0 \advance\arrowtype by 7 \fi
\ifcase\arrowtype \or \advance\xpos by 10
    \put(\xpos,#2){\vector(-1,0){\lengthdash}}
    \advance\xpos by 40
    \put(\xpos,#2){\vector(-1,0){\lengthdash}}
\or \advance \xpos by 10
    \put(\xpos,#2){\vector(-1,0){\lengthdash}}
    \advance\xpos by  \arrowlength
    \advance\xpos by  -50
    \put(\xpos,#2){\vector(-1,0){\lengthdash}}
\or \advance\xpos by 10
    \put(\xpos,#2){\vector(-1,0){\lengthdash}}
\or \advance\xpos by \arrowlength
    \advance\xpos by -\lengthdash
    \put(\xpos,#2){\vector(1,0){\lengthdash}}
\or {\advance\xpos by 10
    \put(\xpos,#2){\vector(1,0){\lengthdash}}}
    \advance\xpos by \arrowlength
    \advance\xpos by -\lengthdash
    \put(\xpos,#2){\vector(1,0){\lengthdash}}
\or \advance\xpos by \arrowlength
    \advance\xpos by -\lengthdash
    \put(\xpos,#2){\vector(1,0){\lengthdash}}
    \advance\xpos by -40
    \put(\xpos,#2){\vector(1,0){\lengthdash}}
   \fi
}}
\def\putdashvvector(#1,#2)#3#4{{%
\arrowlength=#3 \howmanydashes \ypos=#2 \advance\ypos by
-\arrowlength
\multiput(#1,#2)(0,\increment){\numbdashes}%
    {\vrule width .4pt height \lengthdash\unitlength}
\arrowtype=#4 \ypos=#2 \ifnum\arrowtype<0 \advance\arrowtype by 7
\fi \ifcase\arrowtype \or \advance\ypos by \arrowlength
\advance\ypos by -40
    \put(#1,\ypos){\vector(0,1){\lengthdash}}
    \advance\ypos by -40
    \put(#1,\ypos){\vector(0,1){\lengthdash}}
\or \advance\ypos by 10
    \put(#1,\ypos){\vector(0,1){\lengthdash}}
    \advance\ypos by \arrowlength \advance\ypos by -40
    \put(#1,\ypos){\vector(0,1){\lengthdash}}
\or \advance\ypos by \arrowlength \advance\ypos by -40
    \put(#1,\ypos){\vector(0,1){\lengthdash}}
\or \advance\ypos by 10
    \put(#1,\ypos){\vector(0,-1){\lengthdash}}
\or \advance\ypos by 10
    \put(#1,\ypos){\vector(0,-1){\lengthdash}}
    \advance\ypos by \arrowlength \advance\ypos by -40
    \put(#1,\ypos){\vector(0,-1){\lengthdash}}
\or \advance\ypos by 10
    \put(#1,\ypos){\vector(0,-1){\lengthdash}}
    \advance\ypos by 40
    \put(#1,\ypos){\vector(0,-1){\lengthdash}}
\fi }}
\def\puthmorphism(#1,#2)[#3`#4`#5]#6#7#8{{%
\xpos #1 \ypos #2 \width #6 \arrowlength #6 \arrowtype=#7
\putbox(\xpos,\ypos){#3\vphantom{#4}}%
{\advance \xpos by\arrowlength
\putbox(\xpos,\ypos){\vphantom{#3}#4}}%
\horsize{\tempcounta}{#3}%
\horsize{\tempcountb}{#4}%
\divide \tempcounta by2 \divide \tempcountb by2 \advance
\tempcounta by30 \advance \tempcountb by30 \advance \xpos
by\tempcounta \advance \arrowlength by-\tempcounta \advance
\arrowlength by-\tempcountb
\putvector(\xpos,\ypos)(1,0)\arrowlength\arrowtype \divide
\arrowlength by2 \advance \xpos by\arrowlength
\vertsize{\tempcounta}{#5}%
\divide\tempcounta by2 \advance \tempcounta by20
\if a#8 %
   \advance \ypos by\tempcounta
   \putbox(\xpos,\ypos){#5}%
\else
   \advance \ypos by-\tempcounta
   \putbox(\xpos,\ypos){#5}%
\fi}}
\def\putvmorphism(#1,#2)[#3`#4`#5]#6#7#8{{%
\xpos #1 \ypos #2 \arrowlength #6 \arrowtype #7
\settowidth{\xlen}{$#5$}%
\putbox(\xpos,\ypos){#3}%
{\advance \ypos by-\arrowlength
\putbox(\xpos,\ypos){#4}}%
{\advance\arrowlength by-140 \advance \ypos by-70 \ifdim\xlen>0pt
   \if m#8%
      \putsplitvector(\xpos,\ypos)\arrowlength\arrowtype
   \else
   \putvector(\xpos,\ypos)(0,-1)\arrowlength\arrowtype
   \fi
\else
   \putvector(\xpos,\ypos)(0,-1)\arrowlength\arrowtype
\fi}%
\ifdim\xlen>0pt
   \divide \arrowlength by2
   \advance\ypos by-\arrowlength
   \if l#8%
      \advance \xpos by-40
      \putrbox(\xpos,\ypos){#5}%
   \else\if r#8%
      \advance \xpos by40
      \putlbox(\xpos,\ypos){#5}%
   \else
      \putbox(\xpos,\ypos){#5}%
   \fi\fi
\fi }}
\def\putsquarep<#1>(#2)[#3;#4`#5`#6`#7]{{%
\setsqparms[#1]%
\setpos(#2)%
\settokens[#3]%
\puthmorphism(\xpos,\ypos)[\tokenc`\tokend`{#7}]{\width}{\arrowtyped}b%
\advance\ypos by \height
\puthmorphism(\xpos,\ypos)[\tokena`\tokenb`{#4}]{\width}{\arrowtypea}a%
\putvmorphism(\xpos,\ypos)[``{#5}]{\height}{\arrowtypeb}l%
\advance\xpos by \width
\putvmorphism(\xpos,\ypos)[``{#6}]{\height}{\arrowtypec}r%
}}
\def\putsquare{\@ifnextchar <{\putsquarep}{\putsquarep%
   <\arrowtypea`\arrowtypeb`\arrowtypec`\arrowtyped;\width`\height>}}
\def\square{\@ifnextchar< {\squarep}{\squarep
   <\arrowtypea`\arrowtypeb`\arrowtypec`\arrowtyped;\width`\height>}}
\def\squarep<#1>[#2`#3`#4`#5;#6`#7`#8`#9]{{
\setsqparms[#1]
\diagram
\putsquarep<\arrowtypea`\arrowtypeb`\arrowtypec`
\arrowtyped;\width`\height>
(0,0)[#2`#3`#4`{#5};#6`#7`#8`{#9}]
\enddiagram
}}                                                 
\def\putptrianglep<#1>(#2,#3)[#4`#5`#6;#7`#8`#9]{{%
\settriparms[#1]%
\xpos=#2 \ypos=#3 \advance\ypos by \height
\puthmorphism(\xpos,\ypos)[#4`#5`{#7}]{\height}{\arrowtypea}a%
\putvmorphism(\xpos,\ypos)[`#6`{#8}]{\height}{\arrowtypeb}l%
\advance\xpos by\height
\putmorphism(\xpos,\ypos)(-1,-1)[``{#9}]{\height}{\arrowtypec}r%
}}
\def\putptriangle{\@ifnextchar <{\putptrianglep}{\putptrianglep
   <\arrowtypea`\arrowtypeb`\arrowtypec;\height>}}
\def\ptriangle{\@ifnextchar <{\ptrianglep}{\ptrianglep
   <\arrowtypea`\arrowtypeb`\arrowtypec;\height>}}
\def\ptrianglep<#1>[#2`#3`#4;#5`#6`#7]{{
\settriparms[#1]
\diagram
\putptrianglep<\arrowtypea`\arrowtypeb`
\arrowtypec;\height>
(0,0)[#2`#3`#4;#5`#6`{#7}]
\enddiagram
}}                                            
\def\putqtrianglep<#1>(#2,#3)[#4`#5`#6;#7`#8`#9]{{%
\settriparms[#1]%
\xpos=#2 \ypos=#3 \advance\ypos by\height
\puthmorphism(\xpos,\ypos)[#4`#5`{#7}]{\height}{\arrowtypea}a%
\putmorphism(\xpos,\ypos)(1,-1)[``{#8}]{\height}{\arrowtypeb}l%
\advance\xpos by\height
\putvmorphism(\xpos,\ypos)[`#6`{#9}]{\height}{\arrowtypec}r%
}}
\def\putqtriangle{\@ifnextchar <{\putqtrianglep}{\putqtrianglep
   <\arrowtypea`\arrowtypeb`\arrowtypec;\height>}}
\def\qtriangle{\@ifnextchar <{\qtrianglep}{\qtrianglep
   <\arrowtypea`\arrowtypeb`\arrowtypec;\height>}}
\def\qtrianglep<#1>[#2`#3`#4;#5`#6`#7]{{
\settriparms[#1]
\width=\height                                
\diagram
\putqtrianglep<\arrowtypea`\arrowtypeb`
\arrowtypec;\height>
(0,0)[#2`#3`#4;#5`#6`{#7}]
\enddiagram
}}
\def\putdtrianglep<#1>(#2,#3)[#4`#5`#6;#7`#8`#9]{{%
\settriparms[#1]%
\xpos=#2 \ypos=#3
\puthmorphism(\xpos,\ypos)[#5`#6`{#9}]{\height}{\arrowtypec}b%
\advance\xpos by \height \advance\ypos by\height
\putmorphism(\xpos,\ypos)(-1,-1)[``{#7}]{\height}{\arrowtypea}l%
\putvmorphism(\xpos,\ypos)[#4``{#8}]{\height}{\arrowtypeb}r%
}}
\def\putdtriangle{\@ifnextchar <{\putdtrianglep}{\putdtrianglep
   <\arrowtypea`\arrowtypeb`\arrowtypec;\height>}}
\def\dtriangle{\@ifnextchar <{\dtrianglep}{\dtrianglep
   <\arrowtypea`\arrowtypeb`\arrowtypec;\height>}}
\def\dtrianglep<#1>[#2`#3`#4;#5`#6`#7]{{
\settriparms[#1]
\width=\height                                
\diagram
\putdtrianglep<\arrowtypea`\arrowtypeb`
\arrowtypec;\height>
(0,0)[#2`#3`#4;#5`#6`{#7}]
\enddiagram
}}
\def\putbtrianglep<#1>(#2,#3)[#4`#5`#6;#7`#8`#9]{{%
\settriparms[#1]%
\xpos=#2 \ypos=#3
\puthmorphism(\xpos,\ypos)[#5`#6`{#9}]{\height}{\arrowtypec}b%
\advance\ypos by\height
\putmorphism(\xpos,\ypos)(1,-1)[``{#8}]{\height}{\arrowtypeb}r%
\putvmorphism(\xpos,\ypos)[#4``{#7}]{\height}{\arrowtypea}l%
}}
\def\putbtriangle{\@ifnextchar <{\putbtrianglep}{\putbtrianglep
   <\arrowtypea`\arrowtypeb`\arrowtypec;\height>}}
\def\btriangle{\@ifnextchar <{\btrianglep}{\btrianglep
   <\arrowtypea`\arrowtypeb`\arrowtypec;\height>}}
\def\btrianglep<#1>[#2`#3`#4;#5`#6`#7]{{
\settriparms[#1]
\width=\height                               
\diagram
\putbtrianglep<\arrowtypea`\arrowtypeb`
\arrowtypec;\height>
(0,0)[#2`#3`#4;#5`#6`{#7}]
\enddiagram
}}
\def\putAtrianglep<#1>(#2,#3)[#4`#5`#6;#7`#8`#9]{{%
\settriparms[#1]%
\xpos=#2 \ypos=#3 {\multiply \height by2
\puthmorphism(\xpos,\ypos)[#5`#6`{#9}]{\height}{\arrowtypec}b}%
\advance\xpos by\height \advance\ypos by\height
\putmorphism(\xpos,\ypos)(-1,-1)[#4``{#7}]{\height}{\arrowtypea}l%
\putmorphism(\xpos,\ypos)(1,-1)[``{#8}]{\height}{\arrowtypeb}r%
}}
\def\putAtriangle{\@ifnextchar <{\putAtrianglep}{\putAtrianglep
   <\arrowtypea`\arrowtypeb`\arrowtypec;\height>}}
\def\Atriangle{\@ifnextchar <{\Atrianglep}{\Atrianglep
   <\arrowtypea`\arrowtypeb`\arrowtypec;\height>}}
\def\Atrianglep<#1>[#2`#3`#4;#5`#6`#7]{{
\settriparms[#1]
\width=\height                                     
\diagram
\putAtrianglep<\arrowtypea`\arrowtypeb`
\arrowtypec;\height>
(0,0)[#2`#3`#4;#5`#6`{#7}]
\enddiagram
}}
\def\putAtrianglepairp<#1>(#2)[#3;#4`#5`#6`#7`#8]{{%
\settripairparms[#1]%
\setpos(#2)%
\settokens[#3]%
\puthmorphism(\xpos,\ypos)[\tokenb`\tokenc`{#7}]{\height}{\arrowtyped}b%
\advance\xpos by\height
\puthmorphism(\xpos,\ypos)[\phantom{\tokenc}`\tokend`{#8}]%
{\height}{\arrowtypee}b%
\advance\ypos by\height
\putmorphism(\xpos,\ypos)(-1,-1)[\tokena``{#4}]{\height}{\arrowtypea}l%
\putvmorphism(\xpos,\ypos)[``{#5}]{\height}{\arrowtypeb}m%
\putmorphism(\xpos,\ypos)(1,-1)[``{#6}]{\height}{\arrowtypec}r%
}}
\def\putAtrianglepair{\@ifnextchar <{\putAtrianglepairp}{\putAtrianglepairp%
   <\arrowtypea`\arrowtypeb`\arrowtypec`\arrowtyped`\arrowtypee;\height>}}
\def\Atrianglepair{\@ifnextchar <{\Atrianglepairp}{\Atrianglepairp%
   <\arrowtypea`\arrowtypeb`\arrowtypec`\arrowtyped`\arrowtypee;\height>}}
\def\Atrianglepairp<#1>[#2;#3`#4`#5`#6`#7]{{
\settripairparms[#1]
\settokens[#2]
\width=\height                                
\diagram
\putAtrianglepairp                            
<\arrowtypea`\arrowtypeb`\arrowtypec`
\arrowtyped`\arrowtypee;\height>
(0,0)[{#2};#3`#4`#5`#6`{#7}]
\enddiagram
}}
\def\putVtrianglep<#1>(#2,#3)[#4`#5`#6;#7`#8`#9]{{%
\settriparms[#1]%
\xpos=#2 \ypos=#3 \advance\ypos by\height {\multiply\height by2
\puthmorphism(\xpos,\ypos)[#4`#5`{#7}]{\height}{\arrowtypea}a}%
\putmorphism(\xpos,\ypos)(1,-1)[`#6`{#8}]{\height}{\arrowtypeb}l%
\advance\xpos by\height \advance\xpos by\height
\putmorphism(\xpos,\ypos)(-1,-1)[``{#9}]{\height}{\arrowtypec}r%
}}
\def\putVtriangle{\@ifnextchar <{\putVtrianglep}{\putVtrianglep
   <\arrowtypea`\arrowtypeb`\arrowtypec;\height>}}
\def\Vtriangle{\@ifnextchar <{\Vtrianglep}{\Vtrianglep
   <\arrowtypea`\arrowtypeb`\arrowtypec;\height>}}
\def\Vtrianglep<#1>[#2`#3`#4;#5`#6`#7]{{
\settriparms[#1]
\width=\height                                 
\diagram
\putVtrianglep<\arrowtypea`\arrowtypeb`
\arrowtypec;\height>
(0,0)[#2`#3`#4;#5`#6`{#7}]
\enddiagram
}}
\def\putVtrianglepairp<#1>(#2)[#3;#4`#5`#6`#7`#8]{{
\settripairparms[#1]%
\setpos(#2)%
\settokens[#3]%
\advance\ypos by\height
\putmorphism(\xpos,\ypos)(1,-1)[`\tokend`{#6}]{\height}{\arrowtypec}l%
\puthmorphism(\xpos,\ypos)[\tokena`\tokenb`{#4}]{\height}{\arrowtypea}a%
\advance\xpos by\height
\puthmorphism(\xpos,\ypos)[\phantom{\tokenb}`\tokenc`{#5}]%
{\height}{\arrowtypeb}a%
\putvmorphism(\xpos,\ypos)[``{#7}]{\height}{\arrowtyped}m%
\advance\xpos by\height
\putmorphism(\xpos,\ypos)(-1,-1)[``{#8}]{\height}{\arrowtypee}r%
}}
\def\putVtrianglepair{\@ifnextchar <{\putVtrianglepairp}{\putVtrianglepairp%
    <\arrowtypea`\arrowtypeb`\arrowtypec`\arrowtyped`\arrowtypee;\height>}}
\def\Vtrianglepair{\@ifnextchar <{\Vtrianglepairp}{\Vtrianglepairp%
    <\arrowtypea`\arrowtypeb`\arrowtypec`\arrowtyped`\arrowtypee;\height>}}
\def\Vtrianglepairp<#1>[#2;#3`#4`#5`#6`#7]{{
\settripairparms[#1]
\settokens[#2]
\diagram
\putVtrianglepairp                             
<\arrowtypea`\arrowtypeb`\arrowtypec`
\arrowtyped`\arrowtypee;\height>
(0,0)[{#2};#3`#4`#5`#6`{#7}]
\enddiagram
}}

\def\putCtrianglep<#1>(#2,#3)[#4`#5`#6;#7`#8`#9]{{%
\settriparms[#1]%
\xpos=#2 \ypos=#3 \advance\ypos by\height
\putmorphism(\xpos,\ypos)(1,-1)[``{#9}]{\height}{\arrowtypec}l%
\advance\xpos by\height \advance\ypos by\height
\putmorphism(\xpos,\ypos)(-1,-1)[#4`#5`{#7}]{\height}{\arrowtypea}l%
{\multiply\height by 2
\putvmorphism(\xpos,\ypos)[`#6`{#8}]{\height}{\arrowtypeb}r}%
}}
\def\putCtriangle{\@ifnextchar <{\putCtrianglep}{\putCtrianglep
    <\arrowtypea`\arrowtypeb`\arrowtypec;\height>}}
\def\Ctriangle{\@ifnextchar <{\Ctrianglep}{\Ctrianglep
    <\arrowtypea`\arrowtypeb`\arrowtypec;\height>}}
\def\Ctrianglep<#1>[#2`#3`#4;#5`#6`#7]{{
\settriparms[#1]
\width=\height                               
\diagram
\putCtrianglep<\arrowtypea`\arrowtypeb`
\arrowtypec;\height>
(0,0)[#2`#3`#4;#5`#6`{#7}]
\enddiagram
}}                                           
\def\putDtrianglep<#1>(#2,#3)[#4`#5`#6;#7`#8`#9]{{%
\settriparms[#1]%
\xpos=#2 \ypos=#3 \advance\xpos by\height \advance\ypos by\height
\putmorphism(\xpos,\ypos)(-1,-1)[``{#9}]{\height}{\arrowtypec}r%
\advance\xpos by-\height \advance\ypos by\height
\putmorphism(\xpos,\ypos)(1,-1)[`#5`{#8}]{\height}{\arrowtypeb}r%
{\multiply\height by 2
\putvmorphism(\xpos,\ypos)[#4`#6`{#7}]{\height}{\arrowtypea}l}%
}}
\def\putDtriangle{\@ifnextchar <{\putDtrianglep}{\putDtrianglep
    <\arrowtypea`\arrowtypeb`\arrowtypec;\height>}}
\def\Dtriangle{\@ifnextchar <{\Dtrianglep}{\Dtrianglep
   <\arrowtypea`\arrowtypeb`\arrowtypec;\height>}}
\def\Dtrianglep<#1>[#2`#3`#4;#5`#6`#7]{{
\settriparms[#1]
\width=\height                              
\diagram
\putDtrianglep<\arrowtypea`\arrowtypeb`
\arrowtypec;\height>
(0,0)[#2`#3`#4;#5`#6`{#7}]
\enddiagram
}}                                          
\def\setrecparms[#1`#2]{\width=#1 \height=#2}%
\def\recursep<#1`#2>[#3;#4`#5`#6`#7`#8]{{%
\width=#1 \height=#2 \settokens[#3]
\settowidth{\tempdimen}{$\tokena$} \ifdim\tempdimen=0pt
  \savebox{\tempboxa}{\hbox{$\tokenb$}}%
  \savebox{\tempboxb}{\hbox{$\tokend$}}%
  \savebox{\tempboxc}{\hbox{$#6$}}%
\else
  \savebox{\tempboxa}{\hbox{$\hbox{$\tokena$}\times\hbox{$\tokenb$}$}}%
  \savebox{\tempboxb}{\hbox{$\hbox{$\tokena$}\times\hbox{$\tokend$}$}}%
  \savebox{\tempboxc}{\hbox{$\hbox{$\tokena$}\times\hbox{$#6$}$}}%
\fi \ypos=\height \divide\ypos by 2 \xpos=\ypos \advance\xpos by
\width \bfig
\putCtrianglep<-1`1`1;\ypos>(0,0)[`\tokenc`;#5`#6`{#7}]%
\puthmorphism(\ypos,0)[\tokend`\usebox{\tempboxb}`{#8}]{\width}{-1}b%
\puthmorphism(\ypos,\height)[\tokenb`\usebox{\tempboxa}`{#4}]{\width}{-1}a%
\advance\ypos by \width
\putvmorphism(\ypos,\height)[``\usebox{\tempboxc}]{\height}1r%
\efig }}
\def\recurse{\@ifnextchar <{\recursep}{\recursep<\width`\height>}}
\def\puttwohmorphisms(#1,#2)[#3`#4;#5`#6]#7#8#9{{%
\puthmorphism(#1,#2)[#3`#4`]{#7}0a \ypos=#2 \advance\ypos by 20
\puthmorphism(#1,\ypos)[\phantom{#3}`\phantom{#4}`#5]{#7}{#8}a
\advance\ypos by -40
\puthmorphism(#1,\ypos)[\phantom{#3}`\phantom{#4}`#6]{#7}{#9}b }}
\def\puttwovmorphisms(#1,#2)[#3`#4;#5`#6]#7#8#9{{%
\putvmorphism(#1,#2)[#3`#4`]{#7}0a \xpos=#1 \advance\xpos by -20
\putvmorphism(\xpos,#2)[\phantom{#3}`\phantom{#4}`#5]{#7}{#8}l
\advance\xpos by 40
\putvmorphism(\xpos,#2)[\phantom{#3}`\phantom{#4}`#6]{#7}{#9}r }}
\def\puthcoequalizer(#1)[#2`#3`#4;#5`#6`#7]#8#9{{%
\setpos(#1)%
\puttwohmorphisms(\xpos,\ypos)[#2`#3;#5`#6]{#8}11%
\advance\xpos by #8
\puthmorphism(\xpos,\ypos)[\phantom{#3}`#4`#7]{#8}1{#9} }}
\def\putvcoequalizer(#1)[#2`#3`#4;#5`#6`#7]#8#9{{%
\setpos(#1)%
\puttwovmorphisms(\xpos,\ypos)[#2`#3;#5`#6]{#8}11%
\advance\ypos by -#8
\putvmorphism(\xpos,\ypos)[\phantom{#3}`#4`#7]{#8}1{#9} }}
\def\putthreehmorphisms(#1)[#2`#3;#4`#5`#6]#7(#8)#9{{%
\setpos(#1) \settypes(#8)
\if a#9 %
     \vertsize{\tempcounta}{#5}%
     \vertsize{\tempcountb}{#6}%
     \ifnum \tempcounta<\tempcountb \tempcounta=\tempcountb \fi
\else
     \vertsize{\tempcounta}{#4}%
     \vertsize{\tempcountb}{#5}%
     \ifnum \tempcounta<\tempcountb \tempcounta=\tempcountb \fi
\fi \advance \tempcounta by 60
\puthmorphism(\xpos,\ypos)[#2`#3`#5]{#7}{\arrowtypeb}{#9}
\advance\ypos by \tempcounta
\puthmorphism(\xpos,\ypos)[\phantom{#2}`\phantom{#3}`#4]{#7}{\arrowtypea}{#9}
\advance\ypos by -\tempcounta \advance\ypos by -\tempcounta
\puthmorphism(\xpos,\ypos)[\phantom{#2}`\phantom{#3}`#6]{#7}{\arrowtypec}{#9}
}}
\def\setarrowtoks[#1`#2`#3`#4`#5`#6]{%
\def\toka{#1}
\def\tokb{#2}
\def\tokc{#3}
\def\tokd{#4}
\def\toke{#5}
\def\tokf{#6}
}
\def\hex{\@ifnextchar <{\hexp}{\hexp<1000`400>}}
\def\hexp<#1`#2>[#3`#4`#5`#6`#7`#8;#9]{%
\setarrowtoks[#9] \yext=#2 \advance \yext by #2 \xext=#1
\advance\xext by \yext \bfig
\putCtriangle<-1`0`1;#2>(0,0)[`#5`;\tokb``\tokd] \xext=#1
\yext=#2 \advance \yext by #2
\putsquare<1`0`0`1;\xext`\yext>(#2,0)[#3`#4`#7`#8;\toka```\tokf]
\advance \xext by #2
\putDtriangle<0`1`-1;#2>(\xext,0)[`#6`;`\tokc`\toke] \efig }
\begin{document}
\newtheorem{theorem}{Theorem}[section]
\newtheorem{lemma}[theorem]{Lemma}
\newtheorem{corollary}[theorem]{Corollary}
\newtheorem{conjecture}[theorem]{Conjecture}
\newtheorem{remark}[theorem]{Remark}
\newtheorem{definition}[theorem]{Definition}
\newtheorem{problem}[theorem]{Problem}
\newtheorem{example}[theorem]{Example}
\newtheorem{proposition}[theorem]{Proposition}
\title{{\bf Canonical measures and  dynamical systems of Bergman kernels}}
\date{December 19, 2010}
\author{Hajime TSUJI\footnote{Partially supported by Grant-in-Aid for Scientific Reserch (S) 17104001}}
\maketitle
\begin{abstract}
\noindent In this article, we construct the canonical semipositive  current 
or the canonical measure ($=$ the potential of the canonical semipositive current) on a smooth projective variety with  nonnegative Kodaira dimension 
in terms of a dynamical system of Bergman kernels.  This current is considered to be a generalization of a K\"{a}hler-Einstein metric and coincides the one 
considered  independently by J. Song and G. Tian (\cite{s-t}). 
The major difference between \cite{s-t} and the present article is that 
they found the canonical measure in terms of  K\"{a}her-Ricci flows, 
while I found the canonical measure in terms of dynamical systems 
of Bergman kernels.  Hence the present approach can be viewed as  a 
discrete version of a K\"{a}hler-Ricci flow.   

The advantage of the dynamical construction  is two folds.  First, it enables us to deduce the (logarithmic) plurisubharmonic variation propery of the canonical measures on a projective family.  Second, we can overcome the difficulty arising from the 
singularities of the solution of a K\"{a}hler-Ricci flow.  

\end{abstract}
\tableofcontents
\section{Introduction}

In \cite{KE}, I have constructed a canonical K\"{a}hler-Einstein current  on  a smooth projective variety of general type in terms of a dynamical system of Bergman kernels originated  in \cite{tu3}.   This K\"{a}hler-Einstein current is the same one which has been studied in \cite{tu0,s}. 

Since the same kind of dynamical systems has been defined on a smooth 
projective varieties of  nonnegative Kodaira dimension (or even for 
a smooth projective varietiy with pseudoeffective canonical bundle) in \cite{tu3}, 
it is natural to expect that the normalized limit of the dynamical system of  Bergman kernels yields a substitute of a K\"{a}hler-Einstein metric for a smooth projective variety (of non general type) with nonnegative Kodaira dimension.

The purpose of this article is to  prove that the limit satisfies the partial differential equation (see (\ref{LKE})) similar to the K\"{a}hler-Einstein equation on the base 
of the Iitaka fibration (not on the original variety)  and  give  a natural generalization of the notion of K\"{a}hler-Einstein volume form.
We call the (normalized) limit the canonical measure. 
And we call the $-\mbox{Ric}$ of the canonical measure (in the sense of current) the canonical semipositive current. 

There are two major differences between  K\"{a}hler-Einstein metrics 
and the canonical semipositive currents.  

First of all in general the canonical semipositive current is strictly positive not on the original variety but on the base space of the Iitaka fibration. 
In other words, the current is the pullback of  a closed generically strictly positive current on the base space of the Iitaka fibration. 

Secondary although the canonical semipositive current satisfies a similar partial differential equation as a K\"{a}hler-Einstein metric on the base space of the Iitaka
fibration, the equation has an additional term coming from variation of 
Hodge structure on the Iitaka fibration.   
 
The objective of this generalization is to study the deformation of projective 
varieties with nonnegative Kodaira dimension. 
Actually the dynamical construction of the canonical semipositive current 
yields the existence of a closed semipositive current on the family 
which restricts the canonical semipositive current on the general fibers 
(Theorem \ref{relative}).  We discuss the applications 
of Theorem \ref{relative} in \cite{global}.   And we also note that there are similar constructions of 
canonical measures (canonical and supercanonical AZD's) on smooth projective varieties with pseudoeffective canonical bundles (\cite{canAZD}).  

After the completion of this work,  I have noticed a paper of 
Song and Tian (\cite{s-t}) which  also constructed the canonical semipositive 
current from a  different point of view. 
 Actually first they have constructed the canonical semipositive current as the limit of the K\"{a}hler-Ricci flow in the case of semiample canonical bundle.
Then they constructed the current which satisfies the same equation 
without using a K\"{a}hler-Ricci flow on a projective varieites with nonnegative Kodaira dimension
whose canonical bundle is not necessarily semiample. 
In this sense their construcion is modeled after the case of 
semiample canonical bundles. 
But in their construction, the meaning of the canonical semipositive 
current is not clear (although it is apprarently a generalization of 
a K\"{a}hler-Einstein metric). 
 
The main contribution of this article is to give a dynamical construction of the canonical semipositive currents (or the canonical measure in \cite{s-t}) and give the authenticity to the canonical semipositive current.
   
The advantage of the dynamical construction is that we can overcome the difficulty 
arising from the  singularity of the currents.    
For example it seems to be difficult to deduce the plurisubharmonic variation 
property of the canonical measures on a projective family  (Theorem \ref{relative}) by direct calculation.  On the other hand Theorem \ref{relative} is an immediate consequence 
of the dynamical construction by using the logarithmic plurisubharmonic variation properties of Bergman kernels. 
 
Also difficulty  arises to study a K\"{a}hler-Ricci flow, when we consider non minimal 
algebraic varieties.   In this case  the flow of the K\"{a}hler class 
associated with a K\"{a}hler-Ricci flow  
reaches the boundary of the  K\"{a}hler cone in finite time.  Hence  in this 
case it is inevitable  to deal with a singular K\"{a}hler Ricci flow. 
But the dynamical constrution (Theorem \ref{DS}) automatically produces the 
cannocical semipositve current (or the canonical measure) as soon as 
the Kodaira dimension of the variety is nonnegative. 
One may consider Theorem \ref{DS} as  a discretization of a K\"{a}hler-Ricci flow\footnote{In 
\cite{LC}, I have used another discretization of a K\"{a}hler-Ricci flow.} and it overcomes the difficulty arising from  singularities.

After the completion of this paper (math.ArXiv.0805.1829), I extended the 
result to the case of KLT pairs (\cite{LC}) and  thought that the technique used to solve 
the twisted K\"{a}hler-Einstein equation (\ref{LKE}) is standard.  
But recently it was pointed out by S. Boucksom that the existence of the solution of (\ref{LKE}) does not 
follow from the ready made techniques because the corresponding Monge-Amp`{e}re equation has nonalgebraic singularities.  Moreover  I could not follow the proof of 
\cite[Theorem B.2]{s-t} by J. Song and G. Tian which asserts Theorem \ref{main} in this paper independently.  Hence I decieded to publish this paper independently.  
\vspace{3mm} \\
\noindent{\bf Acknowledgement} 
I would like to express to S. Boucksom who pointed out a crucial error in the previous version.  Namely he pointed out the case $D_{X/Y}\neq 0$ (cf. Section 3.2) is lacking in the previous version.   \vspace{3mm}\\
\noindent{\bf Notations}
\begin{itemize}
\item  For a real number $a$, $\lceil a\rceil$ denotes the minimal integer greater than or equal to $a$ and $\lfloor a\rfloor$ denotes the maximal integer smaller than 
or equal to $a$.  
\item Let $X$ be a projective variety and let $D$ be a Weil divsor on $X$.
Let $D = \sum d_{i}D_{i}$ be the irreducible decomposition.  
We set 
\begin{equation}
\lceil D\rceil := \sum \lceil d_{i}\rceil D_{i}
,\lfloor D \rfloor := \sum \lfloor d_{i}\rfloor D_{i}.
\end{equation}
\item Let $L$ be a line bundle on a compact complex manifold $X$. 
A  singular hermitian metric $h$ on $L$ is given by
\[
h = e^{-\varphi}\cdot h_{0},
\]
where $h_{0}$ is a $C^{\infty}$-hermitian metric on $L$ and 
$\varphi\in L^{1}_{loc}(X)$ is an arbitrary function on $X$.
We call $\varphi$ a  weight function of $h$. 

The curvature current $\Theta_{h}$ of the singular hermitian line
bundle $(L,h)$ is defined by
\[
\Theta_{h} := \Theta_{h_{0}} + \partial\bar{\partial}\varphi,
\]
where $\partial\bar{\partial}\varphi$ is taken in the sense of current. 
The $L^{2}$ sheaf ${\cal L}^{2}(L,h)$ of the singular hermitian
line bundle $(L,h)$ is defined by
\[
{\cal L}^{2}(L,h)(U) := \{ \sigma\in\Gamma (U,{\cal O}_{X}(L))\mid 
\, h(\sigma ,\sigma )\in L^{1}_{loc}(U)\} ,
\]
where $U$ runs over the  open subsets of $X$.
In this case there exists an ideal sheaf ${\cal I}(h)$ such that
\[
{\cal L}^{2}(L,h) = {\cal O}_{X}(L)\otimes {\cal I}(h)
\]
holds.  We call ${\cal I}(h)$ the {\bf multiplier ideal sheaf} of $(L,h)$.
\item For a closed positive $(1,1)$ current $T$, $T_{abc}$ denotes 
the abosolutely continuous part of $T$. 
\item For a singular hermitian line bundle $(F,h_{F})$ on a compact complex 
manifold $X$ of dimension $n$.   $K(X,K_{X}+F,h_{F})$ denotes (the diagonal part of) the Bergman kernel of $H^{0}(X,{\cal O}_{X}(K_{X} + F)\otimes {\cal I}(h_{F}))$ 
with respect to the $L^{2}$-inner product: 
\begin{equation}\label{inner}
(\sigma ,\sigma^{\prime}) := (\sqrt{-1})^{n^{2}}\int_{X}h_{F}\cdot\sigma\wedge \bar{\sigma}^{\prime}, 
\end{equation}  
i.e., 
\begin{equation}\label{BergmanK}
K(X,K_{X}+F,h_{F}) = \sum_{i=0}^{N}|\sigma_{i}|^{2}, 
\end{equation}
where $\{\sigma_{0},\cdots ,\sigma_{N}\}$ is a complete orthonormal basis 
of $H^{0}(X,{\cal O}_{X}(K_{X} + F)\otimes {\cal I}(h_{F}))$. 
It is clear that $K(X,K_{X}+F,h_{F})$ is independent of the choice of 
the complete orthonormal basis. 
\item A line bundle $L$ on a compact complex manifold $X$ is said to 
be {\bf pseudoeffective}, if it admits a singular hermitian metric with 
semipositive curvature current.  A singular hermitian line bundle 
$(L,h)$ is said to be pseudoeffective if the curvarure 
current  $\sqrt{-1}\,\Theta_{h}$ is 
semipositive.  
If $X$ is a smooth projective variety,
this is equivalent to the fact that $c_{1}(L)$ is on the closure 
of the effective cone.  
\item Let $(X,D)$ be a pair of a normal variety and a $\mathbb{Q}$-divisor 
on $X$.  Suppose that $K_{X} + D$ is $\mathbb{Q}$-Cartier. 
Let $f : Y \to X$ be a log resolution.  Then we have the formula :
\[
K_{Y} = f^{*}(K_{X}+D) + \sum a_{i}E_{i}, 
\]
where $E_{i}$ is a prime divisor and $a_{i}\in \mathbb{Q}$. 
The pair $(X,D)$ is said to be {\bf subKLT}(resp. {\bf subLC}, if $a_{i} > -1$
(resp. $a_{i} \geqq -1$ holds for every $i$. 
$(X,D)$ is said to be {\bf KLT} (resp. LC), if $(X,D)$ is subKLT(resp. subLC) and $D$ is effective.
\end{itemize}

\subsection{K\"{a}hler-Einstein metrics}
Let $X$ be a compact K\"{a}hler manifold with the K\"{a}hler form 
\begin{equation}
\omega := \frac{\sqrt{-1}}{2}\sum_{i,j}g_{ij}dz^{i}\wedge d\bar{z}^{j}. 
\end{equation}
$(X,\omega)$ is said to be K\"{a}hler-Einstein, if there exists a constant $c$
such that 
\begin{equation}
\mbox{Ric}_{\omega} = c\omega
\end{equation}
holds, where $\mbox{Ric}_{\omega}$ denotes the Ricci form:
\begin{equation}
\mbox{Ric}_{\omega} := -\sqrt{-1}\partial\bar{\partial}\log \det (g_{ij})
\end{equation}
 and we call $\omega$ a K\"{a}hler-Einsein form on $X$. 
If a compact complex manifold $X$ admits a K\"{a}hler-Einstein form, then 
$c_{1}(X)$ is negative or 0 or positive. 
Conversely by the celebrated solution of Calabi's conjecture (\cite{a,y}),
for a compact K\"{a}hler manifold $X$, 
if $c_{1}(X)$ is negative, there exists a unique K\"{a}hler-Einstein form $\omega_{E}$ such that 
\begin{equation}
-\mbox{Ric}_{\omega_{E}} = \omega_{E} 
\end{equation}
(\cite{a,y}) and if $c_{1}(X)$ is $0$ in $H^{2}(X,\mathbb{R})$, then in  every 
K\"{a}hler class on $X$, there exists a unique K\"{a}hler-Einstein form $\omega_{E}$ 
such that  
\begin{equation}
\mbox{Ric}_{\omega_{E}} = 0
\end{equation} 
holds (\cite{y}). 

\subsection{K\"{a}hler-Einstein currents}
There are numerous applications of K\"{a}hler-Einstein metrics.  
But in general, a smooth complex projective variety does not admit 
a K\"{a}hler-Einstein metric, since the first Chern class is not 
definite in general.   
One way to overcome this defect is  to consider K\"{a}hler-Einstein metrics allowing singularities.

In \cite{tu0}, I have constructed a K\"{a}hler-Einstein current $\omega_{E}$ 
on  smooth minimal algebraic variety $X$ of general type. 
More precisely  there exists a unique closed semipositive current $\omega_{E}$
such that  
\begin{enumerate}
\item[(1)] There exists a nonempty Zariski open subset $U$ of $X$ such that  $\omega_{E}$ is a $C^{\infty}$-K\"{a}hler form on $U$. 
\item[(2)] $-\mbox{Ric}_{\omega_{E}} = \omega_{E}$ holds on $U$.  
\item[(3)] $\omega_{E}$ is absolutely continuous on $X$. 
\end{enumerate}
Later K. Sugiyama proved that there exists a K\"{a}hler-Einstein current 
on the canonical model of general type (\cite{s}).  
Also I have constructed a K\"{a}hler-Einstein current on an arbitrary smooth 
projective variety of general type (without using the finite generation of 
canonical rings) in \cite{KE}.
Hence for smooth projective varieties of general type, we have a substitute 
of a K\"{a}hler-Einstein metric.
 
We note that the above K\"{a}hler-Einstein current $\omega_{E}$ 
on a projective variety $X$ of general type 
have  the following properties :  $h_{E}:= n!\,(\omega_{E}^{n})^{-1} (n = \dim X)$ is a singular hermitian metric  on $K_{X}$ such that the curvature current $\sqrt{-1}\,\Theta_{h_{E}}$ is a closed semipositive current and  
\begin{equation}
H^{0}(X,{\cal O}_{X}(mK_{X})\otimes {\cal I}(h_{E}^{m}))
\simeq H^{0}(X,{\cal O}_{X}(mK_{X}))
\end{equation}
holds for every $m \geqq 1$, i.e., $h_{E}$ is an AZD of $K_{X}$ (cf. 
Definition \ref{AZD} below).  
In other words, $h_{E}$ is a singular hermitian metric which extracts all 
the positivity of $K_{X}$. 
\vspace{2mm} \\ 

\noindent But for smooth projective varieties of non general type, the above results 
do not say anything.
 
\subsection{Iitaka fibration}

The simplest way to squeeze out the positivity of canonical bundles is to 
use the pluricanonical systems.

Let $X$ be a smooth projective variety.   The Kodaira dimension $\mbox{Kod}(X)$
is defined by 
\begin{equation}
\mbox{Kod}(X) := \limsup_{m\rightarrow\infty}\frac{\log \dim H^{0}(X,{\cal O}_{X}(mK_{X}))}{\log m}. 
\end{equation}
It is known that $\mbox{Kod}(X)$ is $-\infty$ or a nonnegative integer 
between $0$ and $\dim X$.  \\ 

\noindent Let $X$ be a smooth projective variety with $\mbox{Kod}(X)\geqq 0$. 
Then for a sufficiently large $m > 0$, the complete linear system 
$|m!K_{X}|$ gives a rational fibration (with connected fibers) :
\begin{equation}
f : X -\cdots\rightarrow  Y. 
\end{equation}
We call $f : X -\cdots\rightarrow  Y$ the {\bf Iitaka fibration} of $X$.
  
The Iitaka fibration is independent of the choice of the sufficiently large $m$
up to birational equivalence.  
In this sense the Iitaka fibration is unique. 
By taking a suitable modification, we may assume that 
$f$ is a morphism and $Y$ is smooth.   

The Iitaka fibration $f : X \to Y$ satisfies the following properties: 
\begin{enumerate}
\item[(1)] For a general fiber $F$, $\mbox{Kod}(F) = 0$ holds.   
\item[(2)] $\dim Y = \mbox{Kod}(Y)$. 
\end{enumerate}

\subsection{Analytic Zariski decompositions}\label{AZD}

Let $L$ be a pseudoeffective line bundle on a compact complex manifold $X$.
To analyze the ring :
\begin{equation}
R(X,L) = \oplus_{m=0}^{\infty}H^{0}(X,{\cal O}_{X}(mL)),
\end{equation}
it is useful to  introduce the notion of analytic Zariski decompositions. 

\begin{definition}\label{defAZD}
Let $M$ be a compact complex manifold and let $L$ be a holomorphic line bundle
on $M$.  A singular hermitian metric $h$ on $L$ is said to be 
an {\bf analytic Zariski decomposition}{\em (}{\bf AZD} in short{\em )}, if the followings hold.
\begin{enumerate}
\item[(1)] $\sqrt{-1}\,\Theta_{h}$ is a closed positive current. 
\item[(2)] for every $m\geq 0$, the natural inclusion:
\begin{equation}
H^{0}(M,{\cal O}_{M}(mL)\otimes{\cal I}(h^{m}))\rightarrow
H^{0}(M,{\cal O}_{M}(mL))
\end{equation}
is an isomorphim. \fbox{}
\end{enumerate}
\end{definition}
\begin{remark} If an AZD exists on a line bundle $L$ on a compact complex manifold $M$, $L$ is pseudoeffective by the condition 1 above. \fbox{}
\end{remark}

It is known that for every pseudoeffective line bundle on a compact complex manifold, there exists an AZD on $L$ (cf. \cite{tu,tu2,d-p-s}). 
The advantage of the AZD is that we can handle pseudoeffective line bundle 
$L$ on a compact complex manifold $X$  
as a singular hermitian  line bundle with semipositive curvature current
as long as we consider the ring $R(X,L)$. 
\vspace{1mm}\\
One may construct an AZD for a pseudoeffective line bundle on 
a compact complex manifold as follows.  
Let $L$ be a  pseudoeffective line bundle on a compact complex manifold $X$.  
 Let $h_{0}$ be a $C^{\infty}$-hermitian metric on $L$. 
We set  
\begin{equation}
h_{min} := \inf \{ h\,|\,\mbox{a singular hermitian metric on $L$},\,\, h\geqq h_{0},
\sqrt{-1}\Theta_{h} \geqq 0\}.
\end{equation}
Then $h_{min}$  is an AZD on $L$ with minimal singularities in the 
following sense. 

\begin{definition}\label{minAZD}
Let $L$ be a pseudoeffective line bundle on a compact complex manifold $X$.
An AZD $h$ on $L$  is said to be an {\bf AZD of minimal singularities}, if 
for any AZD $h^{\prime}$ on $L$, there exists a positive constant $C$ such that
\begin{equation}
h \leqq   C \cdot h^{\prime}
\end{equation}
holds. \fbox{} 
\end{definition}
In general an AZD of a pseudoeffective line bundle $L$ on a smooth projective 
variety is not necessarily of minimal singularities. 

\subsection{Requirement of the canonical semipositive current}
Let $f : X \to Y$ be the Iitaka fibration of a smooth projective 
variety $X$ of nonnegative Kodaira dimension.
In this article, we shall consider a canonical semipositive current, say $\omega_{X}$ associated with the Iitaka fibration.  

It is natural to require that $\omega_{X}$ has the following properties :
\begin{enumerate}
\item[(1)] $\omega_{X}$ is unique and birationally invariant, i.e., 
if $X$ is birational to $X^{\prime}$ and let $\mu : X^{\prime\prime} \to X$ 
and $\mu^{\prime} : X^{\prime\prime} \to X^{\prime}$ be modifications 
from a smooth projective variety $X^{\prime\prime}$.  
Then 
\begin{equation}
\omega_{X^{\prime\prime}}  = \mu^{*}\omega_{X} + 2\pi E = (\mu^{\prime})^{*}\omega_{X^{\prime}}+ 2\pi E^{\prime}
\end{equation}
hold, where $E := K_{X^{\prime\prime}} - 
\mu^{*}K_{X}$ and $E^{\prime}:= K_{X^{\prime\prime}}- (\mu^{\prime})^{*}K_{X^{\prime}}$ respectively.
\item[(2)] There exists an AZD $h_{K}$ of $K_{X}$ such that 
$\omega_{X} = \sqrt{-1}\,\Theta_{h_{K}}$. 
\item[(3)] There exists a closed semipositive current $\omega_{Y}$ on $Y$ such that 
$\omega_{X} = f^{*}\omega_{Y}$,  
\item[(4)] There exists a nonempty Zariski open subset $U$ of $Y$ such that 
$\omega_{Y}|U$ is a $C^{\infty}$-K\"{a}hler form. 
\end{enumerate}

\subsection{Canonical bundle formula and the correction term for 
the K\"{a}hler-Einstein equations}\label{CB}

To construct the AZD $h_{K}$ above, we shall solve a partial differential equation.  The equation is similar to the K\"{a}hler-Einstein equation :
$-\mbox{Ric}_{\omega_{E}} = \omega_{E}$, but there are two major differences :
\begin{enumerate}
\item[(1)] The equation is defined on $Y$ not on $X$. 
\item[(2)] The equation has the additional term 
which comes from  variation of Hodge structures.  
\end{enumerate}

For a graded ring $R : = \oplus_{i=0}^{\infty}R_{i}$
 and a positive integer $m$,  we set 
\begin{equation}
R^{(m)} : = \oplus_{i=0}^{\infty}R_{mi}.
\end{equation}
For a KLT pair $(M,D)$,  we set 
\begin{equation}\label{RM}
R(M,K_{M}+D) := \oplus_{m=0}^{\infty}
\Gamma (M,{\cal O}_{M}(\lfloor m(K_{M}+D)\rfloor))
\end{equation}
and 
\begin{equation}
\mbox{Kod}(M,D) := \limsup_{m\rightarrow\infty}
\frac{\log \dim \Gamma (M,{\cal O}_{M}(\lfloor m(K_{M}+D)\rfloor ))}{\log m}. 
\end{equation}
\begin{theorem}(\cite[p.183, Theorem 5.2]{f-m})\label{logcan}
Let $(X,\Delta)$ be a proper KLT pair with 
\begin{equation}
\mbox{\em Kod}(X,K_{X}+\Delta) = n.
\end{equation}
Then there exists a $n$-dimensional KLT pair $(Y^{\prime},\Delta^{\prime})$
with $\mbox{\em Kod}(Y^{\prime},\Delta^{\prime}) = n$, two positive integers
$e,e^{\prime}$ such that 
\begin{equation}
R(X,K_{X}+\Delta)^{(e)}\simeq R(Y^{\prime},K_{Y^{\prime}}+\Delta^{\prime})^{(e^{\prime})}. 
\end{equation}
\fbox{}
\end{theorem}

Let us consider the case that $\Delta = 0$ in Theorem \ref{logcan}. 
Then the  canonical ring $R(X,K_{X})^{(e)}$ is (a subring of) the pullback of the log canonical ring of some  KLT pair $(Y^{\prime},\Delta^{\prime})$ 
of log general type.

Let us explain the equation. 
Let $f : X \to Y$ be an  Iitaka fibration such that 
$f_{*}{\cal O}_{X}(m!K_{X/Y})^{**}$  is locally free on $Y$ for some $m$ (hence for every  sufficiently large $m$), where $**$ denotes the double dual. 
Such $f: X \to Y$ exists by \cite[p.169,Proposition 2.2]{f-m}. 
The divisor $\Delta^{\prime}$ is related  to the $\mathbb{Q}$-line bundle  
\begin{equation}\label{lxy}
L_{X/Y} := \frac{1}{m_{0}!}f_{*}{\cal O}_{X}(m_{0}!K_{X/Y})^{**} 
\end{equation}
on $Y$  in terms of the canonical bundle formula (See \cite{f-m} for detail. In \cite[Section 2]{f-m}), where $m_{0}$ is a sufficiently large positive 
integer.   
We note that $L_{X/Y}$ is independent of a sufficiently large $m_{0}$ (cf. \cite[Section 2]{f-m}).   We call $L_{X/Y}$ the {\bf Hodge $\mathbb{Q}$-line bundle} 
of $f : X\to Y$. 
$L_{X/Y}$  carries a natural singular hermitian metric $h_{L_{X/Y}}$ defined by 
\begin{equation}\label{hxy}
h_{L_{X/Y}}^{m_{0}!}(\sigma,\sigma) = \left(\int_{X_{y}}|\sigma|^{\frac{2}{m_{0}!}}\right)^{m_{0}!},
\end{equation}
where $y \in Y, X_{y}:= f^{-1}(y)$ and $\sigma \in L_{X/Y,y}$. 
It is known that $h_{L_{X/Y}}$ has semipositive curvature in the sense of current
(\cite{ka1},\cite[p.174,Theorem 1.1]{ka3}) by using the variation of Hodge structures.  By a suitable modification of $f : X \to Y$, we may assume that there exists an decomposition:
\begin{equation}\label{MD-dec}
L_{X/Y} = M_{X/Y} + D_{X/Y}
\end{equation}
such that 
\begin{enumerate}
\item[(1)] $M_{X/Y},D_{X/Y}\in \mbox{Div}(Y)\otimes \mathbb{Q}$, 
\item[(2)] $M_{X/Y}$ is nef, $D_{X/Y}$ is effective. 
\item[(3)] $\mbox{Supp}\, D_{X/Y}$ is a divisor with normal crossings on 
$Y$,
\item[(4)] $h_{L_{X/Y}}$ is induces a singular hermitian metric  
$h_{M_{X/Y}}$ on $M_{X/Y}$ with semipositive curvature and is logarithmic growth as Lemma \ref{logg} below. and $D_{X/Y}$ corresponds to the singular part 
of the 1-st Chern current $c_{1}(L_{X/Y},h_{L_{X/Y}}))$.   
\item[(5)] $(Y,D_{X/Y})$ is KLT. 
\end{enumerate}
We call $M_{X/Y}$ the {\em semistable part} and $D_{X/Y}$ the 
{\em discriminant part} of $L_{X/Y}$ respectively (\cite{f-m}). 

Let $\Omega$ be a $C^{\infty}$-volume form on $Y$. 
We shall consider the following equation : 
\begin{equation}\label{LKE}
-\mbox{Ric}_{\omega_{Y}} + \sqrt{-1}\,\Theta_{h_{L_{X/Y}}}  =  \omega_{Y},
\end{equation}
where 
\begin{equation}  
\omega_{Y}  =  -\mbox{Ric}\,\Omega + \sqrt{-1}\,\Theta_{h_{L_{X/Y}}} + 
\sqrt{-1}\partial\bar{\partial}u, 
\end{equation}
and $u$ is the unknown function. 
Here the term $\sqrt{-1}\,\Theta_{h_{L_{X/Y}}}$ corresponds to the boundary divisor 
of the KLT pair $(Y^{\prime},\Delta^{\prime})$ in the canonical bundle formula(\cite[p.183,Theorem 5.2]{f-m}).
The authenticity of the equation (\ref{LKE}) can be 
verified  by checking the fact that the dynamical system of Bergman kernels 
on $X$ as in \cite{KE} yields  the current on $Y$ which satisfies 
the equation (\ref{LKE}) (cf. Theorem \ref{DS}).     

Actually I first constructed the current by using the dynamical system of 
Bergman kernels and then found the equation (\ref{LKE}) inspired by 
the \cite{f-m}. 
   
There are several difficulties to solve the equation (\ref{LKE}).
First of all we cannot expect that there exists a $C^{\infty}$- solution
 $\omega_{Y}$. 
In fact $\Theta_{h_{L_{X/Y}}}$ is not $C^{\infty}$ in general.  And 
if $\Theta_{h_{L_{X/Y}}}$ is not $C^{\infty}$, 
(\ref{LKE}) has no $C^{\infty}$-solution $\omega_{Y}$.  
Moreover $h_{L_{X/Y}}$ is not of algebraic singularities (cf. Definition \ref{algsing}) in general.   And even if $h_{L_{X/Y}}$ is $C^{\infty}$, $u$ is not $C^{\infty}$ in general. 
Secondary the solution $u$ is not unique.  
But if we require that $\Omega^{-1}\cdot h_{L_{X/Y}}\cdot e^{u}$ 
is an AZD of $K_{Y}+ L_{X/Y}$, then the solution $u$ is actually unique 
and the resulting current $\omega_{Y}$ is nothing but the current 
constructed by the dynamical system of Bergman kernels (see Section \ref{Ds} and Section \ref{Dy} below).
 
\subsection{Canonical K\"{a}hler currents}
Now we shall state the existence of the canonical semipositive current
on a smooth projective variety of nonnegative Kodaira dimension. 
\begin{theorem}\label{main}(cf. \cite[Theorem B.2]{s-t}) 
In the above notations, there exists a unique singular hermitian metric 
on $h_{K}$ on $K_{Y}+L_{X/Y}$ and a nonempty Zariski open subset 
$U$ in $Y$ such that 
\begin{enumerate}
\item[(1)] $h_{K}$ is an AZD of $K_{Y} + L_{X/Y}$,
\item[(2)] $f^{*}h_{K}$ is an AZD of $K_{X}$,  
Here we have used the inclusion: 
\[
f^{*}\mathcal{O}_{Y}(m_{0}!(K_{Y} + L_{X/Y}))
\hookrightarrow \mathcal{O}_{X}(m_{0}!K_{X})
\]  
to identify $f^{*}h_{K}$ with a singular hermitian metric on $K_{X}$, 
where $m_{0}$ is the sufficiently large positive integer used in 
(\ref{lxy}), 
\item[(3)] $h_{K}$ is $C^{\infty}$ on $U$,
\item[(4)] $\omega_{Y} = \sqrt{-1}\,\Theta_{h_{K}}$ is a K\"{a}hler form 
on $U$,
\item[(5)] $-\mbox{\em Ric}_{\omega_{Y}} + \sqrt{-1}\,\Theta_{L_{X/Y}} = \omega_{Y}$
holds on $U$. \fbox{} 
\end{enumerate} \vspace{2mm}
\end{theorem}
Although the proof of Theorem \ref{main} is given in \cite{s-t}, 
I shall give an alternative proof in this paper for the completeness, 
since my original proof seems to be different from that in \cite{s-t}.
One can see the KLT version of the above theorem in \cite{LC}. 

\begin{definition}\label{L-K-E}
The current $\omega_{Y}$ on $Y$ is said to be {\bf the canonical K\"{a}hler current}  of the Iitaka fibration $f : X \to Y$.  Also $\omega_{X}:= f^{*}\omega_{Y}$ is said to be {\bf the canonical semipositive current} on $X$. 
We define the measure $d\mu_{can}$ on $X$ by 
\begin{equation}
d\mu_{can}:= \frac{1}{n!}f^{*}\omega_{Y}^{n}\cdot h_{L_{X/Y}}^{-1}
\end{equation}
and is said to be {\bf the canonical measure}, where $n$ denotes 
$\dim Y$.  
\fbox{} 
\end{definition}

\noindent The existence of the canonical K\"{a}hler current is proven 
in terms of solving Monge-Amp\`{e}re equations. 
The proof given here is similar to the one of 
\cite[Section 5.1, Theorem 5.1]{KE}. 
We shall give a proof in Section 2 (see also \cite[Section 4]{s-t}). 

\subsection{Ricci iterations}\label{IT}

Let $f : X \to Y$, $(L_{X/Y},h_{L_{X/Y}})$ be as above.  
We assume that there exists a Zariski decomposition :
\begin{equation}\label{ZD2}
K_{Y} + L_{X/Y} = P + N  \,\,\,\,\,(P,N \in \mbox{Div}(X)\otimes \mathbb{Q})
\end{equation}
i.e.,  $P$ is semiample, $N$ is effective  and 
\begin{equation}
H^{0}(Y,\mathcal{O}_{X}(maP))
\simeq H^{0}(X,\mathcal{O}(ma(K_{X} + L_{X/Y})))
\end{equation}
holds for every $m\geqq 0$, where $a$ is the minimal positive integer 
such that $aL_{X/Y}, aP, aN \in \mbox{Div}(Y)$. 
\vspace{2mm} \\

\noindent Let $h_{P}$ be a $C^{\infty}$-hermitian metric on $P$ with semipositive curvature.  
We set 
\begin{equation}\label{u}
U : = \{ x\in X|\,\,\,\mbox{$|\nu!P|$ is  very ample on a neighbourhood of 
$x$ for every $\nu \gg 1$}\}.
\end{equation} 
Let $a$ be a positive integer such that $aD\in \mbox{Div}(X)$. 
For $m\geqq 0$, we shall define inductively a sequence of closed positive current $\{\omega_{m}\}_{m=0}^{\infty}$ satisfying the following conditions : 
\begin{itemize}
\item[(P1)] $\omega_{0} = a\left(\sqrt{-1}\Theta_{h_{P}} + 2\pi[N]\right)$, 
where $[N]$ denotes the current of integration over $N$.  
\item[(P2)] The cohomology class of $[\omega_{m}]$ of $\omega_{m}$ 
is equal to $2\pi a\cdot c_{1}(K_{Y}+L_{X/Y})$.
\item[(P3)] $\omega_{m}$ is $C^{\infty}$ on $U$.
\item[(P4)] $\{\omega_{m}\}_{m=0}^{\infty}$ satisfies the equation :
\begin{equation}\label{ricciind2}
-\mbox{Ric}_{\omega_{m}} + \frac{a - 1}{a}\,\omega_{m-1} + 
\sqrt{-1}\Theta_{h_{L_{X/Y}}} = \omega_{m}
\end{equation}
on $U$ for every $m\geqq 1$.
\item[(P5)] We define the singular hermitian metric $h_{m}$ by 
\begin{equation}
h_{m}:= \left(\frac{1}{n!}\omega_{m,abc}^{n}\right)^{-\frac{1}{a}}\cdot h_{m-1}^{\frac{a-1}{a}}\cdot h_{L_{X/Y}}^{\frac{1}{a}}
\end{equation}
on $K_{Y}+L_{X/Y}$. 

Then $h_{m}$ is an AZD of $K_{Y}+L_{X/Y}$ for every $m\geqq 1$. 
\end{itemize}

\begin{theorem}\label{iteration}
The dynamical system $\{ \omega_{m}\}_{m=0}^{\infty}$ 
satisfying (P1)-(P5) above exists 
and is unique and the limit 
\begin{equation}
\omega_{K} := \frac{1}{a}\lim_{m\rightarrow\infty}\omega_{m}
\end{equation}
exists in $C^{\infty}$-topology on every compact subset of $U$ 
(cf.{\rm (\ref{u})}) and also in the sense of current on $Y$. 
The closed positive current $\omega_{Y}$ satisfies the equation :
\[
-\mbox{\em Ric}_{\omega_{Y}} + \sqrt{-1}\Theta_{h_{L_{X/Y}}} = \omega_{Y}
\]
on $X$ and  $\omega_{Y}$ is the canonical K\"{a}hler current on 
$Y$ (cf. Definition \ref{LKE}).  \fbox{}
\end{theorem}
The proof of Theorem \ref{iteration} follows from  
the one of \cite[Theorem 3.2]{LC}.  Hence we omit it.   

\subsection{Dynamical construction of the canonical K\"{a}hler currents}\label{Ds}

The canonical K\"{a}hler current  in Theorem \ref{main} 
can be constructed as the limit of a dynamical system as in (\cite{KE,LC}). 
Actually we decompose the Ricci iteration in Theorem \ref{iteration} 
as a sequence of dynamical systems of Bergman kernels.  
The motivation to construct such a dynamical systems is to prove 
the plurisubharmonic variation property of canonical measures (cf. Theorem \ref{relative}).   

Let $X, Y, f : X \to Y, (L_{X/Y},h_{L_{X/Y}})$ as in Section \ref{CB},   
Let $a$ be a positive integer such that $f_{*}{\cal O}_{X}(aK_{X/Y})
\neq 0$. Then  we see that 
\begin{equation}\label{iso}
H^{0}(X,{\cal O}_{X}(maK_{X}))\simeq H^{0}(Y,{\cal O}_{Y}(ma(K_{Y}+L_{X/Y})))
\end{equation}
holds for every $m \geqq 0$.   In particular $\mbox{Kod}(X) = \dim Y$ holds. 
Hence by  (\ref{iso}), we see that $K_{Y} + L_{X/Y}$ is big. 

Let $A$ be a sufficiently ample line bundle on $Y$ and let $h_{A}$ be a $C^{\infty}$-hermitian metric on $A$ with strictly positive 
curvature.   

We shall construct the singular hermitian metrics $\{ h_{m}\}_{m=1}^{\infty}$ 
on $K_{Y} + L_{X/Y}$ inductively as follows.  We shall abuse the 
same notation as in Section \ref{IT}, since $\{ h_{m}\}$ is eventually 
the same as Theorem \ref{iteration} (cf. Theorem \ref{DS} below). 

Let $K_{Y} + L_{X/Y} = P + N$ be the Zariski decomposition as (\ref{ZD2}). 
First we set $h_{1} = h_{P}\cdot h_{N}$ as in Section \ref{IT}, where 
$h_{P}$ is a $C^{\infty}$-hermitian metric on $P$ with semipositive curvature and $h_{N}$ be a metric on $N$ defined as $h_{N}:= |\sigma_{N}|^{-2}$, 
where $\sigma_{N}$ be a multivalued holomorphic section of $N$ with divisor 
$N$.   
 
Suppose that we have constructed the metric 
$h_{m-1}$ on $K_{Y}+L_{X/Y}$ for some $m\geqq 2$. 
To construct $h_{m}$, we shall construct a sequece of 
singular hermitian metrics $\{ h_{\ell,m}\}_{\ell=1}^{\infty}$ as 
follows.  

First we set     
\begin{equation}
K_{1,m} :=  K(Y,A+a(K_{Y} +L_{X/Y}),h_{A}\cdot h_{m-1}^{a-1}\cdot h_{L_{X/Y}}).
\end{equation}
and  
\begin{equation}
h_{1,m} := (K_{1,m})^{-1}. 
\end{equation}
We continue this process. 
Suppose that we have constructed 
$K_{\ell,m}$ and the singular hermitian metric $h_{\ell,m}$ on 
$A+\ell a(K_{Y}+L_{X/Y})$ for some $\ell \geqq 1$. 
\begin{equation}
K_{\ell+1,m}:= K(Y,A+(\ell+1)a(K_{Y}+L_{X/Y}),h_{\ell,m}\cdot h_{L_{X/Y}})
\end{equation}
and 
\begin{equation}
h_{\ell+1,m}:= (K_{\ell+1,m})^{-1}.  
\end{equation}
Thus inductively we construct the sequences $\{ h_{\ell,m}\}_{\ell >0}$
and $\{ K_{\ell,m}\}_{m > 0}$.
This inductive construction is essentially the same one originated by the author in \cite{tu3}.  This is the same construction as in \cite{LC}.  
The following theorem asserts that the above dynamical system yields the 
 Ricci iteration constructed in Section \ref{IT}. 

\begin{theorem}\label{DS}
 Let $X$ be a smooth projective variety of nonnegative Kodaira dimension 
 and let $f : X \to Y$ be the Iitaka fibration as above.  
 Let  $\{ h_{\ell,m}\}_{\ell>0}$ be the sequence 
of hermitian metrics as above and let $n$ denote $\dim Y$.
Let $\omega_{m}$ be the  K\"{a}hler current on $Y$  
defined as  in Theorem \ref{IT}.  
Then 
\begin{equation}\label{hinfty}
h_{m} := \liminf_{\ell\rightarrow\infty} \sqrt[a\ell]{(\ell!)^{n}\cdot h_{\ell,m}}
\end{equation}
is a singular hermitian metric on $K_{Y}+L_{X/Y}$ such that 
\begin{equation}\label{hinftyvol}
h_{m} = \left(\frac{1}{n!}\omega_{m,abc}^{n}\right)^{-\frac{1}{a}}\cdot h_{m-1}^{\frac{a-1}{a}}\cdot h_{L_{X/Y}}^{\frac{1}{a}}
\end{equation}
holds almost everywhere on $Y$ and 
\begin{equation}\label{hinftycurvature}
\omega_{m}:= a\sqrt{-1}\,\Theta_{h_{m}}. 
\end{equation}
holds on $Y$.   
In particular $h_{m}$ (and hence $\omega_{m}$) is unique and is independent of the choice of $A$ and $h_{A}$. 
\fbox{}
\end{theorem}
Here it may be better to replace $h_{m}$ by its lower-semi-continuous 
envelope because of the following classical theorem. 
\begin{theorem}(\cite[p.26, Theorem 5]{l})\label{Lelong}
Let $\{ u_{\alpha}\}_{\alpha\in A}$ be a family of plurisubharmonic function on  a domain $\Omega$ in $\mathbb{C}^{n}$. 
Suppose that $\{ u_{\alpha}\}_{\alpha\in A}$ 
is locally uniformly bounded from above.   
Then  the upper-semi-continuous envelope of $\sup_{\alpha\in A}u_{\alpha}$
is again plurisubharmonic on $\Omega$ \fbox{}
\end{theorem} 
But anyway this adjustment occurs on a set of measure $0$.

\section{Construction of the canonical measure by solving Monge-Amp\`{e}re equations}

In this section we shall prove Theorem \ref{main}.   
Although the proof of Theorem \ref{main} is given in \cite{s-t}
independently, we shall 
give an orginal proof here for the completeness.  
Also the proof in \cite{s-t} is quite close to the proof of 
\cite[Theorem 5.1]{KE}\footnote{The orginal proof of \cite[Theorem 5.1]{KE} has  a gap in  the proof of $C^{2}$-regularity.  I have not checked the proof in \cite{s-t} in full detail.}.   
The present proof is different from that in \cite{s-t}, in the 
following points :
\begin{enumerate}
\item[(1)] We consider a smoothing of the Hodge metric $h_{L_{X/Y}}$ which need not be of algebraic singularities. 
Hence we need to consider the sequence of modified equations.  
\item[(2)] The $C^{0}$-estimate of the solution depends on the estimate of 
the Hodge metric near the discriminant locus and the notion of 
the minimal AZD. 
\end{enumerate}
The techniques used here are quite standard  and have been known 
for more than twenty years (cf. \cite{s,tu0}).   In this sense the proof
of Theorem \ref{main} is not essentially new.   
But as in \cite{s}, we need to require the finite generation 
of canonical ring (\cite{b-c-h-m}) to prove the $C^{2}$-regularity 
of the metrics on a Zariski open subset.  

\subsection{Setup} 
Let $X$ be a smooth projective $n$-fold with $\mbox{Kod}(X) \geqq 0$.
And let 
\begin{equation}
f : X -\cdots\rightarrow Y 
\end{equation}
be the Iitaka fibration associated with the complete linear system 
$|m_{0}!K_{X}|$.  By taking a suitable modifications, we shall assume 
the followings.
\begin{enumerate}
\item[(1)] $f$ is a morphism.
\item[(2)] $Y$ is smooth. 
\item[(3)] $f_{*}{\cal O}_{X}(m_{0}!K_{X/Y})^{**}$ is a line bundle on $Y$. 
\item[(4)] The discriminant locus $D$ of $f$ is a divisor with normal crossings
on $Y$. 
\end{enumerate}  
We  define the $\mathbb{Q}$-line bundle $L_{X/Y}$ on $Y$ by  
\begin{equation}
L_{X/Y} := \frac{1}{m_{0}!}f_{*}{\cal O}_{X}(m_{0}!K_{X/Y})^{**}
\end{equation}
and let  $h_{L_{X/Y}}$ be the singular hermitian metric on $L_{X/Y}$ defined by 
\begin{equation}
h_{L_{X/Y}}^{m_{0}!}(\sigma,\sigma):= \left(\int_{X_{y}}|\sigma|^{\frac{2}{m_{0}!}}\right)^{m_{0}!}. 
\end{equation}
It is clear that $h_{L_{X/Y}}$ is smooth on 
\begin{equation}
Y^{\circ}:= \{ y\in Y| \mbox{$f$ is smooth over $y$}\} = Y \backslash D
\end{equation}
and the singularity of $h_{L_{X/Y}}$  around $D$ is described in terms of  variation of Hodge structures.  And we see that 
\begin{equation}
H^{0}(Y,{\cal O}_{Y}(m!(K_{Y}+L_{X/Y}))\otimes {\cal I}(h_{L_{X/Y}}^{\otimes m!}))
\simeq H^{0}(Y,{\cal O}_{Y}(m!(K_{Y}+L_{X/Y})))
\end{equation}
holds for every sufficiently large $m$,i.e., the $L^{2}$-condition with respect to the singular hermitian metric $h_{L_{X/Y}}$ does not affect  the global 
section of $m!(K_{Y}+L_{X/Y})$.  

Let $\Omega$ be a $C^{\infty}$-volume form on $Y$. 
Let us consider the equation : 
\begin{equation}\label{equ}
-\mbox{Ric}_{\omega_{Y}} + \sqrt{-1}\,\Theta_{h_{L_{X/Y}}} = \omega_{Y}
\end{equation}
on $Y$, where 
\begin{equation}
\omega_{Y} = -\mbox{Ric}\,\Omega + \sqrt{-1}\,\Theta_{h_{L_{X/Y}}} + \sqrt{-1}\partial\bar{\partial}u  
\end{equation}
for some unknown upper-semi-continuous function $u$ bounded from above on $Y$. 
Then the above equation is equivalent to 
\begin{equation}
\log \frac{(\omega + \sqrt{-1}\partial\bar{\partial}u)^{n}}{\Omega} 
= u,
\end{equation}
where $n:= \dim Y$ and 
\begin{equation}
\omega :=  -\mbox{Ric}\,\Omega + \sqrt{-1}\,\Theta_{h_{L_{X/Y}}}.
\end{equation}
We note that since $\mbox{Kod}(X) = \dim Y$, 
$K_{Y} + (L_{X/Y},h_{L_{X/Y}})$ is big, i.e., 
\begin{equation}
\limsup_{m\rightarrow\infty}
(m!)^{-n}\dim H^{0}(Y,{\cal O}_{Y}(m!(K_{Y}+L_{X/Y}))\otimes {\cal I}(h_{L_{X/Y}}^{m!}))
> 0 
\end{equation} 
holds, where $n = \dim Y$. 
The main difficulty for solving the equation  (\ref{equ}) is the fact that 
$h_{L_{X/Y}}$ is not of algebraic singularities in the following sense.
\begin{definition}\label{algsing}
Let $h$ be a singular hermitian metric on a line bundle $L_{X/Y}$. 
We say that $h$ is {of \bf algebraic singularities}, if there exists a positive  integer $m_{0}$,  global holomorphic sections $\sigma_{0},\cdots ,\sigma_{N}$
of $m_{0}L_{X/Y}$ and a $C^{\infty}$-function $\phi$ such that 
\begin{equation}
h= e^{\phi}\cdot(\sum_{i=0}^{N}\mid\sigma_{i}\mid^{2})^{-\frac{1}{m_{0}}}
\end{equation}
holds.   \fbox{}
\end{definition}
There are two ways to treat the singularity of $h_{L_{X/Y}}$ in (\ref{equ}). 
One way is to smooth out the singularities of $h_{L_{X/Y}}$ and the other way 
is to consider the metric of Poincar\'{e} growth on the complement 
of the discriminant locus of $f : X \to Y$. 
The both methods depend on the analysis of the singularties of $h_{L_{X/Y}}$ 
in terms of the theory of variation of Hodge structures.

\subsection{Smoothing of the Hodge metric $h_{L_{X/Y}}$}\label{smoothing}

The singular hermitian metric $h_{L_{X/Y}}$ is generically $C^{\infty}$, but need not be of algebraic singularities. 

The singularity of $h_{L_{X/Y}}$ can be described by using variation of 
Hodge structures.   
Let $a$ be the minimal positive integer such that
$f_{*}{\cal O}_{X}(aK_{X/Y})^{**}$ is not $0$. 
Then the  $a$-th root of  local holomorphic 
section of $f_{*}{\cal O}_{X}(aK_{X/Y})$ can be considered to be 
a family of canonical forms on the family of cyclic $a$-covers of 
the fibers.  
In this way the Hodge metric can be described in terms of the theory 
of  variation of Hodge structures (cf. \cite{sch}). 
Let us assume that the discriminant locus $D$ of $f : X \to Y$ 
is a divisor with normal crossings. 
As in \cite{kawa}, the locally free extension of the Hodge bundle is contorolled  by the monodoromy which is quasi-unipotent.   
And in this setting the local monodoromy is abelian.

\begin{definition}\label{P-growth}
Let $(M,B)$ a pair of a complex maifold $M$ and a divisor $B$ with normal
crossings.   Let $\omega_{P}$ is a K\"{a}hler form on $M- B$. 
$\omega_{P}$ is said to be of Poincar\'{e} growth, if 
for any polydisk $\Delta^{n} := \{ (z_{1},\cdots,z_{n}); 
|z_{i}| < 1, 1\leqq i\leqq n\}$ in $M$ such that 
\begin{equation}
\Delta^{n}\cap B = \{(z_{1},\cdots,z_{n})\in \Delta^{n}|z_{1}\cdots z_{k} = 0\},
\end{equation}
there exist locally bounded positive continuous functions $a \leqq b$ on $\Delta^{n}$
such that  
\begin{equation}
a\left(\sum_{i=1}^{k}\frac{\sqrt{-1}\,dz_{i}\wedge d\bar{z}_{i}}{|z_{i}|^{2}(\log |z_{i}|)^{2}} + \sum_{j=k+1}^{n}\sqrt{-1}\,dz_{j}\wedge d\bar{z}_{j}\right)
\leqq \omega_{P}
\leqq b\left(\sum_{i=1}^{k}\frac{\sqrt{-1}\,dz_{i}\wedge d\bar{z}_{i}}{|z_{i}|^{2}(\log |z_{i}|)^{2}}+ \sum_{j=k+1}^{n}\sqrt{-1}\,dz_{j}\wedge d\bar{z}_{j}\right)
\end{equation}
hold on $\Delta^{n}\cap (M-B)$. 

Let $\Omega_{P}$ be a volume form on $M - B$.  $\Omega_{P}$ is said to be 
of Poincar\'{e} growth, if for any polydisk $\Delta^{n}$ in $M$
such that
\begin{equation}
\Delta^{n}\cap B = \{(z_{1},\cdots,z_{n})\in \Delta^{n}|z_{1}\cdots z_{k} = 0\},
\end{equation}
there exists a locally bounded positive continuous function $c(z)$ on $\Delta^{n}$
such that  
\begin{equation}
\Omega_{P}= 
c(z)\cdot\left(\wedge_{i=1}^{k}\frac{\sqrt{-1}\,dz_{i}\wedge d\bar{z}_{i}}{|z_{i}|^{2}(\log |z_{i}|)^{2}}\right)\wedge \left(\wedge_{j=k+1}^{n}\sqrt{-1}\,dz_{j}\wedge d\bar{z}_{j}\right)
\end{equation}
holds on $\Delta^{n}\cap (M-B)$. \fbox{} 
\end{definition}
\begin{remark}\label{ric}
Let $\Omega_{P}$ be a volume form of Poincar\'{e} growth on 
$(M,B)$ with $M$ compact.  If for every polydisk $\Delta^{n}$ as in Definiton \ref{P-growth}  
the function $c(z)$ above is $C^{2}$ on $\Delta^{n}$, then $-\mbox{\em Ric}\,\Omega_{P}$ is of Poincar\'{e}
growth in the sense that there exists a positive constant $C$ such that 
\begin{equation}
- C\cdot \omega_{P} \leqq -\mbox{\em Ric}\,\Omega_{P} \leqq 
C\cdot \omega_{P},
\end{equation}
where $\omega_{P}$ is a K\"{a}hler form of $M - B$ with Poincar\'{e}
growth.  This is standard and is easily verified by direct calculation.  \fbox{}
\end{remark}
Then we have the following lemma. 

\begin{lemma}\label{Poincare} There exists a positive integer $m_{0}$ such that 
$h_{M_{X/Y}}^{m_{0}!}$ is a singular hermitian metric 
on the line bundle $f_{*}{\cal O}_{Y}(m_{0}!K_{X/Y})^{**}$ such that 
with respect to a local holomorphic frame 
\begin{equation}\label{logg}
h_{M_{X/Y}} =  O((-\log |\sigma_{D}|)^{q})
\end{equation}
holds, where $\sigma_{D}$ is a local defining function of 
$D$ and $q$ is a positive integer.  Hence we have that the estimate
\begin{equation}\label{logg2}
h_{L_{X/Y}} =  O(|\sigma_{D}|^{-2\alpha}(-\log |\sigma_{D}|)^{q})
\end{equation} 
holds for some semiposive rational number $\alpha < 1$. 
 And the curvature $\sqrt{-1}\,\Theta_{h_{M_{X/Y}}}$
is dominated by a constant times a K\"{a}hler form $\omega_{P}$ with Poincar\'{e} growth on $Y \backslash D$.  \fbox{}
\end{lemma}
\begin{remark}\label{belowe}
Besides (\ref{logg}), we see that $h_{L_{X/Y}}$ is bounded from below 
by a smooth metric on $L_{X/Y}$, since $h_{L_{X/Y}}$ has semipositive 
curvature in the sense of current (\cite{ka1},\cite[p.174,Theorem 1.1]{ka3}). 
This fact is essentially due to the theory of variation of Hodge structures.
\fbox{}
\end{remark}
\noindent The estimate of $h_{L_{X/Y}}$ in Lemma \ref{Poincare} follows from \cite{kawa} which uses the  theory variation of Hodge structures due to W. Schmidt (\cite{sch}). 
And the latter estimate of $\sqrt{-1}\,\Theta_{h_{L_{X/Y}}}$ follows from the 
fact that holomorphic sectional curvature in the  horizontal 
direction of the period domain  is dominated by a negative constant (\cite{griff}) and the Yau-Royden Schwarz lemma (\cite{y2,r}).   

In this sense $h_{L_{X/Y}}$ is very close to a smooth metric.   
To smooth out $h_{L_{X/Y}}$, we take a finite open covering ${\cal U} := \{ U_{\alpha}\}$ of $Y$ such that every $U_{\alpha}$ is biholomorphic to the open unit ball in  $\mathbb{C}^{n}$ with ceneter $O$  via the coordinate $z_{\alpha} = (z_{\alpha}^{1},\cdots ,z_{\alpha}^{n})$.   Taking  ${\cal U}$ properly, 
 we may and do assume 
$z_{\alpha}$ is a holomorphic coordinate on a larger open subset $\hat{U}_{\alpha}$ which is biholomorphic to the open ball with radius $2$ in $\mathbb{C}^{n}$ with center $O$ via $z_{\alpha}$.  Let $h_{0}$ be a $C^{\infty}$-hermitian metric on 
the $\mathbb{Q}$ line bundle $L_{X/Y}$.     
We set 
\begin{equation}
\varphi : = \log \frac{h_{L_{X/Y}}}{h_{0}}.
\end{equation}

\noindent Let $\rho$ be a $C^{\infty}$-function on $\mathbb{C}^{n}$ 
 such that $0 \leqq \rho \leqq 1$, $\mbox{supp}\,\rho$ is contained in the unit 
 open ball in $\mathbb{C}^{n}$ with center $O$ 
and 
\begin{equation}
\int_{\mathbb{C}^{n}}\rho(z)\, d\mu(z) = 1, 
\end{equation}
where $d\mu$ is the usual Lebesgue measure on $\mathbb{C}^{n}$. 
For every $0< \delta  < 1$, we set 
\begin{equation}
\rho_{\delta}(z) = \delta^{-2n} \rho(z/\delta).
\end{equation}
We shall take a molification $\varphi_{\alpha,\delta}$ of $\varphi|U_{\alpha}$ 
using the convolution with the molifier $\rho_{\delta}$ as  
\begin{equation}
\varphi_{\alpha,\delta}:= (\varphi|U_{\alpha}) *\rho_{\delta}
\end{equation}
with respect to the coordinate $z_{\alpha}$.  
Since $z_{\alpha}$ is a holomorphic coordinate on $\hat{U}_{\alpha}$, 
$\varphi_{\alpha,\delta}$ is a well defined $C^{\infty}$-function on $U_{\alpha}$ 
for every $0 < \delta < 1$. 
Then $\varphi_{\alpha,\delta}$ converges to $\varphi$ in $L^{1}$-topolgy 
on $U_{\alpha}$ and compact uniformly in 
$C^{\infty}$-topology on $U_{\alpha}\backslash D$  as $\delta \downarrow 0$.

Let $\{\phi_{\alpha}\}$ be a partition of unity subordinate to ${\cal U}$.
We set 
\begin{equation}\label{defsmoothing}
h_{L_{X/Y},\delta} := \exp (\sum_{\alpha}\phi_{\alpha}\cdot \varphi_{\alpha,\delta})\cdot h_{0}. 
\end{equation}
Then $h_{L_{X/Y},\delta}$ is a $C^{\infty}$-hermitian metric on $L_{X/Y}$ and
there exists a positive constant $C$ independent of $\delta > 0$ such that  
\begin{equation}
\sqrt{-1}\,\Theta_{h_{L_{X/Y},\delta}} \leqq C\cdot \omega_{P}
\end{equation}
holds.    In general $h_{L_{X/Y},\delta}$ does not have semipositive 
curvature.    But by the construction, there exists a continuous function $e(\delta )$ on $Y$ 
such that 
\begin{equation}
\sqrt{-1}\,\Theta_{h_{L_{X/Y},\delta}} \geqq - e(\delta)\cdot \omega_{P}
\end{equation}
and 
\begin{equation}
\lim_{\delta\downarrow 0} e(\delta) = 0
\end{equation}
holds uniformly on $Y$. 

\subsection{The construction of the canonical K\"{a}hler currents}\label{exist}

In this subsection, we shall prove the existence of the current $\omega_{Y}$ 
satisfying (\ref{equ}) without assuming the finite generation 
of canonical ring (\cite{b-c-h-m}).  This result is slightly weaker than Theorem \ref{main}. 
But the same strategy works to construct  canonical 
K\"{a}hler-Einstein currents on LC pairs (cf. \cite{LC}). 
Hence the following theorem has independent interest.    

\begin{theorem}\label{current}
In the  notations in Section 2.1, there exists a closed positive current
$\omega_{Y}$ on $Y$ such that  
\begin{enumerate}
\item[(1)] $\omega_{Y}$  represents $2\pi c_{1}(K_{Y} + L_{X/Y})$.
\item[(2)]
$h_{K}:= n!\left(\omega_{Y,abc}^{n}\right)^{-1}\!\!\!\!\cdot h_{L_{X/Y}}
(n = \dim Y)$ 
is an AZD of $K_{Y} + L_{X/Y}$, where $\omega_{Y,abc}$ denotes the 
absolutely continuous part of $\omega_{Y}$.
\item[(3)] $\omega_{Y} = \sqrt{-1}\,\Theta_{h_{K}}$ holds on $Y$.  
\item[(4)] $-\mbox{\em Ric}_{\omega_{Y}} + \sqrt{-1}\,\Theta_{L_{X/Y}} = \omega_{Y}$
holds on $Y$ in the sense of current, where 
$\mbox{\em Ric}_{\omega_{Y}}:= \sqrt{-1}\partial\bar{\partial}\log \omega^{n}_{Y,abc}$. \fbox{} 
\end{enumerate} \vspace{2mm}
\end{theorem}   
The proof of Theorem \ref{current} depends on the monotonicity lemma 
(Lemma \ref{monotonicity}). 

Let $m_{0}$ be a sufficiently large positive integer such that
for every $m \geqq m_{0}$, $m!(K_{Y}+L_{X/Y})$ is  a Cartier divisor on $Y$ and 
$\mid\!\!m!(K_{Y}+L_{X/Y})\!\!\mid$ gives a birational embedding of $Y$.
Let $\pi_{m} : Y_{m}\to Y$ be the resolution of $\mbox{Bs}\mid\!\! m!(K_{Y}+L_{X/Y})\!\!\mid$ such that for every $m > m_{0}$
\begin{equation}\label{pim}
\pi_{m} : Y_{m}\to Y
\end{equation}
factors through $\pi_{m-1} : Y_{m-1}\to Y$.
Let 
\begin{equation}\label{mum}
\mu_{m} : Y_{m} \to Y_{m-1}
\end{equation}
be the natural morphism. 
Here we may and do take $\mu_{m}: Y_{m}\to Y_{m-1}$ such that 
the exceptional divisor of $\pi_{m}$ is contained in  $\pi_{m}^{-1}(V)$. 
This is certainly possible by the definition of $V$.  
Hence we have an (possibly infinite) tower 
\begin{equation}\label{tower}
\cdots\stackrel{\mu_{m+2}}{\to} Y_{m+1}\stackrel{\mu_{m+1}}{\to}Y_{m}\stackrel{\mu_{m}}{\to} Y_{m-1}\stackrel{\mu_{m-1}}{\to}\cdots
\end{equation}
In the following proof, we shall consider this  tower.  
But by the recent result on finite generation of canonical ring
(\cite{b-c-h-m}) we may avoid to consider an infite tower.  Namely 
we just need to consider one sufficiently large $m_{0}$.  This certainly simplifies the proof.   
The reason why we do not use the finite generation of canonical ring  is that it is not essential 
from the analytic  point of view  and one may extend the theory to 
the case of  LC pairs (cf. \cite{LC}). 
Let 
\begin{equation}
\pi_{m}^{*}\mid\!m!(K_{Y}+L_{X/Y})\!\mid = \mid\!\!P_{m}\!\!\mid + F_{m}
\end{equation}
be the decomposition of $\pi_{m}^{*}\mid\! m!(K_{Y}+L_{X/Y})\!\mid$ into the 
free part $\mid\!\!P_{m}\!\!\mid$ and the fixed component $F_{m}$. 
Let $V$ be the analytic subset of $Y$ defined by:  
\begin{eqnarray*}
V:= &\hspace{-30mm}\{ y\in Y^{\circ}\mid y\in \cap_{m>0}\mbox{Bs}|m!(K_{Y}+L_{X/Y})|\,\,\mbox{or} \\
&\mbox{$\Phi_{\mid m!(K_{Y}+L_{X/Y})\mid}$ is not an embedding 
around $y$ for  $m >> m_{0}$}\} \\
\hspace{-90mm} & \cup \{\mbox{the discriminant locus of $f$}\}.
\end{eqnarray*}
By taking a suitable modification of $Y$, we may and do assume that 
$V$ is a divisor with normal crossings. 

There exists an effective $\mathbb{Q}$-divisor 
$E_{m}$ on $Y_{m}$ respectively  such that the followings hold for every $m\geqq m_{0}$.  
\begin{enumerate}
\item[(1)] $P_{m}- E_{m}$ is ample on $Y_{m}$. 
\item[(2)] All the coefficients of $E_{m}$ are less than $1$, i.e., $\lfloor E_{m}\rfloor = 0$.  
\item[(3)] $\mbox{Supp}\,E_{m} = \pi_{m}^{-1}(V)$.  
\item[(4)] $((m+1)!)^{-1}(P_{m+1}- E_{m+1})-\mu_{m+1}^{*}(m!)^{-1}(P_{m}- E_{m}))$
is effective.  
\end{enumerate}
The existence of such $\{E_{m}\}$ follows from the definition of  $V$ and the trivial fact that 
for any composition of successive blowing ups 
\[
\varpi : \tilde{\mathbb{P}^{\nu}}\to \mathbb{P}^{\nu}
\]
of a projective space $\mathbb{P}^{\nu}$ with smooth centers, 
there exists an effective $\mathbb{Q}$-divisor $B$  supported on the exceptional  divisors of $\varpi$ such that $\varpi^{*}\mathcal{O}(1) - B$ 
is ample.  Hence by taking a suitable successive blowing ups 
with smooth centres over $V$, if necessary,  
we may assume the existence of such effective $\mathbb{Q}$-divisors 
$\{ E_{m}\}$ by considering the image of $Y_{m}$ by the morphism 
associated with $|P_{m}|$ from $Y_{m}$ into a projective space. 
 
After taking such a sequence $\{ E_{m}\}$, we replace $\{ E_{m}\}$ by 
$\{ 2^{-m}E_{m}\}$.   
Then it has the same properties as above. 
And we shall denote $\{ 2^{-m}E_{m}\}$ again by $\{ E_{m}\}$. 
Then insead of (4) we have 
\begin{equation}\label{3prime} 
\mbox{(4)}^{\prime} \,\,\,\,((m+1)!)^{-1}(P_{m+1}- E_{m+1})-\mu_{m+1}^{*}(m!)^{-1}(P_{m}- E_{m}))
\,\,\mbox{is effective }
\end{equation}
\vspace{-5mm}
\[
\mbox{and contains} \,\,\varepsilon_{m}(\pi_{m+1}^{-1}V)_{red}\,\, \mbox{for some positive number $\varepsilon_{m}$.} 
\]

\noindent Let $h_{(m)}$ be a $C^{\infty}$-hermitian metric on 
$\pi_{m}^{*}((m!)^{-1}(P_{m}- E_{m}))$ with strictly positive 
curvature.  We note that by (\ref{3prime}) $\{ h_{(m)}\}$ is getting less singular as $m$ tends to infinity and $h_{(m+1)}$ is strictly less singular 
than $h_{(m)}$ along $V$ (if we consider the metrics as 
singular hermitian metrics on $K_{Y} + L_{X/Y}$).   Then 
\begin{equation}\label{Omegadelta}
\Omega_{m,\delta} :=  h_{(m)}^{-1}\cdot (\pi_{m}^{*}h_{L_{X/Y},\delta})
\end{equation}  
is considered as a degenerate volume form on $Y_{m}$, where 
$\{ h_{L_{X/Y},\delta}\}$ is the smoothing of the Hodge metric $h_{L_{X/Y}}$ 
as in Section \ref{smoothing}. 
We note that $\Omega_{m,\delta}^{-1}\cdot (\pi_{m}^{*}h_{L_{X/Y},\delta})^{-1}
= h_{(m)}$ is 
 a metric with algebraic singularities on $K_{Y} + L_{X/Y}$.  

Now we shall consider the equation :
\begin{equation}\label{ml}
(-\mbox{Ric}\,\,\Omega_{m,\delta}+ \sqrt{-1}\,\pi_{m}^{*}\Theta_{h_{L_{X/Y},\delta}} + \sqrt{-1}\partial\bar{\partial}u_{m,\delta})^{n}
= \Omega_{m,\delta}\cdot e^{u_{m,\delta}}. 
\end{equation}
on $Y_{m}$.
Then by the definition (\ref{Omegadelta}) of $\Omega_{m,\delta}$, (\ref{ml}) is  equivalent to 
\begin{equation}
(\sqrt{-1}\Theta_{h_{(m)}} + \sqrt{-1}\partial\bar{\partial}u_{m,\delta})^{n}
= \Omega_{m,\delta}\cdot e^{u_{m,\delta}}. 
\end{equation}

Then since $\sqrt{-1}\Theta_{h_{(m)}}$ is a $C^{\infty}$-K\"{a}hler form 
on $Y_{m}$ by \cite[p.387, Theorem 6]{y}, solving (\ref{ml}),  
we see that there exists a $u_{m,\delta}\in C^{\infty}(Y_{m}\backslash\pi_{m}^{-1}V)$ and the closed positive current: 
\begin{equation}
\omega_{m,\delta} := -\mbox{Ric}\,\Omega_{m,\delta} + \sqrt{-1}\pi_{m}^{*}\Theta_{h_{L_{X/Y},\delta}} +\sqrt{-1}\partial\bar{\partial}u_{m,\delta}
\end{equation}
on $Y_{m}$ such that 
\begin{enumerate}
\item[(1)] $-\mbox{Ric}_{\omega_{m,\delta}} + \sqrt{-1}\,\pi_{m}^{*}\Theta_{h_{L_{X/Y},\delta}}= \omega_{m,\delta}$  holds on $Y_{m}- \mbox{Supp}\, E_{m}$,
\item[(2)] The absolutely continuous part $\omega_{m,\delta,abc}$ of $\omega_{m,\delta}$ is closed and  represents
$2\pi (m!)^{-1}(P_{m}-E_{m})$, 
\item[(3)] $(\pi_{m})_{*}\omega_{m,\delta}$ represents the class $2\pi c_{1}(K_{Y}+L_{X/Y})$,
\item[(4)] There exists a positive constant
$C(m,\delta)$  such that  
\begin{equation}\label{almostbound}
|u_{m,\delta}| \leqq C(m,\delta)
\end{equation}
holds on $Y_{m}$. 
\end{enumerate}
Here we note that 
\begin{equation}
-\mbox{Ric}\,\,\Omega_{m,\delta}+ \sqrt{-1}\,\pi_{m}^{*}\Theta_{h_{L_{X/Y},\delta}}
= \sqrt{-1}\,\Theta_{h_{(m)}}
\end{equation}
holds, hence it is independent of $\delta$. 
If we set 
\begin{equation}\label{omegam}
\omega_{(m)} : = \sqrt{-1}\Theta_{h_{(m)}},
\end{equation}
then the equation is transcripted as :
\begin{equation}\label{ml2} 
\log \frac{\omega_{m,\delta}^{n}}{\Omega_{m,\delta}} = 
\log \frac{(\omega_{(m)} +\sqrt{-1}\partial\bar{\partial}u_{m,\delta})^{n}}
{\Omega_{m,\delta}} = u_{m,\delta}. 
\end{equation}
Let us consider $\{ \omega_{m,\delta}^{n}\}$ as a sequence of volume forms 
on $Y\backslash V$\footnote{Since we consider $\omega_{m,\delta}$ 
as a current, it seems to be more authentic to denote 
$\omega_{m,\delta,abc}^{n}$ instead of $\omega_{m,\delta}$. But we consider the 
eqation (\ref{ml2}) on $Y \backslash V$. }.  And we shall identify $(\pi_{m})_{*}\omega_{m,\delta}$ with  
$\omega_{m,\delta}$ on $Y \backslash V$.   Hereafter we shall identify $Y \backslash  V$ with 
a Zariski open subset of $Y_{m}$ for every $m$ and consider everything 
on $Y \backslash V$ (if without fear of confusion). 
Then by the maximum principle we have the following monotonicity lemma.

\begin{lemma}\label{monotonicity}(Monotonicity Lemma)
\begin{equation}
\omega_{m,\delta}^{n} \leqq \omega_{m+1,\delta}^{n}
\end{equation}
holds on $Y \backslash V$. \fbox{}
\end{lemma}
{\em  Proof of Lemma \ref{monotonicity}}.  
We note that by the construction the followings hold. 
\begin{enumerate}
\item[(1)]  The absolutely continuous parts of $\omega_{m,\delta}$ and $\omega_{m+1,\delta}$ represent
$2\pi (m!)^{-1}(P_{m}-E_{m})$ and $2\pi ((m+1)!)^{-1}(P_{m+1}-E_{m+1})$ respectively.
\item[(2)]  $\mu_{m+1}^{*}((m!)^{-1}(F_{m}+ E_{m}))-((m+1)!)^{-1}(F_{m+1}+E_{m+1})$ is effective and contains $\varepsilon_{m}(\pi_{m+1}^{-1}V)$ for some positive 
number $\varepsilon_{m}$.
\end{enumerate}
We note that by the boundedness of $u_{m,\delta}$ and $u_{m+1,\delta}$ (cf. (\ref{almostbound})) and the equation (\ref{ml2}), we see that 
the asymptotics ot $\omega_{m,\delta}^{n}$ and $\omega_{m+1,\delta}^{n}$
 near $V$ is the same as $\Omega_{m,\delta}$ and $\Omega_{m+1,\delta}$ 
 respectively. Since by the condition (2) above, $\Omega_{m,\delta}/\Omega_{m+1,\delta}$ tends to $0$ toward $V$, the function $\phi_{m,\delta}$  defined by 
\begin{equation}
\phi_{m,\delta}:= \log \frac{\omega_{m,\delta}^{n}}{\omega_{m+1,\delta}^{n}}
\end{equation} 
tends to $-\infty$ toward $V$.   
Hence there exists a point $p_{m}\in Y \backslash V$, where 
$\phi_{m,\delta}$ takes its maximum.  Then  
\begin{equation}\label{atpm}
\sqrt{-1}\partial\bar{\partial}\,\phi_{m,\delta}(p_{m}) \leqq 0
\end{equation}
holds.   By the equation
\begin{equation}
-\mbox{Ric}_{\omega_{k,\delta}} + \sqrt{-1}\,\pi_{k}^{*}\Theta_{h_{L_{X/Y},\delta}}= \omega_{k,\delta},  
(k= m, m+1) 
\end{equation}
(\ref{atpm}) implies that 
\begin{equation}
\omega_{m,\delta}(p_{m}) \leqq \omega_{m+1,\delta}(p_{m})
\end{equation}
holds.  In particulat $\phi_{m,\delta}(p_{m}) \leqq 0$ holds.  
Hence by the definition of$p_{m}$. this implies that $\phi_{m,\delta} \leqq 0$ holds on $Y \backslash V$.  Hence $\omega_{m,\delta}^{n} \leqq \omega_{m+1,\delta}^{n}$ holds
on $Y \backslash V$.  This completes the proof of Lemma \ref{monotonicity}. \fbox{}  
\vspace{5mm}\\ 
Now we shall consider the uniform $C^{0}$-estimate on every compact 
subset of $Y \backslash V$.
Let us fix a positive integer $s \geqq m_{0}$ and let 
\[
\omega_{(s)} = \sqrt{-1}\,\Theta_{h_{(s)}}
\]
be the K\"{a}hler form defined as (\ref{omegam}). 
We shall use $h_{(s)}$ and $\omega_{(s)}$ as standards in the following 
estimate.   
For every $m > s$,  let $v_{m,\delta}$ be a $C^{\infty}$-function on $Y \backslash D$ 
defined by 
\begin{equation}\label{defvmd}
v_{m,\delta} := u_{m,\delta} + \log \frac{h_{(s)}}{h_{(m)}}. 
\end{equation}
Then 
\begin{equation}\label{v-eq0}
\omega_{m,\delta} = \omega_{(s)} + \sqrt{-1}\partial\bar{\partial}v_{m,\delta}
\end{equation}
and 
\begin{equation}\label{v-eq}
\log \frac{(\omega_{(s)} + \sqrt{-1}\partial\bar{\partial}v_{m,\delta})^{n}}{\Omega_{s,\delta}} = v_{m,\delta}
\end{equation}
hold. 
By the condition (\ref{3prime}), we see that for every  $m > s$, 
$\log (h_{(s)}/h_{(m)})$ tends to $+\infty$ toward $V$.  
And  by the boundedness of $u_{m,\delta}$ (cf. (\ref{almostbound})), we have
the estimate: 
\begin{equation}\label{a-bound}
 - C(m,\delta) +\log \frac{h_{(s)}}{h_{(m)}} \leqq v_{m,\delta} \leqq  C(m,\delta) + \log \frac{h_{(s)}}{h_{(m)}},
\end{equation}
 where $C(m,\delta)$ is the positive constant 
as in (\ref{almostbound}). 
Hence  $v_{m,\delta}$ tends to $+\infty$ toward  $V$. 
This implies that there exists a point $p_{0}\in Y \backslash V$, where $v_{m,\delta}$ takes its minimum. 
Now we note that  by (\ref{v-eq}), 
\begin{equation}
\log \frac{(\omega_{(s)} + \sqrt{-1}\partial\bar{\partial}v_{m,\delta})^{n}}{
\omega_{(s)}^{n}} = 
\int_{0}^{1}\Delta_{(s,m,\delta,t)}v_{m,\delta}\,\,dt = v_{m,\delta} - \log \frac{\omega_{(s)}^{n}}{\Omega_{s,\delta}} 
\end{equation}
hold, where $\Delta_{(s,m,\delta,t)}$ denotes the trace of $\sqrt{-1}\partial\bar{\partial}v_{m,\delta}$ with respect to the K\"{a}hler form:
$(1-t)\omega_{m,\delta} + t\omega_{(s)}$.
Hence by the minimum principle, 
\begin{equation}\label{minineq}
v_{m,\delta}(p_{0}) \geqq   \log \frac{\omega_{(s)}^{n}}{\Omega_{s,\delta}}(p_{0})
\end{equation}
holds. 

On the other hand, by (\ref{Omegadelta}) and the  definition of $h_{L_{X/Y},\delta}$(cf. (\ref{defsmoothing})), $\Omega_{s,\delta}$ tends to $0$ toward $V$. 
Hence there exists a positive constant $C_{-}(s)$  independent of $\delta$ 
such that  
\begin{equation}\label{lowermint}
\min_{y\in Y} \log \frac{\omega_{(s)}^{n}}{\Omega_{s,\delta}}(y)
\geqq C_{-}(s) 
\end{equation}
holds.  By  (\ref{lowermint}) and  (\ref{minineq}), we see that  
\begin{equation}
v_{m,\delta}(y) \geqq C_{-}(s)
\end{equation}
holds for every $y\in Y \backslash V$. 
By the definition of $v_{m,\delta}$ (cf. (\ref{defvmd})), we have that 
\begin{equation}\label{lowerbound}
u_{m,\delta}  \geqq \log \frac{h_{(m)}}{h_{(s)}}
+ C_{-}(s)
\end{equation}
holds. 
Replacing $s$ by $t > s$,  for $m > t >  s$ we have
\[
u_{m,\delta}  \geqq \log \frac{h_{(m)}}{h_{(t)}} + C_{-}(t)
\] 
and hence by (\ref{defvmd}), we obtain the following lemma.
\begin{lemma}\label{lowert}
There exists a positive constant $C_{-}(t)$ depending only on $t > s$ such that
 for every $m > t$
\begin{equation}\label{loweq}
v_{m,\delta} \geqq \log \frac{h_{(s)}}{h_{(t)}} + C_{-}(t)
\end{equation}
holds. In particular $v_{m,\delta}$ tends to infinity toward $V$.  \fbox{}
\end{lemma}

On the other hand, we obtain the upper estimate of  $u_{m,\delta}$ as 
follows.
We may and do assume that $V$ is a divisor with normal crossings.  
Let $\Omega_{P}$ be a volume form on $Y \backslash V$ with Poincar\'{e} growth, i.e.,
for every polydisk $\Delta^{n}$ in $Y$ such that 
\begin{equation}
\Delta^{n} \cap V = \{ (z_{1},\cdots,z_{n})\in \Delta^{n}|
z_{1}\cdots z_{k} = 0\},
\end{equation}
\begin{equation}
\Omega_{P} = c \frac{|dz_{1}\wedge\cdots\wedge dz_{n}|^{2}}{\prod_{i=1}^{k}
|z_{i}|(\log |z_{i}|^{2})^{2}}, 
\end{equation} 
where $c$ is a positive $C^{\infty}$-function on $\Delta^{n}$. 
Such a $\Omega_{P}$ can be constructed easily by using a partition of unity. 
We set 
\begin{equation}\label{tilde}
\tilde{u}_{m,\delta}:= u_{m,\delta} + \log \frac{\Omega_{m,\delta}}{\Omega_{P}}.\end{equation}
By the condition (\ref{3prime}) and the boundedness of 
$u_{m,\delta}$,  there exists a point $p_{0}^{\prime}$ on $Y \backslash V$
such that  $\tilde{u}_{m,\delta}$ takes
its  maximum at $p_{0}^{\prime}$.
Then since
\begin{equation}
\log \frac{\omega_{m,\delta}^{n}}{\Omega_{P}} = u_{m,\delta} + \log \frac{\Omega_{m,\delta}}{\Omega_{P}} = \tilde{u}_{m,\delta}
\end{equation}
hold, at $p_{0}^{\prime}$ we have that 
\begin{equation}
\sqrt{-1}\partial\bar{\partial}\log \frac{\omega_{m,\delta}^{n}}{\Omega_{P}}(p_{0}^{\prime}) \leqq 0
\end{equation}
holds.
Hence we have the inequality: 
\begin{equation}
-\mbox{Ric}_{\omega_{m,\delta}} \leqq  (-\mbox{Ric}\,\Omega_{P} )(p_{0}^{\prime}).
\end{equation} 
By the equation:
\begin{equation}
-\mbox{Ric}_{\omega_{m,\delta}} + \sqrt{-1}\,\Theta_{h_{L_{X/Y},\delta}}
= \omega_{m,\delta},
\end{equation}
we see that 
\begin{equation}
\omega_{m,\delta}(p_{0}^{\prime}) \leqq -\mbox{Ric}\,\Omega_{P} +\sqrt{-1}\Theta_{h_{L_{X/Y},\delta}}
\end{equation}
holds.  In particular, 
\begin{equation}
\omega_{m,\delta}^{n}(p_{0}^{\prime}) \leqq (-\mbox{Ric}\,\Omega_{P} +\sqrt{-1}\,\Theta_{h_{L_{X/Y},\delta}})^{n}(p_{0}^{\prime})
\end{equation}
holds.   This implies that 
\begin{equation}\label{maxu}
\tilde{u}_{m,\delta} =  u_{m,\delta} + \log \frac{\Omega_{m,\delta}}{\Omega_{P}}
\leqq \log \frac{(-\mbox{Ric}\,\Omega_{P} + \sqrt{-1}\,\Theta_{h_{L_{X/Y},\delta}})^{n}}{\Omega_{P}}(p_{0}^{\prime})
\end{equation}
holds on $Y$. 
By the construction of $h_{L_{X/Y},\delta}$, Lemma \ref{Poincare} and Remark \ref{ric}, there exists a positive constant 
$C_{+}$ independent of $\delta$ such that  
\begin{equation}\label{ricp}
\frac{(-\mbox{Ric}\,\Omega_{P} + \sqrt{-1}\,\Theta_{h_{L_{X/Y},\delta}})^{n}}{\Omega_{P}} \leqq \exp(C_{+})
\end{equation}
holds on $Y \backslash V$.
Combining (\ref{maxu}) and (\ref{ricp}) we have that 
\begin{equation}\label{upperbound}
u_{m,\delta} \leqq C_{+} - \log \frac{\Omega_{m,\delta}}{\Omega_{P}}
\end{equation}
and 
\begin{equation}\label{uniform}
e^{u_{m,\delta}}\Omega_{m,\delta} \leqq \exp(C_{+})\cdot\Omega_{P}
\end{equation} 
hold on $Y$.

On the other hand by the definition (\ref{Omegadelta}) and the definition 
of the smoothing $\{ h_{L_{X/Y},\delta}\}$ (\ref{defsmoothing}), there 
exists a sequence of positive number $\{ \epsilon (\delta)\}$ such that 
\begin{enumerate}
\item[(1)] $\lim_{\delta\rightarrow 0}\epsilon(\delta) = 0$ holds,  
\item[(2)] For every $0< \lambda < \delta$, the inequality:  
\begin{equation}
\Omega_{m,\delta} \leqq (1+\epsilon(\delta))\Omega_{m,\lambda}
\end{equation}
holds. 
\end{enumerate}
Then for $0 < \lambda < \delta$ by  (\ref{ml2})
\begin{equation}\label{log}
\log \frac{(\omega_{(m)} +\sqrt{-1}\partial\bar{\partial}u_{m,\lambda})^{n}}
{(\omega_{(m)} +\sqrt{-1}\partial\bar{\partial}u_{m,\delta})^{n}} 
= \log \frac{\Omega_{m,\lambda}}{\Omega_{m,\delta}} + (u_{m,\lambda}- u_{m,\delta})
\end{equation}
holds. 
We note that 
\begin{equation}\label{laplacian}
\log \frac{(\omega_{(m)} +\sqrt{-1}\partial\bar{\partial}u_{m,\lambda})^{n}}
{(\omega_{(m)} +\sqrt{-1}\partial\bar{\partial}u_{m,\delta})^{n}} 
= \int_{0}^{1}\tilde{\Delta}_{t}(u_{m,\lambda}-u_{m,\delta})\,dt,  
\end{equation}
where $\tilde{\Delta}_{t}$ denotes the Laplacian with respect to 
$(1 - t)\tilde{\omega}_{m,\lambda} + t\tilde{\omega}_{m,\delta}$. 
Then by (\ref{log}),(\ref{laplacian}) and the maximum principle, we see that 
\begin{equation}\label{lambda}
u_{m,\lambda}- u_{m,\delta} \leqq -\min_{Y}\log \frac{\Omega_{m,\lambda}}{\Omega_{m,\delta}} < \log (1 + \epsilon (\delta))
\end{equation}
holds on $Y$.  
This argument is not quite right, since $u_{m,\lambda},u_{m,\delta}$ are not 
$C^{2}$ on $Y$ (although they are $C^{2}$ bounded). 

To justify the argument we proceed as in the proof of  [p.387,Theorem 6]\cite{y}, i.e., we shall consider the perturbation of the equation (\ref{ml2}).
By the construction of $\Omega_{m,\delta}$ (cf. (\ref{Omegadelta})) 
the $0$-locus of  $\Omega_{m,\delta}$ is the divisor $(\pi_{m}^{*}V)_{red}$ with normal crossings on  $Y_{m}$.  Let $(\pi_{m}^{*}V)_{red} = \sum V^{(m)}_{j}$ be the irreducible decomposition.
 Let us write $\Omega_{m,\delta}$ as 
\begin{equation}
\Omega_{m,\delta} = \left(\prod_{j} \parallel\tau_{j}\parallel^{2a_{j}}\right)
\cdot \tilde{\Omega}_{m,\delta},
\end{equation}  
where $\tilde{\Omega}_{m,\delta}$ is a nondegenerate $C^{\infty}$-volume form on $Y_{m}$
and $\parallel\tau_{j}\parallel$ denotes the hermitian norm of a  
holomorphic section $\tau_{j}$ of the line bundle $\mathcal{O}_{Y_{m}}(V^{(m)}_{j})$
with divisor $V^{(m)}_{j}$ on  $Y_{m}$ with respect to a fixed $C^{\infty}$-hermitian metric and $\{a_{j}\}$ are  positive rational numbers.
We may and do take the factor $\left(\prod_{j} \parallel\tau_{j}\parallel^{2a_{j}}\right)$ independent of $\delta$.   
Now for $0< \varepsilon << 1$, we shall consider the perturbed equation: 
\begin{equation}\label{peq}
\log \frac{(\omega_{(m)} +\sqrt{-1}\partial\bar{\partial}u_{m,\delta}(\varepsilon))^{n}}{\Omega_{m,\delta}(\varepsilon)} = u_{m,\delta}(\varepsilon), 
\end{equation} 
where $\Omega_{m,\delta}(\varepsilon)$ is a $C^{\infty}$ nondegenerate 
volume form on $Y$ 
defined by 
\begin{equation}
\Omega_{m,\delta}(\varepsilon) = \left(\prod_{j}(\parallel\tau_{j}\parallel^{2} + \varepsilon)^{a_{j}}\right)\cdot\tilde{\Omega}_{m,\delta}.
\end{equation}
Then (\ref{peq}) has a unique $C^{\infty}$-solution 
$u_{m,\delta}(\varepsilon)$ and 
\begin{equation}
u_{m,\delta} = \lim_{\varepsilon\downarrow 0}u_{m,\delta}(\varepsilon)
\end{equation}
holds as in \cite[p.387,Theorem 6]{y} in the $C^{2}$-norm with respect to the 
K\"{a}hler form $\omega_{(m)}$ on $Y_{m}$.  
Then replacing $u_{m,\lambda},u_{m,\delta}$ by $u_{m,\lambda}(\varepsilon),u_{m,\delta}(\varepsilon)$ repsectively and letting $\varepsilon\downarrow 0$, 
we may justify (\ref{lambda}). 
Since $\{ u_{m,\delta}\}$ is almost monotone decreasing as $\delta\downarrow 0$ as (\ref{lambda}), we have that 
\begin{equation}
dV_{m}:= \frac{1}{n!}\lim_{\delta\downarrow 0} e^{u_{m,\delta}}\Omega_{m,\delta} 
\end{equation}
exists and  by (\ref{ml2}), we have 
\begin{equation}\label{dvm}
dV_{m} = \frac{1}{n!}\lim_{\delta\downarrow 0}\omega_{m,\delta,abc}^{n}
\end{equation}
holds on $Y$.
We note that 
\[
\int_{Y}\omega_{m,\delta,abc}^{n} = (m!)^{-n}(P_{m}- E_{m})^{n}
\]
holds. 
By (\ref{uniform}) and the Lebesgue's bounded convergence theorem, we  see 
that 
\begin{equation}\label{growthm}
\frac{1}{(2\pi)^{n}}\int_{Y}dV_{m} = \frac{1}{n!(m!)^{n}}(P_{m} - E_{m})^{n}
\end{equation}
holds.   
We set 
\begin{equation}
\omega_{m} := \lim_{\delta\downarrow 0}\omega_{m,\delta}  
\end{equation} 
in the sense of current and 
\begin{equation}
h_{m} : = dV_{m}^{-1}\cdot h_{L_{X/Y}}.
\end{equation}
Then $h_{m}$ is a singular hermitian metric on $K_{Y}+ L_{X/Y}$ 
with semipositive curvature in the sense of current by (\ref{ml}) and (\ref{ml2}). 
But also we may consider $h_{m}$ as a singular hermitian metric $\tilde{h}_{m}$  on $(m!)^{-1}(P_{m} - E_{m})$ by the natural inclusion :
\[
\mathcal{O}_{Y_{m}}(\ell !\,\pi_{m}^{*}(m!(P_{m} - E_{m})))
\hookrightarrow \mathcal{O}_{Y_{m}}(\ell !\,\pi_{m}^{*}(m!(K_{Y} + L_{X/Y}))),
\]
where $\ell$ is a sufficiently large positive integer such that 
$\ell !\,\pi_{m}^{*}(m!(P_{m} - E_{m}))$ is Cartier. 
Now we introduce the following notion. 
\begin{definition}\label{singvolume}
Let $M$ be a projective manifold of dimension $n$  and let $(L,h_{L})$ be a 
pseudoeffective singular hermitian line bundle on $X$.  
We define the number $\mu (L,h_{L})$ by 
\[
\mu (L,h_{L}):= n!\limsup_{m\rightarrow\infty}m^{-n}h^{0}(M,{\cal O}_{M}(mL)
\otimes {\cal I}(h_{L}^{m}))
\]
is called the volume of $(L,h_{L})$. 
 \fbox{}
\end{definition}
\begin{remark}
This definition is easily generalized to the case of 
singular hermitian $\mathbb{Q}$-line bundles. \fbox{}
\end{remark}
\noindent Suppose that  the Lelong number $\nu (\sqrt{-1}\Theta_{\tilde{h}_{m}})$ of $\sqrt{-1}\Theta_{\tilde{h}_{m}}$ satisfies the inequality: 
$\nu (\sqrt{-1}\Theta_{\tilde{h}_{m}},y_{m}) > c$ 
for some  $y_{m}\in Y_{m}$ and a positive number $c$, then 
by the basic property of the Lelong number, we see that 
\begin{equation}\label{maxideal}
\mathcal{I}(\tilde{h}_{m}^{\ell!})_{y_{m}} \subseteq \mathfrak{m}_{y_{m}}^{\lfloor c\ell!\rfloor}
\end{equation}
holds for every sufficiently large $\ell$, where 
$\mathfrak{m}_{y_{m}}$ denotes the maximal ideal at $y_{m}$. 
Hence the strict inequality 
\begin{equation}\label{volineq}
\mu ((m!)^{-1}(P_{m}- E_{m}),\tilde{h}_{m})< (m!)^{-n}(P_{m}- E_{m})^{n} 
\end{equation}
holds, where  $\mu ((m!)^{-1}(P_{m}- E_{m}),\tilde{h}_{m})$ denotes 
the volume of $((m!)^{-1}(P_{m}- E_{m}),\tilde{h}_{m})$ (cf. Definition \ref{volume} below). 
On the other hand 
\begin{equation}\label{bouck}
\frac{1}{(2\pi)^{n}}\int_{Y}\omega_{m,abc}^{n} 
 = \mu ((m!)^{-1}(P_{m}- E_{m}),\tilde{h}_{m})
\end{equation}
holds by  (\ref{growthm}) and a theorem of Boucksom (\cite[Proposition 3.1]{bouk}).
By the equality (\ref{growthm}), (\ref{bouck}) contradicts (\ref{volineq}).   
Hence $\nu (\sqrt{-1}\Theta_{\tilde{h}_{m}})$ is identically $0$ and 
hence $\tilde{h}_{m}$ is an AZD of $(m!)^{-1}(P_{m} - E_{m})$. 
By Lemma \ref{monotonicity}, we see that 
$\{ dV_{m}\}$ is monotone increasing in $m$ on $Y$, hence 
$\{ h_{m}\}$ is getting less singular as $m$ tends to infinity 
as metrics on $K_{Y} + L_{X/Y}$.    
Then by (\ref{uniform}) and Lebesgue's bounded convergence theorem, we have the following lemma. 
\begin{lemma}\label{convergencevol}
\begin{equation}\label{dvy}
dV_{Y}:= \frac{1}{n!}\lim_{m\rightarrow\infty} dV_{m}
\end{equation}
exists as a degenerate volume form on $Y$.  
And if we define the singular hermitian metric $h_{K}$ on
on $K_{Y} + L_{X/Y}$ by 
\begin{equation} 
h_{K}:= dV_{Y}^{-1}\cdot h_{L_{X/Y}}, 
\end{equation}
then $h_{K}$ is an AZD of $K_{Y} + L_{X/Y}$. 
\fbox{}
\end{lemma}
By the construction of $\omega_{m}$,  
\begin{equation}\label{meq}
-\mbox{Ric}_{\omega_{m}} + \sqrt{-1}\pi_{m}^{*}\Theta_{h_{L_{X/Y}}} = \omega_{m}
\end{equation}
holds for every $m \geqq 1$.
Then by (\ref{meq}) and Lemma \ref{convergencevol}, if we set  
\begin{equation}
\omega_{Y}:= \sqrt{-1}\,\Theta_{h_{K}},
\end{equation}
then $\omega_{Y}:= \lim_{m\rightarrow\infty}\omega_{m}$
holds and $\omega_{K}$ is a closed positive current on $Y$. 
Moreover  
\begin{equation}
\omega_{Y}^{n} = \lim_{m\rightarrow\infty}\omega_{m}^{n}
\end{equation} 
holds on $Y \backslash V$ by the definition.  Hence 
$\omega_{Y}$ satisfies  equation:
\begin{equation}
-\mbox{Ric}_{\omega_{Y}} + \sqrt{-1}\Theta_{h_{L_{X/Y}}} = \omega_{Y}.
\end{equation}
This completes the proof of Theorem \ref{current}. \fbox{}

\subsection{Regularity of the  canonical K\"{a}hler current}

Here we shall prove Theorem \ref{main}, by using the recent result on 
the finite generation of canonical ring (\cite{b-c-h-m}). 
By Theorem \ref{current}, we only need to prove the $C^{\infty}$-regularity of $h_{K}$ 
and $\omega_{Y}$ on a nonempty Zariski open subset of $Y$.  

The proof here is more or less parallel to the existence of the singular K\"{a}hler-Einstein metrics in \cite{s,KE} and is based on \cite{y} and the idea in \cite{tu0}.  But since the Hodge metric $h_{L_{X/Y}}$ is  not of algebraic singularities, we need to consider the smoothing of the Hodge metric. 
This is the major difference.   We continue to use the notations  in Section 2.3
Let us start the proof of Theorem \ref{main}.   
By \cite{b-c-h-m}, we see that 
the canonical ring $R(X,K_{X})$ (cf.(\ref{RM})) is finitely generated.
Then by the definition of $L_{X/Y}$, we see that $R(Y,K_{Y}+L_{X/Y})$ 
is finitely generated also.  
Hence this implies that the tower (\ref{tower}) above can be taken 
to be finite.   Here we shall assume that (\ref{tower}) is a finite 
tower. Moreover taking $m_{0}$ sufficiently large, we may assume that $\mu_{m} : Y_{m}\to Y_{m-1}$ (cf.(\ref{mum}))  
is identity for every $m \geqq m_{0}$ and 
\begin{equation}\label{P}
P = \frac{1}{m!}P_{m}
\end{equation}
is independent of $m$. 
Moreover we may and do assume that 
$Y_{m} = Y$ holds for all $m$.    In this case, only $E_{m}$ 
varies and  we may assume  that $E_{m}$ is of the form 
\begin{equation}\label{E}
E_{m} = \frac{1}{2^{m}}E,
\end{equation}
where $E$ is a fixed effective $\mathbb{Q}$-divisor supported on $V$ such that 
$P - E$ is ample. 

Next we shall fix a $C^{\infty}$-hermitian metric $h_{(m)}$ with strictly 
positive curvature on $(m!)^{-1}(P_{m} - E_{m})$ as in Section 2.3. 
Let $h_{P}$ be a $C^{\infty}$-hermitian metric 
defined by the pull back of the Fubini-Study metric on 
the hyperplane bundle on $\mathbb{P}^{N}$ via the morphism 
$\Phi_{|P_{m_{0}}|} : Y \to \mathbb{P}^{N}$ and let   $h_{0}$ be a $C^{\infty}$ hermitian metric on 
$P - E$ with strictly positive curvature on $Y$.   
And we shall take $h_{(m)}$ in the previous section as 
\begin{equation}\label{h(m)}
h_{(m)}:= h_{P}^{1-\frac{1}{2^{m}m!}}\cdot h_{0}^{\frac{1}{2^{m}m!}}. 
\end{equation}

Let us fix $s \geqq m_{0}$ as in the Section 2.3. 

 By (\ref{lowerbound}) and (\ref{upperbound}), we see that
there exists a positive constant $C^{\prime}_{0}$ independent of $m$ and $\delta$
such that  
\begin{equation}\label{mean}
\int_{Y}|u_{m,\delta}|\omega_{(s)}^{n} < C^{\prime}_{0} 
\end{equation}
holds. 
Then since 
\begin{equation}
\omega_{(m)} + \sqrt{-1}\partial\bar{\partial}u_{m,\delta}
\end{equation} 
is a closed positive current on $Y$ and $\omega_{(m)}$ is a $C^{\infty}$-K\"{a}hler form on 
$Y$, $u_{m,\delta}$ is an almost plurisubharmonic funtion on $Y$.
By the sub-mean-value inequality for plurisubharmonic functions, we see that by (\ref{mean})  there exists 
a positive constant $C_{0}$ such that 
\begin{equation}\label{upper2}
u_{m,\delta} \leqq C_{0}
\end{equation}
holds on $Y$.  

By Lemma \ref{lowert}, (\ref{upper2}) and (\ref{defvmd}), we have the following lemma. 
\begin{lemma}\label{c0}
Let $s < t < m$ be as above. 
Then there exists a positive 
constant $C_{0}$ independent of $m$ and  $\delta$ such that
for every $\delta > 0$,
\begin{equation}\label{c03}
-C_{-}(t) + \log \frac{h_{(s)}}{h_{(t)}} \leqq v_{m,\delta}  \leqq C_{0} + \log \frac{h_{(s)}}{h_{(m)}},
\end{equation}
hold on $Y \backslash V$, where $C_{-}(t)$ is the constant as in (\ref{loweq}) in Lemma \ref{lowert}. 
\fbox{}
\end{lemma}

Now we shall estimate the $C^{2}$-norm of $\{ v_{m,\delta}\}_{\delta > 0}$
on every compact subset of $Y \backslash V$. 

\begin{lemma}\label{c2}(\cite[p. 127, Lemma 2.2]{tu0}))
We set 
\begin{equation}
f := \log \frac{\omega_{(s)}^{n}}{\Omega}.
\end{equation}
Let $C$ be a positive number such that 
\begin{equation}\label{bisec}
C + \inf_{\alpha\neq \beta}R_{\alpha\bar{\alpha}\beta\bar{\beta}} > 1
\end{equation}
holds on $Y$, where $R_{\alpha\bar{\alpha}\beta\bar{\beta}}$ denotes the 
bisectional curvature of $\omega_{(s)}$. 

Then 
\begin{equation}
e^{Cv_{m,\delta}}\Delta_{m,\delta}(e^{-Cv_{m,\delta}}(n + \Delta_{s}\,v_{m,\delta}))
\geqq (n + \Delta_{s}\,v_{m,\delta}) 
\end{equation}
\begin{equation*}
 +  \Delta_{s}\left(f+ \log\frac{h_{(m)}}{h_{(s)}}\right) -(n+n^{2}\inf_{\alpha\neq \beta}R_{\alpha\bar{\alpha}\beta\bar{\beta}} ) 
\end{equation*}
\begin{equation*}
- C\cdot n(n+\Delta_{s}\,v_{m,\delta}) + 
(n+\Delta_{s}\,v_{m,\delta})^{\frac{n}{n-1}}\cdot\left(\frac{h_{(m)}}{h_{(s)}}\right)^{-\frac{1}{n-1}}\cdot\exp\left(-\frac{1}{n-1}(v_{m,\delta} + f)\right)
\end{equation*}
holds, where $\Delta_{s}$ denotes the Laplacian with respect to $\omega_{(s)}$
\mbox{\em (}i.e., $\Delta_{s}= \mbox{\em trace}_{\omega_{(s)}}\sqrt{-1}\partial\bar{\partial}$\mbox{\em )} and 
$\Delta_{m,\delta}$ denotes the Laplacian with respect to $\omega_{m,\delta}$. 
\fbox{} \vspace{5mm}
\end{lemma}

\noindent We note that in Lemma \ref{c2}, $C > 0$ does not depend on $m$ and $\delta$,  
but $C$ depends on $s$ (more precisely $\inf_{\alpha\neq \beta}R_{\alpha\bar{\alpha}\beta\bar{\beta}}$).
\begin{lemma}\label{takemax}
For any choice of  $C > 0$ satisfying (\ref{bisec}), $m$ and $\delta$ , there exists a point $y_{0} \in Y \backslash V$ where $e^{-Cv_{m,\delta}}(n + \Delta v_{m,\delta})$ takes its maximum. \fbox{} 
\end{lemma}
{\em Proof}. 
By Lemma \ref{lowert} above (cf. (\ref{c03})) 
\[
v_{m,\delta} \geqq \log \frac{h_{(s)}}{h_{(t)}} + C_{-}(t)
\]
holds.   
We note that for every $t > s$, $h_{(s)}/h_{(t)}$ has pole of 
positive order along $V$ by the condition (\ref{3prime}).
Hence by the lower estimate Lemma \ref{lowert}, $e^{-Cv_{m,\delta}}$ 
tends to $0$ toward $V$.
More precisely 
\begin{equation}\label{decay}
e^{-Cv_{m,\delta}} \leqq \exp(-C_{-}(t))\left(\frac{h_{(t)}}{h_{(s)}}\right)^{C}
\end{equation}
holds on $Y \backslash V$. 
On the other hand, we have that $\sqrt{-1}\partial\bar{\partial}v_{m,\delta}$  
is bounded with respect to the K\"{a}hler form $\omega_{(m)}$ on $Y$ as in \cite[p.387,Theorem 6]{y}, hence also with respect to $\omega_{(s)}$, since both 
$\omega_{(m)}$ and $\omega_{(s)}$ are K\"{a}hler forms on $Y$ by the assumption. Hence for any  $C > 0$ satisfying (\ref{bisec}) and every  $m$ and $\delta$ 
\begin{equation}\label{bounded}
e^{-Cv_{m,\delta}}(n + \Delta_{s}v_{m,\delta}) = O(1)
\end{equation}
holds on $Y$ and 
\begin{equation}
\lim_{y\to V}\left(e^{-Cv_{m,\delta}}(n + \Delta_{s}v_{m,\delta})\right)(y)
= 0
\end{equation}
holds.  Hence for any $C > 0$ satisfying (\ref{bisec}) and  every $m$ and $\delta$, there exists a point $y_{0} \in Y \backslash V$ where $e^{-Cv_{m,\delta}}(n + \Delta v_{m,\delta})$ takes its maximum.
This completes the proof of Lemma \ref{takemax}
\fbox{} \vspace{3mm} \\
Here we have used the fact that the Monge-Amp\`{e}re equation (\ref{ml}) of 
$u_{m,\delta}$ has algebraic singularities\footnote{By the definition (\ref{Omegadelta}), $\Omega_{m,\delta}$ has algebraic singularities}.  
Then by Lemma \ref{c2}, we have the following lemma. 
\begin{lemma}\label{y0}
If we take $C > 0$ satisfying (\ref{bisec}) , then there exists a positive constant $C_{2}$ independent of $m$ and $\delta$
such that 
\begin{equation}
0\leqq e^{-Cv_{m,\delta}(y_{0})}(n + \Delta_{s}\,v_{m,\delta})(y_{0})\leqq C_{2}
\end{equation}
holds. \fbox{}
\end{lemma}
{\em Proof.}
By the maximal principle, we have 
\begin{equation}
 (n + \Delta_{s}\,v_{m,\delta}(y_{0})) 
 +  \Delta_{s}\left(f+ \log\frac{h_{(m)}}{h_{(s)}}\right)(y_{0}) -(n+n^{2}\inf_{\alpha\neq \beta}R_{\alpha\bar{\alpha}\beta\bar{\beta}} ) 
- C\cdot n(n+\Delta_{s}\,v_{m,\delta})(y_{0}) 
\end{equation}
\begin{equation*}
+ (n+\Delta_{s}\,v_{m,\delta}(y_{0}))^{\frac{n}{n-1}}\cdot\left(\frac{h_{(m)}}{h_{(s)}}\right)^{-\frac{1}{n-1}}\!\!\!\!\!\!\!\!\!\!\!\!\!\!(y_{0})\cdot\exp\left(-\frac{1}{n-1}v_{m,\delta} + f\right)(y_{0}) \leqq 0
\end{equation*}
holds.   
Then we see that there exists a positive constant $C_{3}$ independent of 
$m$ and $\delta$ such that  
\begin{equation}\label{maxi}
n+\Delta_{s}\,v_{m,\delta}(y_{0})
\leqq C_{3} \left(1 + \left|\Delta_{s}\left(f+ \log\frac{h_{(m)}}{h_{(s)}}\right)(y_{0})\right|\right)^{n}
\end{equation}
holds.  
Since $f$ is $C^{\infty}$ on $Y$ and  by the definition of $h_{(m)}$ (cf. 
(\ref{h(m)})), there exists a positive constant $C_{4}$ independent of $m$ 
such that  
\begin{equation}
\left|\Delta_{s}\log\frac{h_{(m)}}{h_{(s)}}\right|  < C_{4} 
\end{equation}
holds on $Y$.  Hence by (\ref{maxi}), we see that there exists a positive 
constant $C_{5}$ independent of $m$,$\delta$ and $y_{0}$ such that 
\begin{equation}\label{maxi2}
n+\Delta_{s}\,v_{m,\delta}(y_{0}) \leqq C_{5}
\end{equation} 
holds. 
Next we shall consider the factor $e^{-Cv_{m,\delta}(y_{0})}$. 
We note that for every $t > s$, $h_{(s)}/h_{(t)}$ has pole of 
positive order along $V$ by the condition (\ref{3prime}).
Hence by the lower estimate Lemma \ref{lowert}, $e^{-Cv_{m,\delta}}$ 
tends to $0$ toward $V$.    
Then by the $C^{0}$-estimate Lemma \ref{lowert}, we see that there exists a positive constant $C_{6}$ 
independent of $m$ and $\delta$ 
such that  
\begin{equation}\label{exp}
\, e^{-Cv_{m,\delta}}\left( 1+\left|\Delta_{s}\left(f+ \log\frac{h_{(m)}}{h_{(s)}}\right)\right|\right)^{n}\leqq C_{6}
\end{equation}
holds on $Y \backslash V$.  Combining (\ref{maxi}) and (\ref{exp}),
by the definition of $y_{0}$,  
we complete the proof of Lemma \ref{y0}. \fbox{} \vspace{3mm} \\  

Let us take $C > 0$ satisfying (\ref{bisec}) as in Lemma \ref{c2}. 
By Lemma \ref{y0}  and the defnition of $y_{0}$ 
\begin{equation}\label{c2y}
e^{-Cv_{m,\delta}}(n +\Delta_{s}v_{m,\delta})
\leqq e^{-Cv_{m,\delta}(y_{0})}(n + \Delta_{s}v_{m,\delta}(y_{0}))
\leqq C_{2}
\end{equation}
hold.   Hence by (\ref{c2y}) we have the inequality :
\begin{equation}\label{c2g}
0\leqq n + \Delta_{s}\,v_{m,\delta} \leqq  \exp (C\cdot v_{m,\delta})\cdot C_{2}.
\end{equation}
Estimating $\exp (C\cdot v_{m,\delta})$ from above by Lemma \ref{c0}, (\ref{c2g}) implies that there exists a positive constant $C_{7}$ independent of $m$ and $\delta$ such that 
\begin{equation}
n + \Delta_{s}\,v_{m,\delta} \leqq C_{7}(\frac{h_{(s)}}{h_{(m)}})^{C}
\end{equation}
holds on $Y \backslash V$.
Applying the general theory of fully nonlinear elliptic equations (\cite{tr}), 
to the equation (\ref{v-eq}), we get a uniform higher order estimate of $\{ v_{m,\delta}\}_{\delta > 0}$ on every compact subset of $Y \backslash V$. 
Hence there exists a sequence $\{\delta_{j}\}$ with  $\delta_{j}\downarrow 0$
as $j$ tends to infinity  
such that 
\begin{equation}
\omega_{m}: = \lim_{j\rightarrow\infty}\omega_{m,\delta_{j}}  
\end{equation}
exists in $C^{\infty}$-topology on every compact subset of $Y \backslash V$. 
By using the diagonal argument, we may take $\{ \delta_{j}\}$ 
independent of  $m$.
Then 
\begin{equation}
\omega_{m} := \omega_{(s)} + \sqrt{-1}\partial\bar{\partial}v_{m}
\end{equation}
satisfies the equation :
\begin{equation}
\log \frac{\omega_{m}^{n}}{\Omega_{s}} = v_{m}, 
\end{equation}
on $Y \backslash V$, where 
\begin{equation}
\Omega_{s}:= h_{(s)}^{-1}\cdot h_{L_{X/Y}}
\end{equation} 
and hence  
\begin{equation}\label{m-eq}
-\mbox{Ric}_{\omega_{m}} + \sqrt{-1}\,\Theta_{h_{L_{X/Y}}} = \omega_{m}
\end{equation}
holds. 

Let $h_{min}$ be an AZD of $K_{Y}+L_{X/Y}$ with minimal singularities
as in Section \ref{AZD} (cf. Definition \ref{minAZD}).
We set 
\begin{equation}
\Omega_{min} := h_{min}^{-1}\cdot h_{L_{X/Y}}.
\end{equation}
Then we have the following uniform $C^{0}$-estimate for $\{ v_{m}\}$.  

\begin{lemma}\label{c02}
There exists  a positive constant $C$  such that for every $m > s$ 
\begin{equation}
 v_{m} \leqq C + \log \frac{\Omega_{min}}{\Omega_{s}}
\end{equation}
holds on $Y$. \fbox{}
\end{lemma}
{\em Proof of Lemma \ref{c02}}. 
Since 
\begin{equation}
- \mbox{Ric}_{\omega_{m}} + \sqrt{-1}\Theta_{h_{L_{X/Y}}} = \omega_{m}
\end{equation}
holds, we see that  
\begin{equation}
\omega_{m}^{n} =  O(\Omega_{min})
\end{equation}
holds by Definition \ref{minAZD}.
Hence by the uniform upper bound (\ref{uniform}) and Lemma \ref{monotonicity}, we see that there exists a positive constant $C$ such that 
\begin{equation}
v_{m} \leqq C+ \log \frac{\Omega_{min}}{\Omega_{s}}
\end{equation}
holds for every $m > s$. \fbox{} \vspace{3mm} \\ 
Let $dV_{Y}$ be as in (\ref{dvy}) in Lemma \ref{convergencevol}. 
We set 
\begin{equation}
v := \log \frac{n!\cdot dV_{Y}}{\Omega_{s}}. 
\end{equation}
Then by Lemmas \ref{c2} and \ref{c02}, 
\begin{equation}
\omega_{Y} = \omega_{(s)} + \sqrt{-1}\partial\bar{\partial}v
\end{equation}
satisfies 
\begin{equation}
\omega_{Y}^{n} = e^{v}\cdot \Omega_{s}. 
\end{equation}

By Lemmas \ref{c2}, \ref{monotonicity} and \ref{c02}, we see that taking 
a suitable subsequence $\{ m_{k}\}$, if necessary, we may assume that 
\begin{equation}
v = \lim_{k\rightarrow\infty}v_{m_{k}}
\end{equation}
holds in $C^{\infty}$-topology on every compact subset of 
$Y \backslash V$. 
Then by the above construction (cf. (\ref{m-eq}))
\begin{equation}\label{eqK}
-\mbox{Ric}_{\omega_{Y}} + \sqrt{-1}\,\Theta_{h_{L_{X/Y}}} 
= \omega_{Y}
\end{equation}
holds on $Y \backslash V$.

Let us define the singular hermitian metric on $K_{Y}+L_{X/Y}$ by 
\begin{equation}
h_{K} := (dV_{Y})^{-1}\cdot h_{L_{X/Y}} = n!\cdot (e^{v}\cdot \Omega_{s})^{-1}\cdot
h_{L_{X/Y}}.
\end{equation}
We shall check $h_{K}$ is an AZD of $K_{Y}+L_{X/Y}$. 
First it is clear that $\sqrt{-1}\,\Theta_{h_{K}}$ is a closed semipositive
current by (\ref{eqK}) and the $C^{0}$-estimate: Lemma \ref{c02}. 
We note that all the coefficients of $E_{m}$ is less than $1$ by the construction.  Then  by the construction 
every global holomorphic section of $m!\cdot \pi_{m}^{*}(K_{Y}+L_{X/Y})$
is $L^{2}$ integrable on $Y$ with respect to $h^{m!}_{L_{X/Y}}\cdot (\omega_{m}^{n})^{-(m!-1)}$ by the Monge-Amp\`{e}re equation of $u_{m}$ and the 
almost boundedness of $u_{m}$.    
Then by the monotonicity of  $\{\omega_{m,\delta}^{n}\}$ (Lemma \ref{monotonicity}),  we see that $h_{K} = (dV_{Y})^{-1}\cdot h_{L_{X/Y}}$ is an AZD of 
$K_{Y}+L_{X/Y}$.
This completes the proof of Theorem \ref{main} except the uniqueness of 
the canonical semipositive current. 
The uniqueness is the direct consequence of Theorem \ref{DS}. \fbox{}

\subsection{Generalization to general adjoint line bundles}

In the proof of Theorem \ref{main}, we have not used the property of 
the Hodge bundle $L_{X/Y}$ (cf. \ref{lxy}) except the Poincar\'{e} growth property 
of the curvature of the Hodge metric $h_{X/Y}$ (cf. (\ref{hxy})).   
Hence without changing the proof, we have the following variant of 
Theorem \ref{main}.  

\begin{theorem}
Let $Y$ be a smooth projective $n$-fold and let $(L,h_{L})$ be a 
$\mathbb{Q}$-line bundle 
with $C^{\infty}$-metric $h_{L}$ with semipositive curvature.  Suppose that $K_{Y} + L$ is big.
Let $U$ be the Zariski open subset of $U$ defined by 
\begin{equation}
U := \{ y\in Y|\mbox{$|m!(K_{Y}+L)|$ is very ample around $y$ for $m >> 1$}\}.
\end{equation} 
Then there exist a closed positive current $\omega_{Y}$ on $Y$ 
such that 
\begin{enumerate}
\item[(1)]  There exists a sequence 
of closed positive currents $\{\omega_{m}\}$ such that $\omega_{m}|U$ is 
$C^{\infty}$ and satisfies the equation 
\begin{equation}
-\mbox{\em Ric}_{\omega_{m}} + \sqrt{-1}\Theta_{h_{L}} = \omega_{m}
\end{equation}
holds on $U$ and 
\begin{equation}
\omega_{Y} = \lim_{m\rightarrow\infty}\omega_{m}
\end{equation}
in the sense of currents. 
\item[(2)] 
\begin{equation}
h_{can}:= \left(\frac{1}{n!}\,\omega_{Y,abc}^{n}\right)^{-1}\cdot h_{L}
\end{equation}
is an AZD of $K_{Y} + L$. 
\end{enumerate}  
Moreover if the log canonical ring 
\begin{equation}
R(Y,a(K_{Y}+L)) = \oplus_{m=0}^{\infty}\Gamma(Y,\mathcal{O}_{Y}(ma(K_{Y} +L)))
\end{equation}
is finitely generated, then $\omega_{Y}$ is $C^{\infty}$ on 
the Zariski open subset $U$, where $a$ is the minimal positive integer 
such that $aL$ is a genuine line bundle.  \fbox{}
\end{theorem}

\section{Dynamical systems of Bergman kernels}\label{Dy}
Let $f : X \to Y$ be the Iitaka fibration  and let $(L_{X/Y},h_{L_{X/Y}})$ be the singular hermitian 
$\mathbb{Q}$-line bundle on $Y$ as in Theorem \ref{main}. 
Let $\omega_{Y}$ be the canonical  K\"{a}hler current  on $Y$ (cf. Definition \ref{L-K-E}).  
Then  there exists a nonempty Zariski open subset $U$ of $Y$ such that 
$\omega_{Y}$ is a $C^{\infty}$ on $U$ and 
\begin{equation}
-\mbox{Ric}_{\omega_{Y}} + \sqrt{-1}\,\Theta_{h_{L_{X/Y}}} = \omega_{Y}
\end{equation}
constructed in Theorem \ref{main}. 

Let $L_{X/Y} = M_{X/Y} + D_{X/Y}$ be the decomposition as (\ref{MD-dec}).
In this section we shall consider to describe the canonical K\"{a}hler current 
$\omega_{Y}$ in Therem \ref{main}. First we shall consider the case :  
$D_{X/Y}= 0$.  The general case : $D_{X/Y}\neq 0$ (Theorem \ref{DS}) follows from entirely the same line as \cite[Section 4.2]{LC}.   

\subsection{The case: $D_{X/Y} = 0$}
First we shall consider the case that $D_{X/Y} = 0$. 
The reason why we consider this case is that in this case we do not need to 
use the Ricci iterations (cf. Section \ref{IT}).  
 
Let $m_{0}$ be the sufficiently large positive integer, 
$M$ be a effective Cartier divisor such that  
$A := m_{0}!(K_{Y}+L_{X/Y}) - M$ is sufficiently ample as in Section \ref{Ds}. 
Let   $h_{A}$ be the $C^{\infty}$-hermitian metric on 
$A$.   Hereafter we shall consider $h_{A}$ as a singular hermitian metric 
on $m_{0}!(K_{Y}+L_{X/Y})$ by identifying $h_{A}$ with 
\begin{equation}
h_{A}/|\tau_{M}|^{2},
\end{equation}
where $\tau_{M}$ is a global holomorphic section of ${\cal O}_{Y}(M)$ with divisor $M$.   

We shall construct a sequence of singular hermitian metrics
$\{ h_{m}\}_{m\geqq m_{0}!}$ and a sequence of Bergman kernels $\{ K_{m}\}$ 
as follows. \vspace{3mm} \\

We set $h_{m_{0}!} : = h_{A}$ and 
\begin{equation}
K_{m_{0}!+1} := \left\{\begin{array}{ll} K(Y,K_{Y} +m_{0}!(K_{Y}+L_{X/Y}),h_{m_{0}!}), & \mbox{if}\,\, a > 1 \\ 
& \\ 
& \\
K(Y,K_{Y}+L_{X/Y}+ m_{0}!(K_{Y}+L_{X/Y}),h_{L_{X/Y}}\cdot h_{m_{0}!}), & \mbox{if}\,\, a = 1 
\end{array}\right. 
\end{equation}
Let $\{ h_{m}\}_{m\geqq m_{0}!}$ be the corresponding dynamical 
system of singular hermitian metrics defined by 
\begin{equation}
h_{m} := K_{m}^{-1}
\end{equation}
as in Section \ref{Ds}\footnote{Please do not confuse $h_{m}$ in Section 2}.
\begin{theorem}\label{DS2}
 Let $X$ be a smooth projective variety of nonnegative Kodaira dimension 
 and let $f : X \to Y$ be the Iitaka fibration as above.  
 Let  $\{ h_{m}\}_{m \geq m_{0}!}$ be the sequence 
of hermitian metrics as above and let $n$ denote $\dim Y$.
$\omega_{Y}$ is the canonical K\"{a}hler current on $Y$ 
as in Theorem \ref{main}.  
Then 
\begin{equation}\label{hinfty}
h_{\infty} := \liminf_{m\rightarrow\infty} \sqrt[m]{(m!)^{n}\cdot h_{m}}
\end{equation}
is a singular hermitian metric on $K_{Y}+L_{X/Y}$ such that 
\begin{equation}\label{hinftyvol}
h_{\infty} = (2\pi)^{n}\cdot\left(\frac{1}{n!}\omega_{Y,abc}^{n}\right)^{-1}
\cdot h_{L_{X/Y}}
\end{equation}
holds almost everywhere on $Y$ and 
\begin{equation}\label{hinftycurvature}
\omega_{Y}= \sqrt{-1}\,\Theta_{h_{\infty}}
\end{equation}
holds on $Y$.  
In particular $h_{\infty}$ (and hence $\omega_{Y}$) is unique and  is independent of the choice of $A$ and $h_{A}$. 
\fbox{}
\end{theorem}  
Now we shall prove Theorem \ref{DS2}. 
Let $dV_{Y} = \frac{1}{n!}\omega_{Y,abc}^{n}$ be the volume form associated 
with $(Y,\omega_{Y})$. 
This $dV_{Y}$ is the same as the volume form defined as (\ref{dvy}) 
by the proof of Theorem \ref{main}. 
\begin{lemma}\label{lowers}
\begin{equation}
\limsup_{m\rightarrow\infty}\,h_{L_{X/Y}}\cdot\!\!\sqrt[m]{(m!)^{-n}K_{m}}
\geqq (2\pi)^{-n}dV_{Y}
\end{equation}
holds on $X$. \fbox{} 
\end{lemma}
{\em Proof of Lemma \ref{lowers}}. 
First we shall assume that $L_{X/Y}$ is a genuine line bundle on $Y$ for simplicity. 
If $L_{X/Y}$ is not a genuine line bundle, we tensorize  $(aL_{X/Y},h_{L_{X/Y}}^{a})$ in every $a$-steps instead of tensorize $(L_{X/Y},h_{L_{X/Y}})$ every step, where 
$a$ is the least positive integer such that $aL_{X/Y}$ is Cartier. 
But of course this is a minor technical difference.  Hence we shall give a proof assuming that  $L_{X/Y}$ is Cartier.  The general case is left to readers 
to avoid inessential complication. 
\vspace{3mm} \\ 
Let us consider the (singular) hermitian line bundle $(K_{Y}+ L_{X/Y} ,dV_{Y}^{-1}\cdot h_{L_{X/Y}})$ on $Y$.  Let $U$ be a nonempty Zariski open subset of $Y$ such that 
$\omega_{Y}|U$ is a $C^{\infty}$-K\"{a}hler form.  
Let $p\in U$ be a point.  Then by the equation (\ref{LKE}), there exists a holomorphic normal coordinate
$(U,z_{1},\cdots ,z_{n})$ of $(Y,\omega_{Y})$ around $p$  and a local holomorphic frame $\mbox{\bf e}_{L_{X/Y}}$ of $L_{X/Y}$ on $U$  such that 
\begin{equation}\label{(1)}
dV_{Y}^{-1}\cdot h_{L_{X/Y}}=  \{ \prod_{i=1}^{n}(1 -\mid z_{i}\mid^{2}) + O(\parallel z\parallel^{3})\}\cdot 
2^{n}\cdot\mid dz_{1}\wedge\cdots\wedge dz_{n}\mid^{-2}\cdot |\mbox{\bf e}_{L_{X/Y}}|^{-2}
\end{equation} 
and $h_{L_{X/Y}}(\mbox{\bf e}_{L_{X/Y}},\mbox{\bf e}_{L_{X/Y}})(p) = 1$. 
Suppose that
\begin{equation}\label{assumption}
C_{m-1}\cdot h_{A}^{-1}\cdot dV_{Y}^{m-m_{0}!-1}\cdot h_{L_{X/Y}}^{-(m-m_{0}!-1)}\leqq K_{m-1}
\end{equation}
holds on $Y$ for some positive constant $C_{m-1}$. 
We note that
\begin{equation}\label{extremal}
K_{m}(y) = \sup \{ \mid\sigma\mid^{2}(y) ;  
\sigma \in H^{0}(Y,{\cal O}_{Y}(m(K_{Y}+L_{X/Y}))), (\sqrt{-1})^{n^{2}}\!\int_{X}h_{m-1}\cdot\sigma\wedge
\bar{\sigma} = 1 \}
\end{equation} 
holds for every $y\in Y$, by the extremal property of the Bergman kernel
\footnote{This is well known. See for example, \cite[p.46, Proposition 1.3.16]{kr}.}. 
We note that  for the open unit disk $\Delta = \{ t\in \mathbb{C}\mid \,\,\mid t\mid < 1\}$, 
\begin{equation}\label{disk}
\sqrt{-1}\int_{\Delta}(1 - \mid t\mid^{2})^{m}dt\wedge d\bar{t}
= \frac{2\pi}{m+1} 
\end{equation}
holds.  Then  by H\"{o}rmander's $L^{2}$-estimate of $\bar{\partial}$-operators, we see that there exists a positive constant $\lambda_{m}$ such that 
\begin{equation}\label{induction}
(\lambda_{m}\cdot(2\pi)^{-n}\cdot m^{n})\cdot C_{m-1}\cdot  dV_{Y}^{m-m_{0}!}\leqq h_{L_{X/Y}}^{m-m_{0}!}\cdot h_{A}\cdot K_{m}
\end{equation}
with 
\begin{equation}
\lambda_{m} \geqq 1 - \frac{C}{\sqrt{m}},
\end{equation}
where $C$ is a positive constant independent of $m$. 
\vspace{3mm} \\
In fact this can be verified as follows. 
Let $y\in Y \backslash (\mbox{Supp}\, M \cup V)$ and let $(U,z_{1},\cdots ,z_{n})$
be the normal coordinate as above.  We may assume that 
$U$ is biholomorphic to the polydisk $\Delta^{n}(r)$ of radius $r$ with center $O$
in $\mathbb{C}^{n}$ for some $0 < r < 1$ via $(z_{1},\cdots ,z_{n})$. 

Taking $r < 1$ sufficiently small we may assume that 
there exists a  $C^{\infty}$-function  $\rho$ on $Y$ such that 
\begin{enumerate}
\item[(1)] $\rho$ is identically $1$ on $\Delta^{n}(r/3)$.
\item[(2)] $0\leqq \rho \leqq 1$. 
\item[(3)] $\mbox{Supp}\,\rho \subset\subset U$.
\item[(4)] $\mid d\rho\mid < 3/r$, where $\mid\,\,\,\,\,\mid$ denotes 
the pointwise norm with respect to $\omega_{Y}$. 
\end{enumerate}
We note that by  the equation (\ref{(1)}), the mass of 
$\rho\cdot (dz_{1}\wedge\cdots \wedge dz_{n})^{m}\otimes \mbox{\bf e}_{L_{X/Y}}^{m}$ concentrates 
around the origin as $m$ tends to infinity. 
Hence  by (\ref{disk}) we see that the $L^{2}$-norm 
\begin{equation}
\parallel\rho\cdot(dz_{1}\wedge\cdots \wedge dz_{n})^{m}\otimes \mbox{\bf e}_{L_{X/Y}}^{m}\parallel
\end{equation}
 of $\rho\cdot (dz_{1}\wedge\cdots \wedge dz_{n})^{m}\otimes \mbox{\bf e}_{L_{X/Y}}^{m}$
with respect to $(dV_{Y})^{-m}\cdot h_{L_{X/Y}}^{m}$ and $\omega_{Y}$ is asymptotically
\begin{equation}\label{(2)}
\parallel\rho\cdot (dz_{1}\wedge\cdots \wedge dz_{n})^{m}\otimes \mbox{\bf e}_{L_{X/Y}}^{m}\parallel^{2}
\sim 2^{mn}\left(\frac{2\pi}{m}\right)^{n}
\end{equation}
as $m$ tends to infinity , where $\sim$ means that  the ratio of the both sides 
converges to $1$ as $m$ tends to infinity. 
We set 
\begin{equation}
\phi := n\rho\log \sum_{i=1}^{n}\mid z_{i}\mid^{2}. 
\end{equation}
We may and do assume that $m$ is sufficiently large so that
\begin{equation}
(m-m_{0}!-1)\cdot\omega_{Y} +\sqrt{-1}\,\Theta_{h_{A}} +\sqrt{-1}\Theta_{h_{L_{X/Y}}}
+ \sqrt{-1}\partial\bar{\partial}\phi > 0
\end{equation}
holds on $Y$. 
We note that $\bar{\partial}(\rho\cdot (dz_{1}\wedge\cdots \wedge dz_{n})^{m}\otimes \mbox{\bf e}_{L_{X/Y}}^{m})$  vanishes 
on the polydisc of radius $r/3$ with center $p$ as above. 
Then by (\ref{(2)}),the $L^{2}$-norm 
\[
\parallel\bar{\partial}(\rho\cdot (dz_{1}\wedge\cdots \wedge dz_{n})^{m}\otimes \mbox{\bf e}_{L_{X/Y}}^{m})\parallel_{\phi}
\] 
of 
$\bar{\partial}(\rho\cdot (dz_{1}\wedge\cdots \wedge dz_{n})^{m}\otimes \mbox{\bf e}_{L_{X/Y}}^{m})$ 
with respect to \\ $e^{-\phi}\cdot h_{A}\cdot (dV_{Y})^{-(m-m_{0}!-1)}\otimes h_{L_{X/Y}}^{(m-m_{0}!)}$ and $\omega_{Y}$ 
satisfies the inequality 
\begin{equation}\label{(3)}
\parallel\bar{\partial}(\rho\cdot (dz_{1}\wedge\cdots \wedge dz_{n})^{m}\otimes \mbox{\bf e}_{L_{X/Y}}^{m})\parallel_{\phi}^{2}\leqq C_{0}\cdot\left(\frac{3}{r}\right)^{2n+2}\left(1-\frac{r^{2}}{16}\right)^{m}2^{mn}\left(\frac{2\pi}{m}\right)^{n}
\end{equation}
for every $m$, where $C_{0}$ is a positive constant independent of $m$. 
By  H\"{o}rmander's $L^{2}$-estimate applied to the adjoint line bundle of the hermitian line bundle\footnote{More precisely we apply the $L^{2}$-estimate 
on a complete K\"{a}hler manifold $(U,\omega_{Y} + \epsilon \omega_{U})$, 
where $\omega_{U}$ is a complete K\"{a}hler metric on $U$ 
and $\epsilon > 0$ and then let $\epsilon\downarrow 0$.}:
\begin{equation}
((m-1)(K_{Y}+L_{X/Y})+L_{X/Y},e^{-\phi}\cdot h_{A}\cdot dV_{Y}^{-(m-m_{0}!-1)}\cdot h_{L_{X/Y}}^{m-m_{0}!}),
\end{equation}
 we see that for every 
sufficiently large $m$, there exists a $C^{\infty}$- solution $u$ of 
the equation ; 
\begin{equation}
\bar{\partial}u = \bar{\partial}(\rho\cdot (dz_{1}\wedge\cdots \wedge dz_{n})^{m}\otimes\mbox{\bf e}_{L_{X/Y}}^{m})
\end{equation}
such that 
\begin{equation}
u(p) = 0
\end{equation}
and  
\begin{equation}
\parallel u\parallel_{\phi}^{2} \leqq \frac{2}{m}\parallel\bar{\partial}(\rho\cdot (dz_{1}\wedge\cdots \wedge dz_{n})^{m}\otimes \mbox{\bf e}_{L_{X/Y}}^{m})\parallel_{\phi}^{2}
\end{equation}
hold, where $\parallel\,\,\,\,\parallel_{\phi}$'s  denote the $L^{2}$ norms
with respect to $e^{-\phi}\cdot h_{A}\cdot dV_{Y}^{-(m-m_{0}!-1)}\otimes h_{L_{X/Y}}^{m-m_{0}!}$ and $\omega_{Y}$
respectively.   
Then $\rho\cdot (dz_{1}\wedge\cdots \wedge dz_{n})^{m}\otimes \mbox{\bf e}_{L_{X/Y}}^{m}- u$
is a holomorphic section of $m(K_{Y}+L_{X/Y})$ such that 
\begin{equation}
(\rho\cdot (dz_{1}\wedge\cdots \wedge dz_{n})^{m}\otimes \mbox{\bf e}_{L_{X/Y}}^{m}- u)(p)
= \left((dz_{1}\wedge\cdots \wedge dz_{n})^{m}\otimes \mbox{\bf e}_{L_{X/Y}}^{m}\right)(p)
\end{equation}
and
\begin{equation}
\!\!\!\!\!\!\!\!\!\!\!\!\!\!\!\!\!\!\!\!\!\!\!\!\!\!\!\!\!
\parallel\rho\cdot (dz_{1}\wedge\cdots \wedge dz_{n})^{m}\otimes \mbox{\bf e}_{L_{X/Y}}^{m}- u\parallel^{2}
\end{equation}
\vspace{-5mm}
\[
\,\,\,\,\,\,\,\,\,\,\,\,\,\,\,\,\,\,\,\,\,\,\,\,\,\,\,\,\,\,\,\,\,\,\,\,\,
\leqq \left(1+ C_{0}\cdot \left(\frac{3}{r}\right)^{2n+2}\cdot\sqrt{\frac{2}{m}}\cdot\left(1 -\frac{r^{2}}{16}\right)^{m}\right)\cdot 2^{mn}\cdot\left(\frac{2\pi}{m}\right)^{n}. 
\]
Hence by the assumption of the induction (\ref{assumption}) and the extremal propety 
of Bergman kernels, this implies that there exists a positive constant 
$C$ independent of $m$ such that  
\begin{equation}
K_{m}(p) \geqq \left(1 - \frac{C}{\sqrt{m}}\right)\cdot m^{n}\cdot (2\pi)^{-n}\cdot C_{m-1}\cdot\left( h_{A}^{-1}\cdot h_{L_{X/Y}}^{-(m-m_{0}!)}\cdot dV_{Y}^{m-m_{0}!}\right)(p)  
\end{equation} 
holds, since the point norm of $(dz_{1}\wedge\cdots \wedge dz_{n})^{\otimes m}
\otimes {\bf e}_{A}$ at $p$ (with respect to $h_{A}\cdot dV_{Y}^{-(m-m_{0}!)}\cdot h_{L_{X/Y}}^{(m-m_{0}!)}$) is asymptotically equal to $2^{mn}$.
Then by induction on $m$, using (\ref{extremal}) and 
(\ref{induction}), we see that there exist a positive  constant $C^{\prime}$
and a positive intger $m_{1} > m_{0}!$  such that $C/\sqrt{m_{1}}< 1$ and for every $m >  m_{1}$ 
\begin{equation}\label{induc}
K_{m} \geqq C ^{\prime}\left(\prod_{k=m_{1}}^{m}\left(1 -\frac{C}{\sqrt{k}}\right)\right)\cdot (m!)^{n}\cdot (2\pi )^{-mn}\cdot h_{A}^{-1}\cdot h_{L_{X/Y}}^{-(m-m_{0}!)}\cdot dV_{Y}^{m-m_{0}!}
\end{equation}
holds at $p$.  
Moving $p$,  this implies that  
\begin{equation}
\limsup_{m\rightarrow\infty}\,\,h_{L_{X/Y}}\cdot \sqrt[m]{(m!)^{-n}K_{m}}
\geqq (2\pi)^{-n}dV_{Y}
\end{equation}
holds on $Y$.
\fbox{} \vspace{5mm} \\
On the other hand, we obtain the upper estimate by the following lemma. 
\begin{lemma}\label{holder}
\begin{equation}
\int_{Y}h_{L_{X/Y}}\cdot K_{m}^{\frac{1}{m}} \leqq (\prod_{k=m_{0}}^{m}(N_{k}+1))^{\frac{1}{m}}\cdot \left(\int_{Y}h_{L_{X/Y}}\cdot\sqrt[m_{0}]{K_{m_{0}}}\right)^{\frac{m_{0}}{m}}
\end{equation}
holds, where $N_{k} := \dim \mid k(K_{Y}+L_{X/Y})\mid =\dim H^{0}(Y,{\cal O}_{Y}(k(K_{Y}+L_{X/Y}))) -1$. \fbox{}
\end{lemma}
{\em Proof}.
First we note that  the trivial equality:
\begin{equation}
\int_{Y}h_{L_{X/Y}}\cdot\frac{K_{m}}{K_{m-1}} = N_{m} + 1
\end{equation}
holds by the definition of $K_{m}$ and the equality  $h_{m-1} = 1/K_{m-1}$. 
Then by H\"{o}lder's ineqality, we have 
\begin{eqnarray*}
\int_{Y}h_{L_{X/Y}}\cdot K_{m}^{\frac{1}{m}} & = & \int_{Y}h_{L_{X/Y}}\cdot\frac{K_{m}^{\frac{1}{m}}}
{h_{L_{X/Y}}\cdot K_{m-1}^{\frac{1}{m-1}}}\cdot h_{L_{X/Y}}\cdot K_{m-1}^{\frac{1}{m-1}}
\\
& \leqq  & \left(\int_{Y}h_{L_{X/Y}}^{m}\cdot \frac{K_{m}} 
{h_{L_{X/Y}}^{m}\cdot K_{m-1}^{\frac{m}{m-1}}}\cdot  (h_{L_{X/Y}}\cdot K_{m-1}^{\frac{1}{m-1}})\right)^{\frac{1}{m}}
\cdot \left(\int_{Y}h_{L_{X/Y}}\cdot K_{m-1}^{\frac{1}{m-1}}\right)^{\frac{m-1}{m}}
\\
& = & \left(\int_{Y}h_{L_{X/Y}}\cdot \frac{K_{m}} 
{K_{m-1}}\right)^{\frac{1}{m}}
\cdot \left(\int_{Y}h_{L_{X/Y}}\cdot K_{m-1}^{\frac{1}{m-1}}\right)^{\frac{m-1}{m}}
\\
& = & (N_{m}+1)^{\frac{1}{m}}
\cdot \left(\int_{Y}h_{L_{X/Y}}\cdot K_{m-1}^{\frac{1}{m-1}}\right)^{\frac{m-1}{m}}. 
\end{eqnarray*}
Hence we obtain the inequality:
\begin{equation}
\int_{Y}h_{L_{X/Y}}\cdot K_{m}^{\frac{1}{m}} \leqq (N_{m}+1)^{\frac{1}{m}}
\cdot \left(\int_{Y}h_{L_{X/Y}}\cdot K_{m-1}^{\frac{1}{m-1}}\right)^{\frac{m-1}{m}}.
\end{equation}
Continuing this process, by using 
\begin{equation}
\int_{Y}h_{L_{X/Y}}\cdot K_{m-1}^{\frac{1}{m-1}} \leqq (N_{m-1}+1)^{\frac{1}{m-1}}\cdot
\left(\int_{Y}h_{L_{X/Y}}\cdot K_{m-2}^{\frac{1}{m-2}}\right)^{\frac{m-2}{m-1}},
\end{equation}
we have that 
\begin{equation}
\int_{Y}h_{L_{X/Y}}\cdot (K_{m})^{\frac{1}{m}} \leqq \{(N_{m}+1)\cdot (N_{m-1}+1)\}^{\frac{1}{m}}
\cdot \left(\int_{Y}h_{L_{X/Y}}\cdot (K_{m-2})^{\frac{1}{m-2}}\right)^{\frac{m-2}{m}}
\end{equation}
holds.
Continueing this process we obtain the lemma. \fbox{} \vspace{3mm} \\ 

To estimate the growth of $\{ N_{m}\}_{m\geq m_{0}!}$, we  
introduce the following notion. 

\begin{definition}\label{volume} Let $L$ be a line bundle on a compact complex 
manifold $M$ of dimension $n$. 
We define the {\bf volume} $\mu (M,L)$ of $M$ with respect to $L$ by
\begin{equation}
\mu (M,L) := n!\cdot\limsup_{m\rightarrow\infty}m^{-n}
\dim H^{0}(M,{\cal O}_{M}(mL)).
\end{equation}
\fbox{} \end{definition}

\noindent We note that Definition \ref{volume} can be generalized to 
the case of $\mathbb{Q}$-line bundles in an obvious way. 
 Using Lemma \ref{holder}, we obtain the following lemma. 
\begin{lemma}\label{uppers}
\begin{equation}
\limsup_{m\rightarrow\infty}\frac{1}{(m!)^{\frac{n}{m}}}\int_{Y}h_{L_{X/Y}}\cdot (K_{m})^{\frac{1}{m}} \leqq \frac{\mu (Y,K_{Y}+L_{X/Y})}{n!}
\end{equation}
holds. \fbox{}
\end{lemma}
{\em Proof}.
By the definition of the volume $\mu (Y,K_{Y}+L_{X/Y})$,  
\begin{equation}
N_{m}+1 = \frac{\mu (Y,K_{Y}+L_{X/Y})}{n!}m^{n} + o(m^{n})
\end{equation}
holds. Then by Lemma \ref{holder}, we see that 
\begin{equation}
\limsup_{m\rightarrow\infty}\frac{1}{(m!)^{\frac{n}{m}}}\int_{Y}h_{L_{X/Y}}\cdot (K_{m})^{\frac{1}{m}} \leqq \frac{\mu (Y,K_{Y}+L_{X/Y})}{n!}
\end{equation}
holds. \fbox{}

\begin{lemma}\label{volume growth}
\begin{equation}
\frac{1}{(2\pi)^{n}}\int_{Y}dV_{Y} = \frac{1}{n!}\int_{Y}\left(\frac{1}{2\pi}\omega_{Y,abc}\right)^{n}
= \frac{1}{n!}\,\mu (Y,K_{Y}+L_{X/Y})
\end{equation}\,
holds. \fbox{}
\end{lemma}
{\em Proof of Lemma \ref{volume growth}}.
Let  $|P_{m}|$ be the free part of $|\pi_{m}^{*}m!(K_{Y}+L_{X/Y})|$, 
where $\pi_{m}$ is the resolution of $\mbox{Bs}|m!(K_{Y}+L_{X/Y})|$ 
as in the last section (cf. (\ref{pim})). 
By Fujita's theorem (\cite[p.1,Theorem]{f}), we see that 
\begin{equation}\label{feq}
\lim_{m\rightarrow\infty}(m!)^{-n}P_{m}^{n} = \mu (Y,K_{Y}+L_{X/Y})
\end{equation}  
Then by (\ref{uniform}),(\ref{bouck}), Lemma \ref{monotonicity} and (\ref{feq}), Lebesgue's bounded convergence theorem implies that 
\begin{equation}
\mu (Y,K_{Y}+L_{X/Y})= \lim_{m\rightarrow\infty}\int_{Y}\left(\frac{1}{2\pi}\omega_{m,abc}\right)^{n} = \int_{Y}\left(\frac{1}{2\pi}\omega_{Y,abc}\right)^{n}
\end{equation}
hold.  This implies the lemma. 
\fbox{}  \vspace{3mm} \\

\noindent  We note that by Lemma \ref{uppers} and the submeanvalue inequality for plurisubharmonic functions,  $\{ h_{L_{X/Y}}\cdot \sqrt[m]{(m!)^{-n}K_{m}}\}$ 
is a family of uniformly bounded semipositive $(n,n)$ forms on $Y$.
  Then by Lebesgue's bounded convergence theorem, we see that 
\begin{equation}
\limsup_{m\rightarrow\infty}\int_{Y}h_{L_{X/Y}}\cdot \sqrt[m]{(m!)^{-n}K_{m}} 
= 
\int_{Y}\limsup_{m\rightarrow\infty}\,\,h_{L_{X/Y}}\cdot \sqrt[m]{(m!)^{-n}K_{m}}
\end{equation}
holds. 
Combining Lemmas \ref{lowers}, \ref{uppers} and  \ref{volume growth}, we have the equality:
\begin{equation}
\limsup_{m\rightarrow\infty}\frac{1}{(m!)^{\frac{n}{m}}}\,h_{L_{X/Y}}\cdot K_{m}^{\frac{1}{m}}
= (2\pi)^{-n}dV_{Y},
\end{equation}
holds almost everywhere on $Y$.
Hence by the definition of $h_{\infty}$ (cf. (\ref{hinfty}))
\begin{equation}
h_{\infty} = \left(\limsup_{m\rightarrow\infty}\frac{1}{(m!)^{\frac{n}{m}}}\, K_{m}^{\frac{1}{m}}\right)^{-1} = (2\pi)^{n}\!\!\cdot dV_{Y}^{-1}\cdot h_{L_{X/Y}}
\end{equation}
hold almost everywhere on $Y$. 
This implies the equality (\ref{hinftyvol}) in Theorem \ref{DS}.  Then by the equation (\ref{LKE})
we have the equality (\ref{hinftycurvature}) in Theorem \ref{DS}:
\[
 \omega_{Y} = \sqrt{-1}\Theta_{h_{\infty}}. 
\]
This completes the proof of Theorem \ref{DS} assuming that $L_{X/Y}$ is Cartier. 
The proof of the general case can be obtained by entirely the same estimates. 
More precisely if $L_{X/Y}$ is not a genuine line bundle,  we may have small ripple in  
the $L^{2}$-estimates in the proof of Lemma \ref{lowers},
since we tensorize $(L_{X/Y}^{\otimes a},h_{L_{X/Y}}^{\otimes a})$ every 
$a$ steps.  But the ripple disappears when we take the normalized limit
as is easily be seen.  
This completes the proof of Theorem \ref{DS2}. \fbox{}

\subsection{The case: $D_{X/Y}\neq 0$}

In this section we shall give the proof of Theorem \ref{DS}. The proof given below is essentially the same as the one of \cite[Theorem 4.8]{LC}.  

Let $K_{Y} + L_{X/Y} = P + N$ be the Zariski decomposition as (\ref{ZD2}). 
Let $a$ be the positive integer such that $aP,aL_{X/Y},aN \in \mbox{Div}(Y)$ as 
in Theorem \ref{DS}.  
Then  we have the following lemma. 
\begin{lemma}\label{MI}
Let $h$ be an arbitrary AZD of  $K_{Y} + L_{X/Y}$.
Then for every positive integer $m$,   
\[
\mathcal{I}(h^{ma}) = \mathcal{I}(h^{ma-1}\cdot h_{L_{X/Y}})
\]
holds. \fbox{}
\end{lemma}
{\em Proof}.  
We note that the followings holds:
\begin{enumerate}
\item[(1)] $(Y,D_{X/Y})$ is KLT. 
\item[(2)] $h_{L_{X/Y}}$ induces 
a metric on $M_{X/Y}$ with logarithmic growth (Lemma \ref{logg}).  
\item[(3)] $aP,aN,aL_{X/Y}\in \mbox{Div}(Y)$. 
\item[(4)] $P$ is semiample (\cite{b-c-h-m}). 
\end{enumerate}
Hence we see that 
\[
\mathcal{O}_{Y}(ma(K_{Y}+L_{X/Y}))\otimes \mathcal{I}(h^{ma}) 
\simeq \mathcal{O}_{Y}(maP)
\]
and 
\[
\mathcal{O}_{Y}(ma(K_{Y}+L_{X/Y}))\otimes \mathcal{I}(h^{ma-1}\cdot h_{L_{X/Y}})\simeq
\mathcal{O}_{Y}(maP)
\]
hold.  This implies the lemma, \fbox{} \vspace{3mm} \\

\noindent Let $\{ h_{\ell,m}\}$ be the sequence of singular hermitian metrics 
constructed as Theorem \ref{DS}. Then by Lemma \ref{MI}, we have the following lemma.   
\begin{lemma} 
If we take $A$ to be sufficiently ample, then we have the followings. 
\begin{enumerate}
\item[(1)] $\mathcal{O}_{Y}(A+ (\ell+1)a(K_{Y}+L_{X/Y}))\otimes\mathcal{I}(h_{\ell,m}\cdot h_{L_{X/Y}})$ is globally generated for every $\ell,m\geqq 1$.  
\item[(2)] Let $h_{P}$ is a $C^{\infty}$-hermitian metric 
 on $P$ and $h_{N} = |\sigma_{N}|^{-2}$ where $\sigma_{N}$ is a multivalued 
 holomorphic section on $N$ with divisor $N$. Then 
\[
h_{\ell,m} = O(h_{A}\cdot (h_{P}\cdot h_{N})^{\ell a})
\]
holds for every $\ell,m \geqq 1$. 
\end{enumerate}
\fbox{}
\end{lemma} 
Then the rest of the proof follows from the parallel argument as in Section 3.1. And it is entirely the same as the one of \cite[Theorem 4.8,Section 4.2]{LC}.
\fbox{} 

\subsection{Uniqueness of canonical measures}
For the uniqueness of the canonical measure, we have the following uniqueness. 

\begin{corollary}\label{bir}
$d\mu_{can}$ is birationally invariant. \fbox{} 
\end{corollary}
{\em Proof of Corollary \ref{bir}}. 
Let $f : X \to Y$ be as above and let us consider 
the following commutative diagram :
\begin{equation*}
\begin{picture}(400,400)
\setsqparms[1`1`1`1;350`350]
\putsquare(0,00)%
[\tilde{X}` X  ` \tilde{Y}`Y;
\pi`\tilde{f}`f`\varpi]
\end{picture}
\end{equation*}  
where $\pi : \tilde{X}\to X, \varpi : \tilde{Y}\to 
Y$ are modifications.  Let $\tilde{A}$ and $A$ be ample line bundles
on $\tilde{Y}$ and $Y$ respectively.

Then $\varpi^{*}A$ is nef and big on $\tilde{Y}$.  
Hence by Kodaira's lemma, we have that there exists a positive integers $a_{1},a_{2}$ such that 
\begin{equation}
a_{1}\tilde{A} - \varpi^{*}A, a_{2}\varpi^{*}A - \tilde{A} 
\end{equation}  
are $\mathbb{Q}$-effective. 
Let $d\mu_{\tilde{X},can}, d\mu_{X,can}$ be canonical measures on 
$\tilde{X}$ and $X$ respectively. 
Then since $a_{1}\tilde{A} - \varpi^{*}A$ is $\mathbb{Q}$-effective, by Theorem \ref{DS} and its proof, we see that 
\begin{equation}
d\mu_{\tilde{X},can} \geqq \pi^{*}d\mu_{X,can}
\end{equation}
holds.    
In fact this can be verified as follows.   
Let us fix  $C^{\infty}$-hermitian metrics $h_{A}$ and $h_{\tilde{A}}$
on $A$ and $\tilde{A}$ respectively.  
Let  $b$ be a sufficiently large positive integer and let 
\begin{equation}
\tau \in H^{0}(\tilde{Y},{\cal O}_{\tilde{Y}}(b(a_{1}\tilde{A}-\varpi^{*}A)))
\end{equation}
be a nonzero section such that 
\begin{equation}
(h_{\tilde{A}}^{a_{1}}\cdot \varpi^{*}h_{A}^{-1})^{b}(\tau,\tau) \leqq 1. 
\end{equation}
Let $\{\tilde{K}_{m}\}_{m\geqq 0}, \{ K_{m}\}_{m\geqq 0}$ 
be the dynamical systems of Bergman kernels as in Theorem \ref{DS}, starting from $(ba_{1}\tilde{A},h_{\tilde{A}}^{ba_{1}})$ and $(bA,h_{A}^{b})$ on $\tilde{Y}$ and $Y$ respectively. 
Then by the extremal property of Bergman kernels, we see that 
\begin{equation}
\tilde{K}_{\ell,1} \geqq |\tau|^{2}\cdot \varpi^{*}K_{\ell,1}
\end{equation}
holds for every $\ell \geqq 1$. 
Hence by Theorem \ref{DS}, we see that 
\begin{equation}
d\mu_{\tilde{X},can} \geqq \pi^{*}d\mu_{X,can}
\end{equation}
holds.  

Similarly since $a_{2}\varpi^{*}A - \tilde{A}$ is $\mathbb{Q}$
-effective, we have the opposite inequality :
\begin{equation}
d\mu_{\tilde{X},can} \leqq \pi^{*}d\mu_{X,can}
\end{equation}
Hence we have that the equality 
\begin{equation}
d\mu_{\tilde{X},can} = \pi^{*}d\mu_{X,can}
\end{equation}
holds.  This completes the poof of Corollary \ref{bir}.  \fbox{}

\section{Relative version of Theorems \ref{main} and \ref{DS}}

In this section we shall consider  variation of the canonical measures on projective families.   Our result is as follows.  

\begin{theorem}\label{relative}
Let $f : X \to S$ be a projective family such that $X,S$ are 
smooth and $f$ has connected fibers. 
Suppose that $f_{*}{\cal O}_{S}(mK_{X/S}) \neq 0$ for some $m > 0$.
There exists a relative measure $d\mu_{can,X/S}$ such that  
the singular hermitian metric $h_{X/S} := d\mu_{can,X/S}^{-1}$ on  $K_{X/S}$ 
satisfies: 
\begin{enumerate}
\item[(1)] $\omega_{X/S}:= \sqrt{-1}\,\Theta_{h_{X/S}}$ is semipositive on $X$.
\item[(2)] For every smooth fiber $X_{s} := f^{-1}(s)$, 
$h_{X/S}|X_{s}$ is well defined and is an AZD of $K_{X_{s}}$.  
\item[(3)] There exists a set $T$ of measure $0$ on $S$ such that 
for every $s\in S \backslash T$, $X_{s}$ is smooth and  $\omega_{X/S}|X_{s}$ is the canonical semipositive current on $X_{s}$ constructed as in Theorems \ref{main} and \ref{DS}. \fbox{} 
\end{enumerate} 
\end{theorem}
\begin{remark}
Even for $s \in S \backslash T$, $d\mu_{can,X/S}|X_{s}$ may not be precisely 
equal to the canonical measure $d\mu_{can,s}$ on $X_{s}$ 
as a degenerate volume form on $X_{s}$.  But as a measure $d\mu_{can,X/S}|X_{s}
= d\mu_{can,s}$ holds in exact sense. \fbox{}   
\end{remark}
We call $d\mu_{can,X/S}$ {\bf the relative canonical measure} of 
$f : X \to S$. 
Combining the logarithmic plurisubharmonicity of Bergman kernels (\cite{b3,b-p} and \cite[Theorem 3.4]{KE}), this theorem  strengthens the following famous result due to Y. Kawamata. 
\begin{theorem}(\cite[p.57,Theorem 1]{ka1})\label{kawamata}
Let $f : X \to S$ be an algebraic fiber space.   
Suppose that $\dim S = 1$.   Then  for every positive integer $m$,  
$f_{*}\mathcal{O}_{X}(mK_{X/S})$ is a semipositive vector bundle on $Y$, in the  sense that every quotient $\mathcal{Q}$ of $f_{*}\mathcal{O}_{X}(mK_{X/S})$, 
 $\deg \mathcal{Q} \geqq 0$ holds. 
\fbox{}
\end{theorem} 
The main difference between Theorems \ref{relative} and \ref{kawamata} 
is that the semipositivity is on the total space in Theorem \ref{relative},
while the semipositity is on the direct image of the relative pluricanonical 
systems in Theorem \ref{kawamata}.  In \cite{LC}, we consider the relative log canonical bundle of a family of log canonical pairs.  In the case of log canonical pairs, this difference becomes an essential one.  \vspace{3mm} \\ 

\noindent{\em Proof of Theorems \ref{relative}}. 
Since the assertion is local, we may assume that $S$ is the unit open 
polydisk in $\mathbb{C}^{n}$. 
Let $m_{0}$ be a sufficiently large positive integer and let 
\begin{equation}
F_{m_{0}} := f_{*}{\cal O}_{X}(m_{0}!K_{X/S})   
\end{equation}
and (shrinking $S$, if necessary) let $\sigma_{0},\cdots ,\sigma_{N(m)}$ be a set of  global generators of 
$F_{m_{0}}$ on $S$. 
We take the image $Y$ of the rational map
\begin{equation}
\Phi_{m_{0}} : X -\cdots\rightarrow \mathbb{P}^{N(m_{0})}_{S}. 
\end{equation}
If we take $m_{0}$ sufficienly large,  taking modifications of $\hat{X}$ of $X$ and $\hat{Y}$ of $Y$ respectively, 
we have the relative Iitaka fibration 
\begin{equation}
\begin{picture}(400,400)
\setsqparms[1`1`1`1;350`350]
\putVtriangle(0,00)%
[\hat{X}`\hat{Y} ` S;g `\hat{f} `h]
\end{picture}
\end{equation}  
such that $\hat{X}$ and $\hat{Y}$ are smooth and 
$g_{*}{\cal O}_{\hat{X}}(m!K_{\hat{X}/\hat{Y}})^{**}$ is a line bundle 
on $\hat{Y}$.  We define the $\mathbb{Q}$-line bundle $L_{X/Y}$ on $\hat{Y}$
by 
\begin{equation}
L_{X/Y} = \frac{1}{m_{0}!}g_{*}{\cal O}_{\hat{X}}(m_{0}!K_{\hat{X}/\hat{Y}})^{**} 
\end{equation}
and let $a$ be the least positive integer such that 
$\hat{f}_{*}{\cal O}_{\hat{X}}(aK_{\hat{X}/\hat{Y}})\neq 0$. 
Hereafter we shall replace $f : X \to S$ by 
$\hat{f} : \hat{X} \to S$ and replace $X$ and $Y$ by 
$\hat{X}$ and $\hat{Y}$ respectively.  This does not affect the proof 
of Theorem \ref{relative} by the birational invariance of 
cananonical measures. 

Let $S^{\circ}$ be the locus of $S$ such that $f$ is smooth over $S^{\circ}$. 
Let $A$ be a sufficienly ample line bundle on $Y$ and let $h_{A}$ be a $C^{\infty}$-metric on $A$ with strictly positive curvature. 
Then as in Section \ref{Dy}, for every $s\in S^{\circ}$, we define the 
dynamical system of Bergman kernels $\{K_{\ell,m,s}\}_{\ell,m\geq 1}$ as in 
Section \ref{Dy}. 
We note that the Hodge metric $h_{L_{X/Y}}$ on $L_{X/Y}$ defined as in Section \ref{Dy}
has semipositive curvature in the sense of current on $Y$ (not on every fiber $Y_{s}$) by (\cite{ka1},\cite[p.174,Theorem 1.1]{ka3}). 
Then by the plurisubharmonicity  of the Bergman kernel(\cite{b3,b-p},\cite[Theorem 3.4]{KE}) of the adjoint line bundle of singular hermitian line bundle of  semipositive curvature current,  by induction on $m$, we see that 
\begin{equation}
h_{\ell,m} := (K^{*}_{\ell,m})^{-1}
\end{equation}
extends to a singular hermitian metric on 
\begin{equation}
A+ a\ell (K_{Y/S} + L_{X/Y})
\end{equation}
on $Y$ and the extended metric has semipositive curvature in the sense of current, i.e. $\log K^{*}_{\ell,m}$ 
is plurisubharmonic on $Y$ by Theorem \ref{Lelong}.  Then by Theorem \ref{DS}, 
\begin{equation}
K_{m}:= \mbox{the upper semicontinuous envelope of}\,\,\,\, \limsup_{m\rightarrow\infty} \sqrt[am]{(m!)^{-n}K_{m}^{*}}
\end{equation}
exists as a nontrivial $L_{X/Y}$-valued relative volume form on $Y$ and 
\begin{equation}
h_{m} := K_{m}^{-1}
\end{equation}
is a singular hermitian metric on $K_{Y/S} + L_{X/Y}$ with 
semipositive curvature current. 
We set 
\begin{equation}
h_{m} : = g^{*}h_{m}.  
\end{equation}
Then as before we may consider $h_{X/S}$ as a singular hermitian metric on $K_{X/S}$ 
with semipositive curvature current, i.e.,
\begin{equation}
\omega_{X/S} :=  \sqrt{-1}\,\Theta_{h_{X/S}}
\end{equation} 
is semipositive on $X$. 
By Theorem \ref{DS} and the birational invariance of the canonical semipositive current (Corollary \ref{bir}), there exists a subset $T$ of measure $0$ 
on $S$ such that $S \backslash T$ is contained in $S^{\circ}$ and  for every $s\in S \backslash T$, $\omega_{X/S}|X_{s}$ is the canonical semipositive current 
on $X_{s}$.   Moreover for $s\in T \cap S^{\circ}$, we see that 
$h_{X/S}|X_{s}$ is an AZD of $K_{X_{s}}$ by the very definition of the 
upper-semi-continuous envelope.  
This completes the proof of Theorem \ref{relative}. \fbox{}

Author's address\\
Hajime Tsuji\\
Department of Mathematics\\
Sophia University\\
7-1 Kioicho, Chiyoda-ku 102-8554\\
Japan \\
e-mail address: tsuji@mm.sophia.ac.jp  or h-tsuji@h03.itscom.net


\begin{thebibliography}{99}
\bibitem[A]{a} Aubin, T.: Equation du type Monge-Amp\`{e}re sur les variet\'{e}
 k\"{a}hlerienne compactes, C.R. Acad. Paris {\bf 283} (1976), 459-464.
\bibitem[B1]{b1} Berndtsson, B.: Subharmonicity properties of the Bergman kernel and some other functions associated to pseudoconvex domains, math.CV/0505469 (2005). 
\bibitem[B2]{b2} Berndtsson, B.: Curvature of vector bundles and subharmonicity
of vector bundles, math.CV/050570 (2005).
\bibitem[B3]{b3} Berndtsson, B.:  Curvature of vector bundles associated to holomorphic fibrations, math.CV/0511225 (2005).
\bibitem[B-P]{b-p} Berndtsson, B. and Paun, M. : 
Bergman kernels and the pseudoeffectivity of relative canonical bundles, 
math.AG/0703344 (2007).
\bibitem[B-C-H-M]{b-c-h-m} Birkar, C.-Cascini, P.-Hacon,C.-McKernan, J.:
  Existence of minimal models for varieties of log general type, arXiv:math/0610203.
\bibitem[Bo]{bouk} Boucksom, S.: On the volume of a line bundle, Internat. J. Math. {\bf 13} (2002), 1043-1063.     
\bibitem[D-P-S]{d-p-s}  Demailly, J.P.- Peternell, T.-Schneider, M. : 
Pseudo-effective line bundles on compact K\"{a}hler manifolds, 
International Jour. of Math. {\bf 12} (2001), 689-742. 
\bibitem[F-M]{f-m} Fujino, O. and Mori, S.: Canonical bundle formula, J. Diff. Geom. {\bf 56} (2000), 167-188.
\bibitem[F]{f} Fujita, T.: Approximating Zariski deecomposition
of big line bundle , Kodai Math. J. {\bf 17} (1994), 1-4.
\bibitem[G]{griff} Griffiths, Ph.: Periods of integrals on algebraic
manifolds III: Some global differential-geometric properties of the
period mapping,  Publ.\ Math., Inst.\ Hautes Etud. Sci. {\bf 38} (1970)
, 125--180 .
\bibitem[Ka1]{kawa} Kawamata, Y.: Characterization of Abelian Varieties,
Compos.\ Math. {\bf 43} 253--276 (1981). 
\bibitem[Ka2]{ka1} Kawamata, Y.: Kodaira dimension of Algebraic fiber spaces over curves, Invent. Math. {\bf 66} (1982), 57-71.
\bibitem[Ka3] {ka3}Kawamata, Y.: On effective non-vanishing and base-point-freeness. Kodaira's issue. Asian J. Math. {\bf 4}(2000), no. 1, 173--181.
\bibitem[Kr]{kr} Krantz, S.: Function theory of several complex variables, 
John Wiley and Sons (1982).
\bibitem[L]{l} Lelong, P.: Fonctions Plurisousharmoniques et Formes 
Differentielles Positives, Gordon and Breach (1968).  
\bibitem[N]{n}Nadel, A.M.: Multiplier ideal sheaves and existence of K\"{a}hler-Einstein
metrics of positive scalar curvature, Ann. of Math. {\bf 132}(1990),549-596.
\bibitem[R]{r} Royden, H.L.: The Ahlfors Schwarz lemma in several complex 
variables, Comment. Math. Helv. {\bf 55} (1980),547-558.     
\bibitem[Sch]{sch} Schmid, W.: Variation of Hodge structure: the
singularities of the period mapping. Invent.\ math. {\bf 22}, 211--319
(1973). 
\bibitem[S-T]{s-t} Song, J. and Tian, G. : Canonical measures and K\"{a}hler-Ricci flow, arXiv:0802.2570 (2008). 
\bibitem[Su]{s} Sugiyama, K.: Einstein-K\"{a}hler metrics on minimal varieties of general type and an inequality between Chern numbers. 
Recent topics in differential and analytic geometry, 417--433, 
Adv. Stud. Pure Math., {\bf 18}-{\rm I}, 
Academic Press, Boston, MA (1990).
\bibitem[Tr]{tr} Trudinger, N.S.: Fully nonlinear elliptic equation under 
natural structure conditions, Trans. A.M.S. {\bf 272} (1983), 751-769. 
\bibitem[T0]{tu0} Tsuji H.: Existence and degeneration of K\"{a}hler-Einstein metrics on minimal algebraic varieties of general type. Math. Ann. {\bf 281} (1988), no. 1, 123--133. 
\bibitem[T1]{tu}Tsuji H.: Analytic Zariski decomposition, Proc. of Japan Acad.
{\bf 61}(1992), 161-163.
\bibitem[T2]{tu2} Tsuji, H.:  Existence and Applications of Analytic Zariski Decompositions, Trends in Math., Analysis and Geometry in Several Complex Variables(Katata 1997), Birkh\"{a}user Boston, Boston MA.(1999), 253-272.
\bibitem[T3]{tu3}Tsuji, H.: Deformation invariance of plurigenera, Nagoya Math. J. {\bf 166} (2002), 117-134.
\bibitem[T4]{KE} Tsuji, H.: Dynamical construction of K\"{a}hler-Einstein metrics, math.AG/0606023 (2006), Nagoya Maht. J. {\bf 199}, 107-122 (2010).
\bibitem[T5]{canAZD} Tsuji, H.: Canonical singular hermitian metrics on 
relative canonical bundles, math.ArXiv0704.0566 (2007).
\bibitem[T6]{LC} Tsuji, H.: Ricci iterations and canonical K\"{a}hler currents on LC pairs, math.arXiv.0903.5445(2009).  
\bibitem[T7]{global} Tsuji, H.: Global generation of the direct images of pluri log canonical systems, arXiv.math.1012.0884 (2010).  
\bibitem[Y1]{y} Yau, S.-T.: On the Ricci curvature of a compact K\"{a}hler manifold and the complex Monge-Amp\`{e}re equation,  Comm. Pure  Appl. Math. {\bf 31} (1978),339-441.
\bibitem[Y2]{y2} Yau, S.-T.: A general Schwarz lemma for K\"{a}hler manifolds, Amer. J. of Math. {\bf 100} (1978), 197-303. 
\end{thebibliography}
\end{document}